\documentclass{article}
\usepackage{liyuancao}
\input{preamble.tex}

\title{The Error in Multivariate Linear Extrapolation with Applications to Derivative-Free Optimization}

\author{Liyuan Cao \thanks{School of Mathematics, Nanjing University, Nanjing, China
(caoliyuan@nju.edu.cn).}
        \and Zaiwen Wen \thanks{Beijing International Center for Mathematical Research, Peking University, Beijing, China
(wenzw@pku.edu.cn).}
        \and Ya-xiang Yuan \thanks{State Key Laboratory of Scientific and Engineering Computing, Academy of Mathematics and Systems Science, Chinese Academy of Sciences, Beijing, China (yyx@lsec.cc.ac.cn) }
}

\begin{document}

\maketitle

\begin{abstract}
  We study in this paper the function approximation error of multivariate linear extrapolation. 
The sharp error bound of linear interpolation already exists in the literature. 
However, linear extrapolation is used far more often in applications such as derivative-free optimization, while its error is not well-studied. 
We introduce in this paper a method to numerically compute the sharp bound on the error, and then present several analytical bounds along with the conditions under which they are sharp. 
We also provide a complexity analysis of a basic simplicial search method to illustrate an application of these error bounds in derivative-free optimization. 
All results are under the assumptions that the function being interpolated has Lipschitz continuous gradient and is interpolated on an affinely independent sample set. 
\end{abstract}

\section{Introduction} \label{sec:intro}
Polynomial interpolation is one of the most basic techniques for approximating functions and plays an essential role in applications such as finite element methods and derivative-free optimization. 
This led to a large amount of literature concerning its approximation error. 
This paper contributes to this area of study by analyzing the function approximation error of linear interpolation and extrapolation. 
Specifically, given a function $f: \R^n \rightarrow \R$ and an affinely independent sample set $\Theta:= \{\xx_1,\xx_2,\dots,\xx_{n+1}\} \subset \R^n$, one can find a unique affine function $\hat{f}: \R^n \rightarrow \R$ such that $\hat{f}(\xx_i) = f(\xx_i)$ for all $i \in \{1,\dots,n+1\}$. 
We investigate in this paper the (sharp) upper bound on the approximation error $|\hat{f}(\xx_0) - f(\xx_0)|$ at any given $\xx_0\in\R^n$ under the assumption that $f \in C_\nu^{1,1}(\R^n)$, where the class $C_\nu^{1,1}(\R^n)$ represents the continuously differentiable functions defined on $\R^n$ with their first derivative $Df$ being $\nu$-Lipschitz continuous, i.e., 
\begin{equation} \label{eq:Lipschitz}
    \|Df(\u) - Df(\v)\| \le \nu \|\u - \v\| \quad \text{for all } \u,\v \in \R^n, 
\end{equation}
where $\nu>0$ is the Lipschitz constant, and the norms are Euclidean. 

The main motivation behind this investigation is model-based derivative-free optimization (DFO) methods, which iteratively employ linear interpolation to approximate the black-box objective function and combine the resulting approximation model with trust region methods for optimization (see, e.g., \cite{powell1994direct} and \cite{DFO_book}). 
A large gap currently exists between the theory and the practice of these algorithms. 
The practical algorithms, particularly the ones developed by Powell, gain their efficiency by reusing points evaluated in previous iterations when constructing new interpolation models. 
As a poorly positioned interpolation point may adversely affect the model's accuracy, these algorithms usually contain an intricate mechanism to choose which old points to keep in the interpolation set and which new ones to add. 
Meanwhile, the existing theoretical analyses (e.g. \cite{gratton2018complexity,chen2018stochastic,blanchet2019convergence,cao2023first}) rely on varying assumptions that the approximation models are sufficiently accurate (with sufficiently high probability). 
Under these assumptions, the complexities of the trust region methods are established, but the methods to obtain such models are either simple (and inefficient in practice) or left unexplained. 
This gap prevents the existing theoretical results from accurately representing the efficiency of the practical algorithms. 

To bridge this gap, at least one essential building block is missing: a theory that relates the sample set to the interpolation model's accuracy. 
We aim to build such a block in this paper by analyzing the approximation error of linear interpolation and extrapolation. 
It is worth to point out here that the sharp upper bound on $|\hat{f}(\xx) - f(\xx)|$ is already discovered and proved in \cite{waldron1998error} for linear interpolation, but only for the case when the word ``interpolation'' is used in its narrow sense, i.e., when $\xx \in \conv(\Theta)$, the convex hull of $\Theta$. 
However, the approximation model $\hat{f}$ is used more often than not in model-based DFO methods to estimate the function value at points outside $\conv(\Theta)$. 
As illustrated in Figure~\ref{fig:DFO-TR}, these optimization algorithms attempt to minimize the objective function by alternately constructing a linear interpolation model and minimizing the model inside a trust region. 
The minimizer of the model inside the trust region would then have its function value evaluated and become a part of the sample set for constructing the next iteration's model.  
In practice, this minimizer is rarely located inside $\conv(\Theta)$, necessitating further analysis of the extrapolation case.

\begin{figure}[tbhp]
    \centering
    \subfloat[\raggedright Two consecutive iterations of a DFO algorithm based on linear interpolation and trust region method. The circle represents the trust region, which changes center and expands after finding the minimizer $\xx_4$ that has a lower function value than the current center $\xx_3$. ]{\label{fig:DFO-TR}\resizebox{0.7\linewidth}{!}{\begin{tikzpicture}
    \filldraw[black] (-0.3,1) circle (2pt) node[anchor=east] {$\xx_1$}; 
    \filldraw[black] (-1.1,-0.5) circle (2pt) node[anchor=south east] {$\xx_2$}; 
    \filldraw[black] (1,0) circle (2pt) node[anchor=west] {$\xx_3$}; 
    \filldraw[black] (1+0.4,0-0.9165) circle (2pt) node[anchor=north west] {$\xx_4$}; 
    
    \draw[thick] (1,0) -- (-0.3,1);
    \draw[thick] (-0.3,1) -- (-1.1,-0.5);
    \draw[thick] (-1.1,-0.5) -- (1,0);
    
    \draw[fill=none](1,0) circle (1.0);

    \node at(2.8,0) [below] {\Large$\Rightarrow$};

    \filldraw[black] (-0.3+5,1) circle (2pt) node[anchor=east] {$\xx_1$}; 
    \filldraw[black] (-1.1+5,-0.5) circle (2pt) node[anchor=south east] {$\xx_2$}; 
    \filldraw[black] (1+5,0) circle (2pt) node[anchor=west] {$\xx_3$}; 
    \filldraw[black] (1+0.4+5,0-0.9165) circle (2pt) node[anchor=north west] {$\xx_4$}; 
    
    \draw[thick] (1+5,0) -- (1+0.4+5,0-0.9165);
    \draw[thick] (1+0.4+5,0-0.9165) -- (-1.1+5,-0.5);
    \draw[thick] (-1.1+5,-0.5) -- (1+5,0);
    
    \draw[fill=none](1+0.4+5,0-0.9165) circle (1.15); 
\end{tikzpicture}}}
    \hfill
    \subfloat[\raggedright One iteration of a simplex method. The next simplex will be formed by $\{\xx_2,\xx_3,\xx_4\}$. ]{\label{fig:NelderMead}\resizebox{0.25\linewidth}{!}{\begin{tikzpicture}
    \filldraw[black] (-0.3,1) circle (2pt) node[anchor=east] {$\xx_1$}; 
    \filldraw[black] (-1.1,-0.5) circle (2pt) node[anchor=south east] {$\xx_2$}; 
    \filldraw[black] (1,0) circle (2pt) node[anchor=north] {$\xx_3$}; 
    
    \draw[thick] (1,0) -- (-0.3,1);
    \draw[thick] (-0.3,1) -- (-1.1,-0.5);
    \draw[thick] (-1.1,-0.5) -- (1,0);
    
    \filldraw[black] (0.25, -1.75) circle (2pt) node[anchor=north] {$\xx_4$}; 
    \draw[dashed] (-0.3,1) -- (0.25, -1.75); 
\end{tikzpicture}}}
   \caption{An illustration of two DFO algorithms when minimizing a bivariate function, where $f(\xx_1) >$ $f(\xx_2) >$ $f(\xx_3)$.  The vertices of the triangles represents $\Theta$. This figure only illustrates the algorithms' behavior when the trial point $\xx_4$ satisfies $f(\xx_4) < f(\xx_3)$.}
   \label{fig:DFO}
\end{figure}

The theories developed in this paper are intended for linear interpolation in general and can be applied beyond model-based DFO methods. 
In some cases, they can even be used for analyzing numerical methods that do not involve linear interpolation. 
One such case is a class of DFO methods known as the {\it simplicial (direct) search methods} or {\it simplex methods}, with the Nelder-Mead method \cite{nelder1965simplex} being its most famous variant.
As illustrated in Figure~\ref{fig:NelderMead}, the main routine of these algorithms involves taking a set of $n+1$ affinely independent points $\Theta$ (the vertices of a simplex) and reflecting the one with the largest function value through the hyperplane defined by the rest. 
While linear interpolation is not used in these algorithms, the range of the function value at the reflection point (see $\xx_4$ in Figure~\ref{fig:NelderMead} which lies outside $\conv(\Theta)$) can be determined by the sum of the value estimated by the affine function interpolating $\Theta$ and the bound on the estimation error. 
We provide a concrete example of such application in Section~\ref{sec:app} by deriving the complexity of a basic variant of the simplicial search method that only performs the reflection operation. 
To the best of our knowledge, this is the first time the complexity of any of these methods has ever been determined. 

The main contributions of this paper are as follows. 
\begin{enumerate}
    \item We formulate the problem of finding the sharp error bound as a nonlinear programming problem and show that it can be solved numerically to obtain the desired bound. 
    \item An analytical bound on the function approximation error is derived and proved to be sharp for interpolation and, under certain conditions, for extrapolation. 
    \item We present the largest function approximation error that is achievable by quadratic functions in $C_\nu^{1,1}(\R^n)$, as well as the condition under which it is an upper bound on the error achievable by all functions in $C_\nu^{1,1}(\R^n)$. 
    \item For the bivariate ($n=2$) case, we study linear interpolation's function approximation error under all possible geometric configurations of $\Theta$ and $\xx_0$ and provide the corresponding sharp bounds. We also prove piecewise quadratic functions can achieve the largest approximation error in this two-dimensional case. 
    \item Using one of the error bounds developed in this paper, the complexity of a simplicial search method is determined for the first time. 
\end{enumerate}

The paper is organized as follows. 
We review the existing results and the related literature in the rest of this section. 
Our notation and the preliminary knowledge are introduced in Section~\ref{sec:preliminaries}. 
The nonlinear programming problem is present in Section~\ref{sec:numerical}. 
In Section~\ref{sec:phase1}, we generalize and improve an existing analytical bound. 
In Section~\ref{sec:phase2}, we study the error in approximating quadratic functions. 
In Section~\ref{sec:phase3}, we show how to calculate the sharp bound on function approximation error of bivariate linear interpolation. 
The complexity analysis of a basic simplicial search method is presented in Section~\ref{sec:app}. 
We conclude the paper in Section~\ref{sec:discussion} by discussing our findings and some open questions. 

\subsection{Related Work}

The function approximation error of univariate ($n=1$) interpolation using polynomials of any degree is already well-studied, and the results can be found in classical literature such as \cite{davis1975book}. 
If a $(d+1)$-times differentiable function $f$ defined on $\R$ is interpolated by a polynomial of degree $d$ on $d+1$ unique points $\{x_1, x_2, \dots, x_{d+1}\} \subset \R$, then the resulting polynomial has the approximation error 
\begin{equation} \label{eq:Cauchy remainder}
    \frac{(x-x_1)(x-x_2)\cdots(x-x_{d+1})}{(d+1)!} D^{n+1} f(\xi) \quad \text{for all } x \in \R
\end{equation} 
for some $\xi$ with $\min(x,x_1,,\dots,x_{d+1}) < \xi < \max(x,x_1,\dots,x_{d+1})$. 
Unfortunately this result cannot be extended to the multivariate ($n>1$) case directly, even if the polynomial is linear ($d=1$). 

The function approximation error of multivariate polynomial interpolation has been studied by researchers from multiple research fields. 
Motivated by their application in finite element methods, formulae for the errors in both Lagrange and Hermite interpolation with polynomials of any degree are derived in \cite{ciarlet1972general}. 
As part of an effort to develop derivative-free optimization algorithms, a bound on the error of quadratic interpolation is provided in \cite{powell2001lagrange}. 
The sharp error bound for linear interpolation is found by researchers in approximation theory in \cite{waldron1998error} for the case when $\xx \in \conv(\Theta)$ using the unique Euclidean sphere that contains $\Theta$. 
Following \cite{waldron1998error}, a number of sharp error bounds are derived in \cite{stampfle2000optimal} for linear interpolation under several different smoothness or continuity assumptions in addition to \eqref{eq:Lipschitz}. 

Following Powell's development of model-based DFO methods, researchers have been studying their convergence properties. 
An early investigation is presented in \cite{conn1997convergence}, where the Newton polynomials are used as a tool to maintain the geometry of the interpolation set and to bound the approximation error. 
This tool is later reverted in \cite{conn2008geometry} to the Lagrange polynomials, which Powell originally used in his DFO methods. 
It is proven in \cite{conn2008geometry} that, with a proper maintenance of the interpolation set, the model's error in approximating the objective function's zeroth-, first-, and second-order derivatives can be bounded by polynomials of the trust region radius. 
Models satisfying such error bounds are referred to as {\it fully-linear} or {\it fully-quadratic} models, and trust region DFO methods are shown to converge globally in \cite{conn2009global} while assuming these model qualities. 
Constructing interpolation models using random samples is proposed and analyzed in \cite{bandeira2012computation} and \cite{bandeira2014convergence}, where the models are only considered to be fully-linear or fully-quadratic with sufficiently high probability in each iteration of the DFO algorithm. 
Later theoretical works, including but not limited to \cite{gratton2018complexity}, \cite{chen2018stochastic}, \cite{blanchet2019convergence}, and \cite{cao2023first}, manage to obtain better convergence results of the trust region methods under weaker assumptions on the approximation models but have limited discussion on how these models can be obtained from sample sets. 

Compared with model-based DFO methods, theories regarding the performance of simplicial search methods are less developed.
The Nelder-Mead method is shown to fail to converge to a stationary point in some examples in \cite{woods1985interactive} and \cite{mckinnon1998convergence}. 
However, by forbidding the expansion operation, it is shown in \cite{lagarias2012convergence} to converge to the minimizer of any two-dimensional twice continuously differentiable function with positive definite Hessian and bounded level sets. 
A method with a fortified descent condition is proposed in \cite{tseng1999fortified} and every cluster point of the sequence of its generated candidate points is proved to be stationary. 
There are several more convergent variants of the simplicial direct search methods, but most of them are embedded within other algorithmic frameworks.
For a comprehensive review of this subject, please refer to Section~2.1.1 of \cite{larson2019derivative}. 

Our theoretical analysis uses results from the study of the \textit{performance estimation problem} (PEP). 
First proposed in \cite{drori2014performance}, a PEP is a nonlinear programming formulation of the problem of finding an optimization algorithm's worst-case performance over a set of possible objective functions. 
It involves maximizing a performance measure of the given algorithm (the larger the measure, the worse the performance) over all possible objective functions and is an infinite-dimensional problem. 
However, with some algorithms and functions, particularly first-order nonlinear optimization methods and convex functions, the PEP has been shown to have finite-dimensional equivalents that can be solved numerically \cite{taylor2017smooth,taylor2017exact}, thus providing a computer-aided analysis tool for estimating an algorithm's worst-case performance. 
For the function class $C_\nu^{1,1}(\R^n)$, the following proposition is the key to find the finite-dimensional equivalents of PEP and will be used in our analysis. 

\begin{proposition}[{\cite[Theorem 3.10]{taylor2017exact}}]\label{prop:functional}
Let $\nu > 0$ and $\cI$ be an index set, and consider a set of triples $\{(\xx_i,\g_i,y_i)\}_{i\in\cI}$ where $\xx_i\in\R^n$, $\g\in\R^n$, and $y_i\in\R$ for all $i\in\cI$. 
There exists a function $f\in C_\nu^{1,1}(\R^n)$ such that both $\g_i = Df(\xx_i)$ and $y_i = f(\xx_i)$ hold for all $i\in\cI$ if and only if the following inequality holds for all $i,j\in\cI$: 
\begin{equation} \label{eq:Lipschitz stronger ij} \begin{aligned}
    y_j \le y_i + \frac{1}{2} (\g_i + \g_j) \cdot (\xx_j - \xx_i) + \frac{\nu}{4} \|\xx_j-\xx_i\|^2 - \frac{1}{4\nu} \|\g_j - \g_i\|^2. 
\end{aligned} \end{equation} 
\end{proposition}

\section{Notation and Preliminaries} \label{sec:preliminaries}
Since this paper involves approximation theory and optimization, to appeal to audiences from both research fields, we provide a detailed introduction to our notation and the preliminary knowledge. 

We use boldface lowercase letters for vectors and uppercase letters for matrices. 
The element in the $i$th row and $j$th column of a matrix $\bU$ will be denoted by $[\bU]_{ij}$ or $[\bU]_{i,j}$. 
All one and all zero vectors are written as $\mathbf{1}$ and $\mathbf{0}$, while $\e_i$ is the unit vector along the $i$th coordinate. 
The identity matrix and the all one square matrix are denoted by $\bI$ and $\bJ$, respectively.
Subscripts are sometimes added to $\mathbf{1}$, $\mathbf{0}$, $\bI$, or $\bJ$ to mark their dimension. 
For example, $\bJ_{n+1}$ would be the all one square matrix of order $n+1$. 
We denote by $\|\cdot\|$ the Euclidean norm. 
The dot product between vectors or matrices of the same dimension, $\u\cdot\v$ or $\bU \cdot \bV$, is the summation of the entry-wise products, which are customarily denoted by $\u^T\v$ and Tr$(\bU^T \bV)$ in optimization literature. 
 
Given $\xx_0\in\R^n$ and an affinely independent sample set $\Theta:= \{\xx_1,\xx_2,\dots,\xx_{n+1}\} \subset \R^n$, we define matrices 
\begin{equation} \label{eq:Y Phi}
    \bY = \begin{bmatrix}
        (\xx_1-\xx_0)^T \\
        (\xx_2-\xx_0)^T \\
        \vdots \\
        (\xx_{n+1}-\xx_0)^T 
    \end{bmatrix} \in \R^{(n+1)\times n}
    \text{ and } 
    \bPhi = \begin{bmatrix} \mathbf{1} &\bY \end{bmatrix} \in \R^{(n+1)\times (n+1)}. 
\end{equation}
The affine function $\hat{f}(\bu) = c + \bg \cdot (\bu - \xx_0)$ that interpolates $f$ on $\Theta$, where $c\in\R$ and $\bg \in \R^n$, can be determined by solving the linear system 
\begin{equation}
    \bPhi \begin{pmatrix}
        c\\ \bg 
    \end{pmatrix}
    = \begin{pmatrix}
        f(\xx_1)\\ \vdots\\ f(\xx_{n+1})
    \end{pmatrix}, 
\end{equation}
where the left-hand side is essentially $\begin{pmatrix} \hat{f}(\xx_1) &\cdots &\hat{f}(\xx_{n+1}) \end{pmatrix}^T$. 
The affine independence of $\Theta$ implies the nonsingularity of $\bPhi$, which further implies the uniqueness of $\hat{f}$. 

With the $\Theta$, $Y$, $\Phi$ and $\hat{f}$ defined above, let $\ell_1, \dots, \ell_{n+1}$ be the {\it Lagrange polynomials}, i.e., the unique set of polynomials such that $\ell_i(\xx_j) = 1$ if $i=j$, and $\ell_i(\xx_j) = 0$ if $i\neq j$. 
From a linear algebra perspective, this simply means they are the set of affine functions defined by $\Theta$ as 
\begin{equation} \begin{aligned} 
    \ell_i(\u) 
    = [\bPhi^{-1}]_{1i} + \left(\begin{bmatrix} \mathbf{0} &\bI_n \end{bmatrix} \bPhi^{-1} \e_i\right) \cdot (\bu - \xx_0) \text{ for } i = 1,\dots,n+1. 
\end{aligned} \end{equation}
We additionally define a constant function $\ell_0$ such that $\ell_0(\u) = -1$ for all $\u\in\R^n$.
Then, the most essential properties of the Lagrange polynomials can be written as 
\begin{align}
    \sum_{i=1}^{n+1} \ell_i(\xx_0) f(\xx_i) &= \hat{f}(\xx_0), \label{eq:Lagrange m}  \\
    \sum_{i=1}^{n+1} \ell_i(\xx_0) &= 1 \quad (\text{or, equivalently,} \sum_{i=0}^{n+1} \ell_i(\xx_0) = 0), \label{eq:Lagrange 0} \\
    \text{and } \sum_{i=1}^{n+1} \ell_i(\xx_0) \xx_i &= \xx_0 \quad (\text{or, equivalently,} \sum_{i=0}^{n+1} \ell_i(\xx_0) \xx_i = \mathbf{0}) \label{eq:Lagrange Y}, 
\end{align}
where \eqref{eq:Lagrange Y} can also be written as $Y^T \begin{bmatrix} \ell_1(\xx_0) &\cdots &\ell_{n+1}(\xx_0)\end{bmatrix}^T = \mathbf{0}$ in matrix form. 
The concept of Lagrange polynomials is fundamental to approximation theory. 
The book \cite{DFO_book} offers a comprehensive introduction to them in the context of derivative-free optimization. 

Without loss of generality, we assume the set $\Theta = \{\xx_1,\xx_2,\dots,\xx_{n+1}\}$ is ordered in a way such that $\ell_1(\xx_0) \ge \ell_2(\xx_0) \ge \cdots \ge \ell_{n+1}(\xx_0)$. 
We define the following two sets of indices: 
\begin{subequations} \begin{align} 
\cI_+ &= \{i\in \{0,1,\dots,n+1\}:~ \ell_i(\xx_0)>0\} = \{1,2,\dots, |\cI_+|\}, \\
\cI_- &= \{i\in \{0,1,\dots,n+1\}:~ \ell_i(\xx_0)<0\} \\
    &= \{0,~ n-|\cI_-|+3,~ n-|\cI_-|+4, \dots,~ n+1\}. \nonumber 
\end{align} \end{subequations}
Notice \eqref{eq:Lagrange 0} implies $\cI_+ \neq \emptyset$, and  $\ell_0(\xx_0)=-1$ implies $\cI_- \neq \emptyset$. 
Furthermore, $\cI_+ \cap \cI_- = \emptyset$ and $\cI_+ \cup \cI_- \subseteq \{0,1,\dots,n+1\}$. 
The two sets $\cI_+$ and $\cI_-$ do not always define a partition of $\{0,1,\dots,n+1\}$ since $\ell_i(\xx_0) = 0$ can hold for one or more $i \in \{1,2,\dots,n+1\}$. 
For example, if $(\ell_0(\xx_0),\ell_1(\xx_0),\ell_2(\xx_0),\ell_3(\xx_0),\ell_4(\xx_0),\ell_5(\xx_0)) = (-1,4,0,0,-1,-2)$, then $n = 4$, $\cI_+  = \{1\}$, and $\cI_- = \{0,4,5\}$.

The class of functions $C_\nu^{1,1}(\R^n)$ is ubiquitous in the research of nonlinear optimization. 
It is well-known (see, e.g., Section 1.2.2 of the textbook \cite{Nesterov_book}) that the inclusion $f \in C_\nu^{1,1}(\R^n)$ implies 
\begin{equation} \label{eq:Lipschitz quadratic}
    |f(\v) - f(\u) - Df(\u) \cdot (\v - \u)| \le \frac{\nu}{2} \|\v - \u\|^2 \text{ for all } \u,\v \in \R^n, 
\end{equation}
and that if $f$ is twice differentiable on $\R^n$, \eqref{eq:Lipschitz} and \eqref{eq:Lipschitz quadratic} are equivalent to 
\begin{equation} \label{eq:Lipschitz Hessian}
    -\nu I \preceq D^2 f(\u) \preceq \nu I \text{ for all } \u \in \R^n, 
\end{equation}
where the condition \eqref{eq:Lipschitz Hessian} is often written as $\|~|D^2 f|~\|_{L_\infty(\R^n)} \le \nu$ in approximation theory literature. 
What is less well-known about the class $C_\nu^{1,1}(\R^n)$ is that $f \in C_\nu^{1,1}(\R^n)$ also implies 
\begin{equation} \label{eq:Lipschitz stronger} \begin{aligned}
    f(\v) \le{} &f(\u) + \frac{1}{2} (Df(\u) + Df(\v)) \cdot (\v - \u) \\ 
    &+ \frac{\nu}{4} \|\v-\u\|^2 - \frac{1}{4\nu} \|Df(\v) - Df(\u)\|^2 \text{ for all } \u,\v \in \R^n. 
\end{aligned} \end{equation} 
For differentiable functions, \eqref{eq:Lipschitz}, \eqref{eq:Lipschitz quadratic}, and \eqref{eq:Lipschitz stronger} are equivalent, as shown in the following proposition. 
\begin{proposition}
    Assume $f:\R^n \rightarrow \R$ is differentiable. Then, the conditions \eqref{eq:Lipschitz}, \eqref{eq:Lipschitz quadratic}, and \eqref{eq:Lipschitz stronger} are equivalent. 
\end{proposition}
\begin{proof}
The implication \eqref{eq:Lipschitz}$\rightarrow$\eqref{eq:Lipschitz quadratic} is already established in \cite[Section 1.2.2]{Nesterov_book}. 
Now we prove \eqref{eq:Lipschitz quadratic}$\rightarrow$\eqref{eq:Lipschitz stronger}.
Assume \eqref{eq:Lipschitz quadratic} holds. 
For any $\u,\v,\w \in \R^n$, there are 
\begin{align*}
    f(\w) &\stackrel{\eqref{eq:Lipschitz quadratic}}{\ge} f(\v) + Df(\v)\cdot(\w-\v) - \frac{\nu}{2} \|\w-\v\|^2 \\
    \text{and } f(\w) &\stackrel{\eqref{eq:Lipschitz quadratic}}{\le} f(\u) + Df(\u)\cdot(\w-\u) + \frac{\nu}{2} \|\w-\u\|^2.  
\end{align*}
Then, the right-hand side of the first inequality is less than or equal to that of the second one: 
\[ f(\v) - f(\u) + Df(\v)\cdot(\w-\v) - Df(\u)\cdot(\w-\u) - \frac{\nu}{2}\left(||\w-\v||^2+||\w-\u||^2\right) \le 0. 
\]
Setting $\w=(\v+\u)/2 + (Df(\v) - Df(\u)) / (2\nu)$ in the above inequality yields \eqref{eq:Lipschitz stronger}. 

Finally, \eqref{eq:Lipschitz} can be obtained from \eqref{eq:Lipschitz stronger} by adding together \eqref{eq:Lipschitz stronger} and another copy of \eqref{eq:Lipschitz stronger} with $\u$ and $\v$ swapped, multiplying both sides by $2\nu$, and then taking their square roots. 
\end{proof}


\section{Error Estimation Problem} \label{sec:numerical}
In this section, we formulate the problem of finding the sharp error bound as a numerically solvable nonlinear optimization problem. 
We first make the important observation that the problem of finding the sharp upper bound on the error is the same as asking for the largest error that a function from $C_\nu^{1,1}(\R^n)$ can achieve. 
Thus, it can be formulated as the following problem of maximizing the approximation error over the functions in $C_\nu^{1,1}(\R^n)$: 
\begin{equation} \label{prob:D} \tag{EEP} \everymath{\displaystyle} \begin{aligned} 
    &z(\Theta,\xx_0) = &\max_f &&|\hat{f}(\xx_0) - f(\xx_0)| \\
    &&\text{s.t. } &&f \in C_\nu^{1,1}(\R^n), 
\end{aligned} 
\end{equation}
where $\hat{f}$ is the affine function that interpolates $f$ on the given set of $n+1$ affinely independent points $\Theta = \{\xx_1,\dots,\xx_{n+1}\}$. 
We call this problem the \textit{error estimation problem} (EEP), a name inspired by the \textit{performance estimation problem} (PEP). 

By \eqref{eq:Lagrange m}, the objective function in \eqref{prob:D} can be written as $|\sum_{i=0}^{n+1} \ell_i(\xx_0) f(\xx_i)|$. 
Moreover, the absolute sign can be dropped thanks to the symmetry of \eqref{eq:Lipschitz}, that is, $-f \in C_\nu^{1,1}(\R^n)$ if and only if $f \in C_\nu^{1,1}(\R^n)$, and the approximation error on the two functions $f$ and $-f$ are negatives of each other. 
This simplifies the objective to $\sum_{i=0}^{n+1} \ell_i(\xx_0) f(\xx_i)$. 
However, \eqref{prob:D} still involves a functional variable $f$ and is an infinite-dimensional problem. 
This is where we need to use Proposition~\ref{prop:functional}, according to which the constraint $f \in C_\nu^{1,1}(\R^n)$ in \ref{prob:D} holds if and only if \eqref{eq:Lipschitz stronger} holds for every pair of points in $\Theta \cup \{\xx_0\}$. 
Thus, by letting $y_i = f(\xx_i)$ and $\g_i = D f(\xx_i)$ for $i = 0,1,\dots,n+1$, we obtain the following finite-dimensional equivalent form of \eqref{prob:D}: 
\begin{equation} \label{prob:f-D} \tag{f-EEP} \everymath{\displaystyle} \begin{aligned}
    z(\Theta,\xx_0) = &\max_{y_i,\g_i} &&\sum_{i=0}^{n+1}\ell_i(\xx_0) y_i \\
    &\text{s.t. } &&y_j \le y_i + \frac{1}{2} (\g_i + \g_j) \cdot (\xx_j - \xx_i) + \frac{\nu}{4} \|\xx_j-\xx_i\|^2 \\
    &~ &&\qquad - \frac{1}{4\nu} \|\g_j - \g_i\|^2 \quad \forall i,j\in\{0,\dots,n+1\}. 
\end{aligned} 
\end{equation}

Similar to PEP, \eqref{prob:D} is an analytical tool. It is meant to be solved analytically (rather than numerically) to obtain the closed form solution of the sharp error bound $z(\Theta,\xx_0)$. 
However, such a solution is too difficult to find directly.
Therefore, solving its finite-dimensional equivalent \eqref{prob:f-D} numerically becomes a crucial step in obtaining the insight that will guide the subsequent discovery of the final analytical solution. 
Currently, we have \eqref{prob:f-D}, a convex quadratically constrained quadratic program (QCQP). 
While this type of problem is generally tractable, \eqref{prob:f-D} contains $n+1$ redundant degrees of freedom, resulting in infinitely many optimal solutions that can challenge standard solvers.

There are many approaches to eliminate these degrees of freedom. 
For example, one can fix $\{y_i\}_{i=1}^{n+1}$ in \eqref{prob:f-D} to their observed values. 
Indeed, these function values are needed for constructing the affine approximation $\hat{f}$, so it is natural to assume they are known. 
However, we note that the optimal value $z(\Theta,\xx_0)$ is affected by $\Theta$ but not by the observed function values at these points. 
Thus, for the purpose of solving \eqref{prob:f-D}, it is justified to simply set $y_i = 0$ for all $i=1,\dots,n+1$. 
Alternatively, one can also fix $(\g_i,y_i)$ to $(\mathbf{0}, 0)$ for any $i \in \{0,1,\dots,n+1\}$. 
We formally prove in the following proposition the $n+1$ degrees of freedom can be removed in these two ways. 

\begin{proposition}
    The following statements are true. 
    \begin{enumerate}
        \item If any function $f$ is optimal to \eqref{prob:D}, then, with any $c\in\R$ and $\g\in\R^n$, the function $h(\u) = f(\u) + c + \g\cdot\u$ is also optimal. 
        \item The optimal value of \eqref{prob:f-D} with additional constraints that fix $\{y_i\}_{i=1}^{n+1}$ to any values is still $z(\Theta,\xx_0)$. 
        \item The optimal value of \eqref{prob:f-D} with additional constraints that, for one and only one of $k \in \{0,1,\dots,n+1\}$, fix the pair $(\g_k,y_k)$ to any values is still $z(\Theta,\xx_0)$. 
    \end{enumerate}
\end{proposition}

\begin{proof}
    By the definition \eqref{eq:Lipschitz} in the first paragraph of this paper, it is easy to see $h \in C_\nu^{1,1}(\R^n)$ whenever $f\in C_\nu^{1,1}(\R^n)$. The two objective values can also be shown to be the same using \eqref{eq:Lagrange m}, \eqref{eq:Lagrange 0}, and \eqref{eq:Lagrange Y}: 
    \[ \sum_{i=0}^{n+1} \ell_i(\xx_0) h(\xx_i) = \sum_{i=0}^{n+1} \ell_i(\xx_0) [f(\xx_i) + c + \g\cdot\xx_i] = \sum_{i=0}^{n+1} \ell_i(\xx_0) f(\xx_i). \] 
    The first statement thus holds true. 

    To prove the second statement, we first assume \eqref{prob:f-D} has an optimal solution $\{y_i^\star, \g_i^\star\}_{i=0}^{n+1}$. 
    Now suppose the problem has an additional set of constraints $y_i = \bar{y}_i$ for $i=1,\dots,n+1$, where $\{\bar y_i\}_{i=1}^{n+1}$ is an arbitrary set of real numbers. 
    Let $(c,\g)$ be the solution to the linear system $c + \g\cdot(\xx_i - \xx_0) = \bar y_i - y_i^\star, i=1,\dots,n+1$, which is uniquely solvable due to the nonsingularity of $\Phi$ (defined in \eqref{eq:Y Phi}). 
    Then, this new problem has an optimal solution $\g_i = \g_i^\star + \g$ for all $i=0,1,\dots,n+1$ and $y_0 = \bar{y}_0 \stackrel{\rm def}{=} y_0^\star + c$, and the corresponding optimal value is the same as that of the original \eqref{prob:f-D}. 
    Indeed, the constraints of this new problem are satisfied as 
    \[ \begin{aligned} 
        &- y_j + y_i + \frac{1}{2} (\g_i + \g_j) \cdot (\xx_j - \xx_i) + \frac{\nu}{4} \|\xx_j-\xx_i\|^2 - \frac{1}{4\nu} \|\g_j - \g_i\|^2 \\
        &= - \bar{y}_j + \bar{y}_i + \frac{1}{2} (\g_i^\star + \g_j^\star + 2\g) \cdot (\xx_j - \xx_i) + \frac{\nu}{4} \|\xx_j-\xx_i\|^2 - \frac{1}{4\nu} \|\g_j^\star - \g_i^\star\|^2 \\
        &= - y_j^\star + y_i^\star + \frac{1}{2} (\g_i^\star + \g_j^\star) \cdot (\xx_j - \xx_i) + \frac{\nu}{4} \|\xx_j-\xx_i\|^2 - \frac{1}{4\nu} \|\g_j^\star - \g_i^\star\|^2 
        \ge 0
    \end{aligned} \]
    for all $i,j = 0,1,\dots,n+1$, 
    where the second equality is true because $\g\cdot\xx_j - \g\cdot\xx_i = (\bar{y}_j - y_j^\star - c + \bg \cdot \xx_0) - (\bar{y}_i - y_i^\star - c + \bg \cdot \xx_0) = (\bar{y}_j - y_j^\star) - (\bar{y}_i - y_i^\star)$; and the objective function 
    \[ \begin{aligned}  
        \sum_{i=0}^{n+1} \ell_i(\xx_0) y_i 
        &= -y_0^\star - c + \sum_{i=1}^{n+1} \ell_i(\xx_0) \bar{y}_i 
        \stackrel{\eqref{eq:Lagrange 0}}{=} -y_0^\star + \sum_{i=1}^{n+1} \ell_i(\xx_0) [\bar{y}_i - c] \\
        &\leftstackrel{\eqref{eq:Lagrange 0}\eqref{eq:Lagrange Y}}{=} -y_0^\star + \sum_{i=1}^{n+1} \ell_i(\xx_0) [\bar{y}_i - c - \g\cdot(\xx_i-\xx_0)]
        = \sum_{i=0}^{n+1} \ell_i(\xx_0) y_i^\star. 
    \end{aligned} \]
    Similarly, \eqref{prob:f-D} with the additional constraint $(\g_k,y_k) = (\bar{\g}_k, \bar{y}_k)$ for some $k\in\{0,1,\dots,n+1\}$ and any $(\bar{\g}_k, \bar{y}_k) \in \R^n\times\R$ also has the same optimal value as \eqref{prob:f-D}, and its optimal solution satisfies $\g_i = \g_i^\star - \g_k^\star + \bar \g_k$ and $y_i = y_i^\star - y_k^\star + \bar y_k + (\bar \g_k - \g_k^\star)\cdot(\xx_i-\xx_k)$ for all $i = 0,1,\dots,n+1$. 
    The constraints are satisfied as 
    \[ \begin{aligned} 
        &- y_j + y_i + \frac{1}{2} (\g_i + \g_j) \cdot (\xx_j - \xx_i) + \frac{\nu}{4} \|\xx_j-\xx_i\|^2 - \frac{1}{4\nu} \|\g_j - \g_i\|^2 \\
        &= - [y_j^\star - y_k^\star + \bar y_k + (\bar \g_k - \g_k^\star)\cdot(\xx_j-\xx_k)] + [y_i^\star - y_k^\star + \bar y_k + (\bar \g_k - \g_k^\star) \cdot (\xx_i-\xx_k)] \\
        &\quad + \frac{1}{2} (\g_i^\star + \g_j^\star - 2\g_k^\star + 2\bar \g_k) \cdot (\xx_j - \xx_i) + \frac{\nu}{4} \|\xx_j-\xx_i\|^2 - \frac{1}{4\nu} \|\g_j^\star - \g_i^\star\|^2 \\
        &= - y_j^\star + y_i^\star + \frac{1}{2} (\g_i^\star + \g_j^\star) \cdot (\xx_j - \xx_i) + \frac{\nu}{4} \|\xx_j-\xx_i\|^2 - \frac{1}{4\nu} \|\g_j^\star - \g_i^\star\|^2 
        \ge 0
    \end{aligned} \]
    for all $i,j = 0,1,\dots,n+1$; 
    and the objective function 
    \[  \sum_{i=0}^{n+1} \ell_i(\xx_0) y_i 
        = \sum_{i=0}^{n+1} \ell_i(\xx_0) [y_i^\star - y_k^\star + \bar y_k + (\bar \g_k - \g_k^\star)\cdot(\xx_i-\xx_k)]
        \stackrel{\eqref{eq:Lagrange 0}\eqref{eq:Lagrange Y}}{=} \sum_{i=0}^{n+1} \ell_i(\xx_0) y_i^\star. 
    \]
\end{proof}

\begin{figure}[ht]
    \centering
    \includegraphics[width=\linewidth]{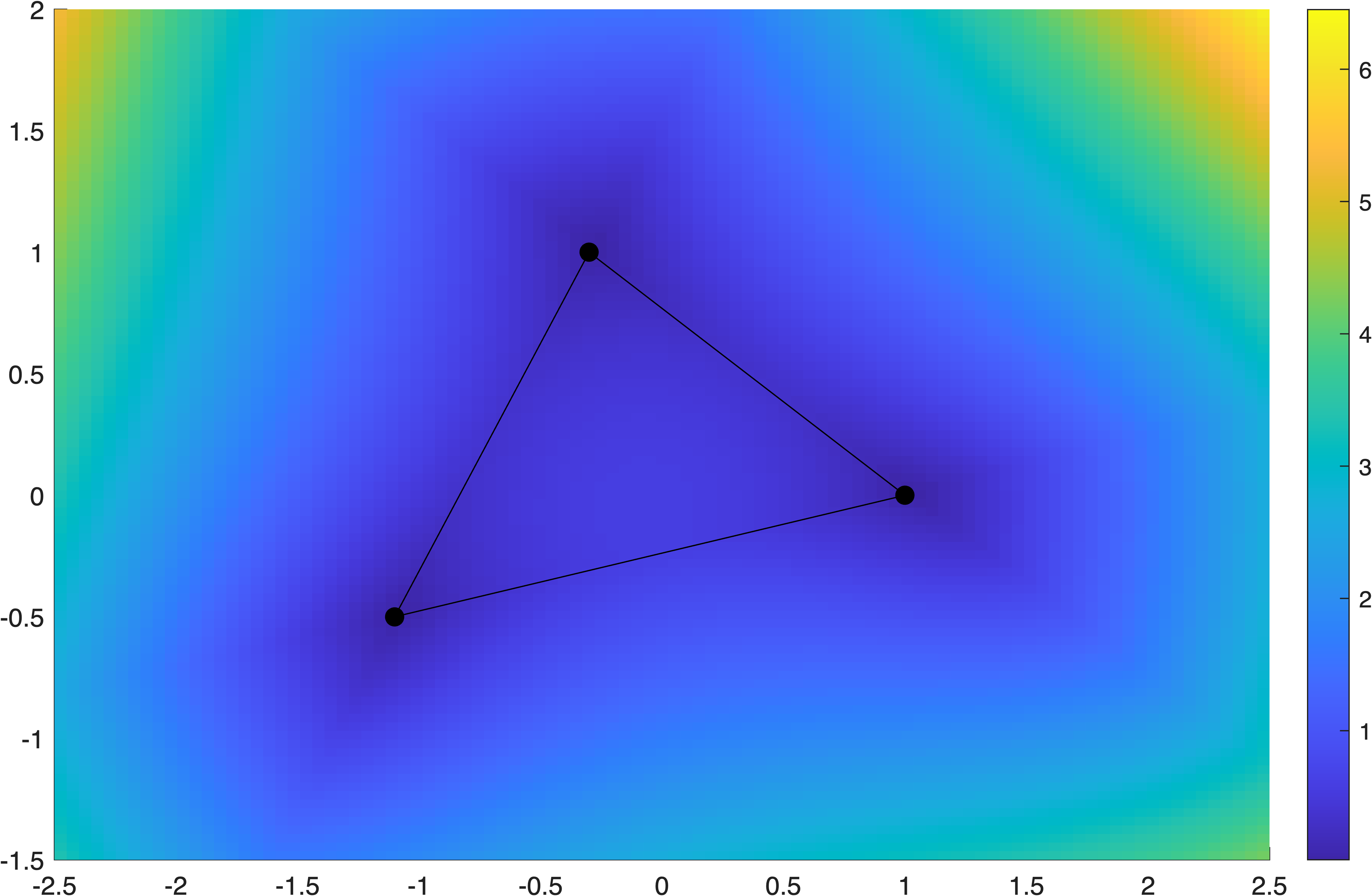}
    \caption{The sharp error bound $z(\Theta,\xx_0)$ as a function of $\xx_0$.The sample set and the Lipschitz constant are chosen as $\Theta = \{(-0.3,1), (-1.1,-0.5), (1,0)\}$ and $\nu = 1$.}
    \label{fig:numerical}
\end{figure}

Utilizing the above proposition, \eqref{prob:f-D} can be solved by most off-the-shelf convex optimization solvers. 
In our experiments, we first set $\g_0 = \mathbf{0}$ and $y_0 = 0$ to eliminate the extra degrees of freedom, and then solve \eqref{prob:f-D} with the {\tt fmincon} function in Matlab with the default {\tt interior-point} solver. 
With the ability to solve \eqref{prob:f-D} numerically, we now use it to visualize the sharp error bound for bivariate linear interpolation. 
Firstly, the affinely independent set $\Theta \subset \R^2$ is fixed to $\{(-0.3,1), (-1.1,-0.5), (1,0)\}$. 
Then, the problem \eqref{prob:f-D} is repeatedly solved to obtain $z(\Theta,\xx_0)$ for each point $\xx_0$ on a $100\times100$ grid that covers the area $[-2.5,2.5]\times[-1.5,2.5]$ evenly. 
This process gives us the values of $z(\Theta,\xx_0)$ (as a function of $\xx_0$) at 10000 different points. 
We plot the result in Figure~\ref{fig:numerical} using the {\tt surf} function in Matlab, where value of the sharp error bound $z(\Theta,\xx_0)$ at each point $\xx_0$ is indicated by the color of the pixel at $\xx_0$.  

\section{An Improved Upper Bound} \label{sec:phase1}

We now begin our attempt at finding the analytical form of the bound $z(\Theta,\xx_0)$, where both $\Theta$ and $\xx_0$ are considered to be fixed. 
The theoretical results in \cite{ciarlet1972general} and \cite{powell2001lagrange} are obtained by comparing $f$ against its Taylor expansion at $\xx_0$. 
We generalize their approach in Theorem~\ref{thm:phase1} by using the Taylor expansion of $f$ at an arbitrary $\u \in \R^n$. 
\begin{theorem}\label{thm:phase1}
    Assume $f \in C^{1,1}_\nu(\R^n)$. 
    Let $\xx_0$ be a fixed vector in $\R^n$, and $\Theta = \{\xx_1,\dots,\xx_{n+1}\}\subset \R^n$ be a fixed set of $n+1$ affinely independent vectors. 
    Let $\hat{f}$ be the affine function that interpolates $f$ on $\Theta$. 
    The inequality  
    \begin{equation} \label{eq:phase1 u}
        |\hat{f}(\xx_0) - f(\xx_0)| \le \frac{\nu}{2} \sum_{i=0}^{n+1} |\ell_i(\xx_0)| \|\xx_i-\u\|^2, 
    \end{equation}
    holds for all $\u \in \R^n$. 
\end{theorem}

\begin{proof}
By \eqref{eq:Lipschitz quadratic}, we have for any $\u \in \R^n$ 
\begin{subequations} \label{phase1 set of inequalities} \begin{align}
\ell_i(\xx_0) [f(\xx_i) - f(\u) - Df(\u)\cdot (\xx_i-\u)] &\le \ell_i(\xx_0) \frac{\nu}{2} \|\xx_i-\u\|^2, \forall i\in \cI_+, \\
-\ell_i(\xx_0) [-f(\xx_i) + f(\u) + Df(\u)\cdot(\xx_i-\u)] &\le -\ell_i(\xx_0) \frac{\nu}{2} \|\xx_i-\u\|^2, \forall i\in \cI_-. 
\end{align} \end{subequations}
Now add all inequalities above together. The sum of the left-hand sides is 
\[ \begin{aligned} 
    &\sum_{i=0}^{n+1} \ell_i(\xx_0) [f(\xx_i) - f(\u)] + Df(\u) \cdot \sum_{i=0}^{n+1} \ell_i(\xx_0) [\u-\xx_i] \\
    &\stackrel{\eqref{eq:Lagrange 0}}{=} \sum_{i=0}^{n+1} \ell_i(\xx_0) f(\xx_i) - Df(\u) \cdot \sum_{i=0}^{n+1} \ell_i(\xx_0) \xx_i 
    \stackrel{\eqref{eq:Lagrange Y}}{=} \sum_{i=0}^{n+1} \ell_i(\xx_0) f(\xx_i)
    \stackrel{\eqref{eq:Lagrange m}}{=}  \hat{f}(\xx_0) - f(\xx_0), 
\end{aligned} \] 
while the sum of the right-hand sides is $\nu/2 \sum_{i=0}^{n+1} |\ell_i(\xx_0)| \|\xx_i-\u\|^2$. 
Thus the sum of the inequalities in \eqref{phase1 set of inequalities} is \eqref{eq:phase1 u} when $\hat{f}(\xx_0)-f(\xx_0) \ge 0$. 
If the inequalities in \eqref{phase1 set of inequalities} have their left-hand sides multiplied by $-1$, they would still hold according to \eqref{eq:Lipschitz quadratic}, and their summation would be \eqref{eq:phase1 u} for the $\hat{f}(\xx_0)-f(\xx_0) < 0$ case. 
\end{proof}

A result similar to Theorem~\ref{thm:phase1} has been presented in \cite{ciarlet1972general}, but it is only proved there that \eqref{eq:phase1 u} holds when $\u=\xx_0$. 
In comparison, the new bound provides more convenience in analyzing DFO algorithms that use trust region methods, since the free point $\u$ can be set to the center of the trust region. 
Another advantage of the new bound is that it can be minimized with respect of $\u$ to obtain the improved bound \eqref{eq:phase1}. 
\begin{corollary} \label{cor:phase1}
    Under the setting of Theorem~\ref{thm:phase1}, the function approximation error $|\hat{f}(\xx_0) - f(\xx_0)|$ is bounded from above by 
    \begin{equation} \label{eq:phase1}
        \bar{z}(\Theta,\xx_0) := \frac{\nu}{2} \sum_{i=0}^{n+1} |\ell_i(\xx_0)| \|\xx_i-\w\|^2, \
        \text{ where }
        \w = \frac{\sum_{i=0}^{n+1} |\ell_i(\xx_0)| \xx_i}{\sum_{i=0}^{n+1} |\ell_i(\xx_0)|}. 
    \end{equation}
\end{corollary}
\begin{proof} 
    Notice the right-hand side of \eqref{eq:phase1 u} is a convex differentiable function of $\u$ defined on $\R^n$ with the first derivative being $\nu \sum_{i=0}^{n+1} |\ell_i(\xx_0)| (\u - \xx_i)$. 
    It attains its minimum when the derivative is the zero vector, which occurs at $\u = \w$. 
    Thus, we can set $\u=\w$ in \eqref{eq:phase1 u} to obtain the improved upper bound \eqref{eq:phase1}. 
\end{proof}

\begin{figure}[ht]
    \centering
    \includegraphics[width=\linewidth]{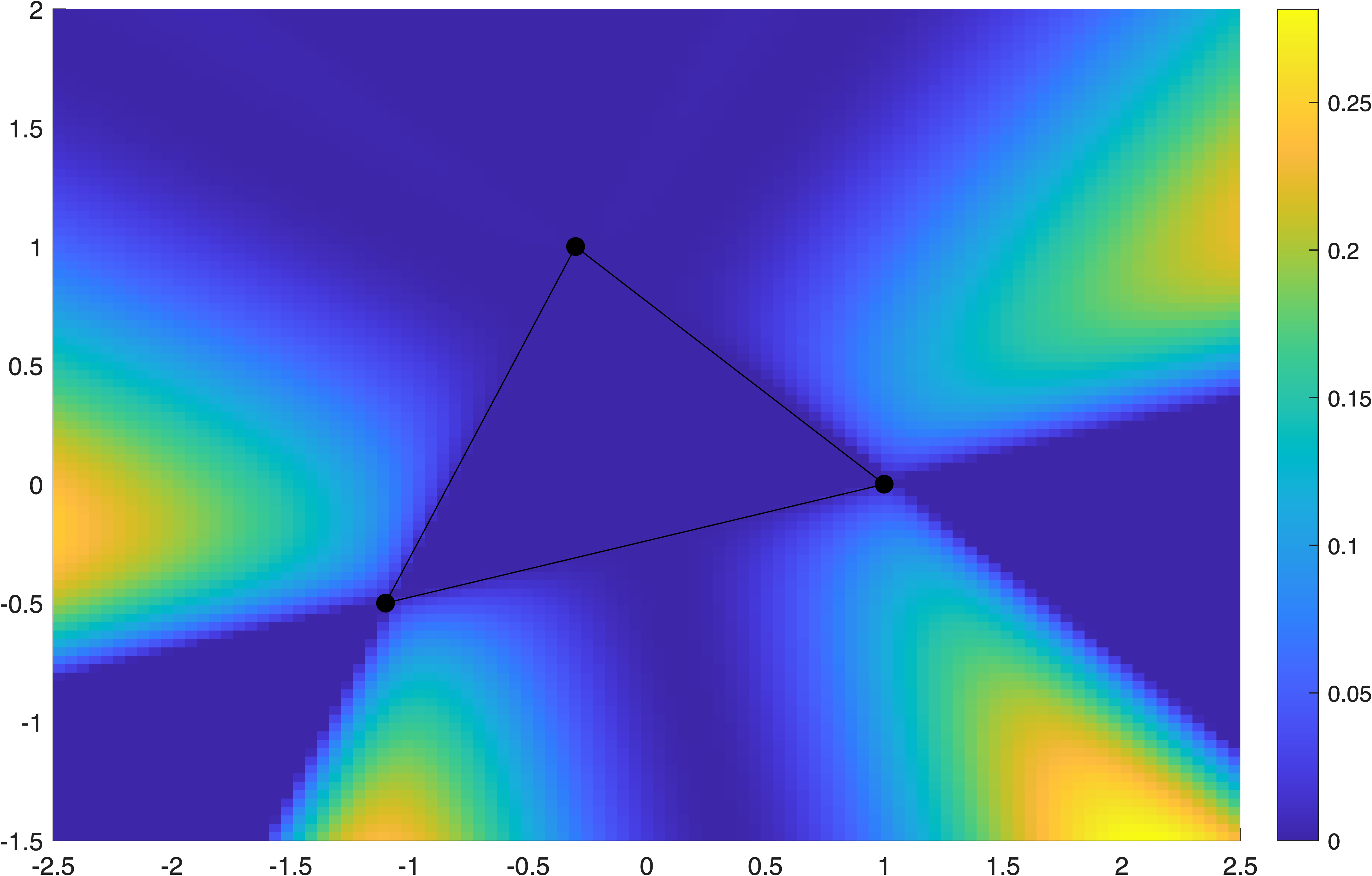}
    \caption{The difference $\bar{z}(\Theta,\xx_0) - z(\Theta,\xx_0)$ as a function of $\xx_0$. The sample set and the Lipschitz constant are chosen as $\Theta = \{(-0.3,1), (-1.1,-0.5), (1,0)\}$ and $\nu = 1$.}
    \label{fig:numerical_sec3}
\end{figure}

Since $z(\Theta,\xx_0)$ is the sharp upper bound on $|\hat{f}(\xx_0) - f(\xx_0)|$, while $\bar{z}(\Theta,\xx_0)$ is only an upper bound, we have $z(\Theta,\xx_0) \le \bar{z}(\Theta,\xx_0)$. 
To check the sharpness of $\bar{z}(\Theta,\xx_0)$, we first compare it against $z(\Theta,\xx_0)$ numerically. 
Similar to how Figure~\ref{fig:numerical} is generated, we fix $\Theta$ to $\{(-0.3,1), (-1.1,-0.5), (1,0)\}$ and evaluate $\bar{z}(\Theta,\xx_0) - z(\Theta,\xx_0)$ for each $\xx_0$ on a $100\times100$ grid that covers the area $[-2.5,2.5]\times[-1.5,2.5]$ evenly. 
The result is plotted in Figure~\ref{fig:numerical_sec3}, where $\bar{z}(\Theta,\xx_0) - z(\Theta,\xx_0)$ for each point $\xx_0$ is indicated by the color of the pixel at $\xx_0$. 
By checking the figure visually, it appears that for this particular example the bound \eqref{eq:phase1} is sharp if $\xx_0$ is located in $\conv(\Theta)$  (the area shown in Figure~\ref{fig:phase1 hull}) or in one of the translated cones 
\begin{equation} \label{eq:cone}
    \left\{\xx_i + \sum_{j=1}^{n+1} \alpha_j(\xx_i - \xx_j):~ \alpha_j \ge 0 \text{ for all } j = 1,2,3 \right\} 
\end{equation}
for $i\in\{1,2,3\}$ (the areas shown in Figure~\ref{fig:phase1 negative}). 
When $\xx_0$ falls in any remaining areas, as shown in Figure~\ref{fig:phase1 not covered}, the numerical results indicate that $z(\Theta,\xx_0) < \bar{z}(\Theta,\xx_0)$. 

Next, we are going to prove $\bar{z}(\Theta,\xx_0)$ is indeed the sharp bound when $\xx_0$ is in areas similar to the ones shown in Figures~\ref{fig:phase1 hull} or \ref{fig:phase1 negative}. 
These areas can be classified mathematically using the signs of the values of the Lagrange functions at $\xx_0$. 
The point $\xx_0 \in \conv(\Theta)$ if and only if $\ell_i(\xx_0) \ge 0$ for all $i=1,\dots,n+1$; and $\xx_0$ is in the translated cones like \eqref{eq:cone} if and only if there is only one positive element in $\{\ell_i(\xx_0)\}_{i=1}^{n+1}$. 

When $f \in C^{1,1}_\nu(\R^n)$, the proof of Theorem 3.1 in \cite{waldron1998error} essentially states that 
\begin{equation} \label{eq:Waldron} 
    |\hat{f}(\xx_0) - f(\xx_0)| \le \frac{\nu}{2} \left(\sum_{i=1}^{n+1} \ell_i(\xx_0) \|\xx_i\|^2 - \|\xx_0\|^2\right)
    = \frac{\nu}{2} \sum_{i=0}^{n+1} \ell_i(\xx_0) \|\xx_i\|^2, 
\end{equation}
holds for all $\xx_0 \in \conv(\Theta)$, and that \eqref{eq:Waldron} a sharp upper bound, as $|\hat{f}(\xx_0) - f(\xx_0)| = \frac{\nu}{2} \sum_{i=0}^{n+1} \ell_i(\xx_0) \|\xx_i\|^2$ when $f$ is the quadratic function $f(\u) = \nu \|\u\|^2/2$. 
We show in Theorem~\ref{thm:phase1 convex hull} that, when $\xx_0 \in \conv(\Theta)$, $\bar{z}(\Theta,\xx_0)$ from \eqref{eq:phase1} is indeed the same as the sharp bound \eqref{eq:Waldron}. 

\begin{theorem} \label{thm:phase1 convex hull}
    Let $\Theta = \{\xx_1,\dots,\xx_{n+1}\}\subset \R^n$ be a fixed set of $n+1$ affinely independent vectors, and $\xx_0$ an arbitrary point inside $\conv(\Theta)$. 
    The vector $\w$ and the bound $\bar{z}(\Theta,\xx_0)$ defined in \eqref{eq:phase1} satisfy $\w=\xx_0$ and $\bar{z}(\Theta,\xx_0) = \frac{\nu}{2} \sum_{i=0}^{n+1} \ell_i(\xx_0) \|\xx_i\|^2$. 
\end{theorem}
\begin{proof}
When $\xx \in \conv(\Theta)$, $\ell_i(\xx_0) \ge 0$ for all $i\in\{1,2,\dots,n+1\}$ and $\ell_0(\xx_0) = -1$. Then, 
\[  \w = \frac{\sum_{i=0}^{n+1} |\ell_i(\xx_0)| \xx_i}{\sum_{i=0}^{n+1} |\ell_i(\xx_0)|}
    = \frac{\xx_0 + \sum_{i=1}^{n+1} \ell_i(\xx_0) \xx_i}{1 +\sum_{i=1}^{n+1} \ell_i(\xx_0)}
    \stackrel{\eqref{eq:Lagrange Y}\eqref{eq:Lagrange 0}}{=} \frac{\xx_0 + \xx_0}{1 + 1}
= \xx_0. 
\]
The following chain of equalities shows the bounds in \eqref{eq:phase1} and \eqref{eq:Waldron} are equal: 
\[ \begin{aligned}
    \bar{z}(\Theta,\xx_0)
    &= \frac{\nu}{2} \sum_{i=0}^{n+1} |\ell_i(\xx_0)| \|\xx_i-\w\|^2 
    = \frac{\nu}{2} \sum_{i=0}^{n+1} \ell_i(\xx_0) \|\xx_i-\xx_0\|^2 \\
    &= \frac{\nu}{2} \sum_{i=0}^{n+1} \ell_i(\xx_0) \left(\|\xx_i\|^2 - 2\xx_i\cdot\xx_0 + \|\xx_0\|^2\right) \\
    &\leftstackrel{\eqref{eq:Lagrange 0}}{=} \frac{\nu}{2} \sum_{i=0}^{n+1} \ell_i(\xx_0) \left(\|\xx_i\|^2 - 2\xx_i\cdot\xx_0 \right) 
    \stackrel{\eqref{eq:Lagrange Y}}{=} \frac{\nu}{2} \sum_{i=0}^{n+1} \ell_i(\xx_0) \|\xx_i\|^2. 
\end{aligned} \]
\end{proof}

Consider the following matrix $G\in\R^{n\times n}$:  
\begin{equation} \label{eq:G}
    G := \sum_{i=0}^{n+1} \ell_i(\xx_0) \xx_i \xx_i^T, 
\end{equation}
which has the property that for any $\u,\v \in \R^n$, 
\begin{equation} \label{eq:G recenter} \begin{aligned}
    \sum_{i=0}^{n+1} \ell_i(\xx_0) [\xx_i-\u] [\xx_i-\v]^T 
    &= \sum_{i=0}^{n+1} \ell_i(\xx_0) \left[\xx_i\xx_i^T - \u\xx_i^T - \xx_i\v^T + \u\v^T \right] \\
    &\leftstackrel{\eqref{eq:Lagrange Y}}{=} \sum_{i=0}^{n+1} \ell_i(\xx_0) \left[\xx_i\xx_i^T + \u\v^T\right] 
    \stackrel{\eqref{eq:Lagrange 0}}{=} \sum_{i=0}^{n+1} \ell_i(\xx_0) \xx_i\xx_i^T = G. 
\end{aligned} \end{equation} 
Recall the set $\Theta = \{\xx_1,\xx_2,\dots,\xx_{n+1}\}$ is ordered in a way such that $\ell_1(\xx_0) \ge \ell_2(\xx_0) \ge \cdots \ge \ell_{n+1}(\xx_0)$. With matrix $G$, we verify in Theorem~\ref{thm:phase1 negative} that
the improved bound $\bar{z}(\Theta,\xx_0)$ is sharp for linear extrapolation when $\xx_0$ is in one of the translated cones indicated by \eqref{eq:cone} and depicted in Figure~\ref{fig:phase1 negative}. 
\begin{theorem} \label{thm:phase1 negative}
    Let $\Theta = \{\xx_1,\dots,\xx_{n+1}\}\subset \R^n$ be a fixed set of $n+1$ affinely independent vectors. 
    Let $\xx_0\in\R^n$ be a point such that $\ell_1(\xx_0) > 0$, and $\ell_i(\xx_0) \le 0$ for all $i\in\{2,\dots,n+1\}$. 
    Then, the vector $\w$ and the bound $\bar{z}(\Theta,\xx_0)$ defined in \eqref{eq:phase1} satisfy $\w=\xx_1$ and $\bar{z}(\Theta,\xx_0) = z(\Theta,\xx_0)$.
\end{theorem}
\begin{proof}
Since $\ell_1(\xx_0) > 0$ and $\ell_i(\xx_0) \le 0$ for all $i > 1$, 
\[  \w = \frac{\sum_{i=0}^{n+1} |\ell_i(\xx_0)| \xx_i}{\sum_{i=0}^{n+1} |\ell_i(\xx_0)|}
    = \frac{2\ell_1(\xx_0)\, \xx_1 - \sum_{i=0}^{n+1} \ell_i(\xx_0)\, \xx_i}{2\ell_1(\xx_0) -\sum_{i=0}^{n+1} \ell_i(\xx_0)} 
    \stackrel{\eqref{eq:Lagrange 0}\eqref{eq:Lagrange Y}}{=} \frac{2\ell_1(\xx_0)\, \xx_1}{2\ell_1(\xx_0)} 
    = \xx_1. 
\] 
The bound $\bar{z}(\Theta,\xx_0)$ equals $\nu/2$ multiplied by 
\[ \begin{aligned}
    \sum_{i=0}^{n+1} &|\ell_i(\xx_0)| \|\xx_i-\w\|^2 
    = - \sum_{i=0}^{n+1} \ell_i(\xx_0) \|\xx_i-\xx_1\|^2 \\
    &= -\text{Tr}\bigg(\sum_{i=0}^{n+1} \ell_i(\xx_0) [\xx_i-\xx_1][\xx_i-\xx_1]^T \bigg) 
    \leftstackrel{\eqref{eq:G recenter}}{=} - \text{Tr}(G)
    =-\sum_{i=0}^{n+1} \ell_i(\xx_0) \|\xx_i\|^2. 
\end{aligned} \]
Now, consider the function $f(\u) = -\frac{\nu}{2} \|\u\|^2 \stackrel{\eqref{eq:Lipschitz Hessian}}{\in} C^{1,1}_\nu(\R^n)$. 
It achieves an approximation error 
\[\hat{f}(\xx_0) - f(\xx_0) 
\stackrel{\eqref{eq:Lagrange m}}{=} \sum_{i=0}^{n+1} \ell_i(\xx_0) f(\xx_i)  
= -\sum_{i=0}^{n+1} \ell_i(\xx_0) \frac{\nu}{2}\|\xx_i\|^2, 
\]
which matches $\bar{z}(\Theta,\xx_0)$, meaning $\bar{z}(\Theta,\xx_0) = z(\Theta,\xx_0)$. 
\end{proof}

\begin{figure}[tbhp]
     \centering
     \subfloat[\raggedright The convex hull covered by Theorem~\ref{thm:phase1 convex hull}]{\label{fig:phase1 hull}\resizebox{0.3\linewidth}{!}{\begin{tikzpicture}
\filldraw[black] (1,0) circle (2pt) node[anchor=north] {$\xx_1$}; 
\filldraw[black] (-0.3,1) circle (2pt) node[anchor=east] {$\xx_2$}; 
\filldraw[black] (-1.1,-0.5) circle (2pt) node[anchor=south east] {$\xx_3$}; 

\fill[color=gray, opacity=0.7] (1,0) -- (-0.3,1) -- (-1.1,-0.5);

\draw[thick] (-1.99,2.3) -- (3,-1.5385);
\draw[thick] (0.5,2.5) -- (-1.7933,-1.8);
\draw[thick] (-3.5,-1.0714) -- (3,0.4762);

\end{tikzpicture}}}
     \hfill
     \subfloat[\raggedright The cones covered by Theorem~\ref{thm:phase1 negative}]{\label{fig:phase1 negative}\resizebox{0.3\linewidth}{!}{\begin{tikzpicture}
\filldraw[black] (1,0) circle (2pt) node[anchor=north] {$\xx_1$}; 
\filldraw[black] (-0.3,1) circle (2pt) node[anchor=east] {$\xx_2$}; 
\filldraw[black] (-1.1,-0.5) circle (2pt) node[anchor=south east] {$\xx_3$}; 

\shade[shading angle=80] (3,0.4762) -- (1,0) -- (3,-1.5385);
\shade[shading angle=-175] (-1.99,2.3) -- (-0.3,1) -- (0.5,2.5);
\shade[shading angle=-40] (-3.5,-1.0714) -- (-1.1,-0.5) -- (-1.7933,-1.8);

\draw[thick] (-1.99,2.3) -- (3,-1.5385);
\draw[thick] (0.5,2.5) -- (-1.7933,-1.8);
\draw[thick] (-3.5,-1.0714) -- (3,0.4762);

\end{tikzpicture}}}
     \hfill
     \subfloat[The areas where \eqref{eq:phase1} holds but is not sharp]{\label{fig:phase1 not covered} \resizebox{0.3\linewidth}{!}{\begin{tikzpicture}
\shade[shading angle=135] (0.5,2.5) -- (-0.3,1) -- (1,0) -- (3,0.4762);
\shade[shading angle=-115] (-1.99,2.3) -- (-0.3,1) -- (-1.1,-0.5) -- (-3.5,-1.0714);
\shade[shading angle=7] (-1.7933,-1.8) -- (-1.1,-0.5) -- (1,0) -- (3,-1.5385);

\filldraw[black] (1,0) circle (2pt) node[anchor=north] {$\xx_1$}; 
\filldraw[black] (-0.3,1) circle (2pt) node[anchor=east] {$\xx_2$}; 
\filldraw[black] (-1.1,-0.5) circle (2pt) node[anchor=south east] {$\xx_3$}; 

\draw[thick] (-1.99,2.3) -- (3,-1.5385);
\draw[thick] (0.5,2.5) -- (-1.7933,-1.8);
\draw[thick] (-3.5,-1.0714) -- (3,0.4762);

\end{tikzpicture}}}
    \caption{A visualization of results in Section~\ref{sec:phase1} for bivariate interpolation. The ordering of the points in $\Theta = \{\xx_1,\xx_2,\xx_3\}$ in this figure and all figures hereafter is arbitrary and not determined by the values of the Lagrange polynomials at $\xx_0$. }
    \label{fig:phase1}
\end{figure}

\section{The Worst Quadratic Function} \label{sec:phase2}
We have derived an improved error bound in the previous section and showed when it is sharp. 
In this section, we try to find the mathematical formula for the piecewise smooth function in the remaining areas indicated in Figure~\ref{fig:phase1 not covered}. 
Instead of attempting to improve another existing upper bound, we take the opposite approach by trying to find the function that can achieve the maximum error. 
Considering quadratic functions are easier to analyze as they share a general closed-form formula and, under the settings of both Theorem~\ref{thm:phase1 convex hull} and Theorem~\ref{thm:phase1 negative}, the optimal set of \eqref{prob:D} contains at least one quadratic function, we investigate whether \eqref{prob:D} has an analytical solution when $f$ is restricted to be quadratic. 

Let $f$ be a quadratic function of the form $f(\u) = c + \g \cdot \u + \frac{1}{2} H\u \cdot \u$ with $c\in\R, \g\in\R^n$, and symmetric $H \in \R^{n \times n}$. 
Recall the matrix $G$ defined in \eqref{eq:G}. 
Because of \eqref{eq:Lipschitz Hessian} and  
\[ \begin{aligned}
\hat{f}(\xx_0) - f(\xx_0) ~
&\leftstackrel{\eqref{eq:Lagrange m}}{=} \sum_{i=0}^{n+1} \ell_i(\xx_0) f(\xx_i) 
= \sum_{i=0}^{n+1} \ell_i(\xx_0) \left[c + \g \cdot \xx_i + \frac{1}{2} H\xx_i \cdot \xx_i \right] \\
&\leftstackrel{\eqref{eq:Lagrange Y}}{=} \sum_{i=0}^{n+1} \ell_i(\xx_0) \left[c + \frac{1}{2} H\xx_i \cdot \xx_i \right] 
\stackrel{\eqref{eq:Lagrange 0}}{=} \sum_{i=0}^{n+1} \ell_i(\xx_0) \left[\frac{1}{2} H\xx_i \cdot \xx_i \right] \\
&= \frac{1}{2} H \cdot \sum_{i=0}^{n+1} \ell_i(\xx_0) \xx_i \xx_i^T
\stackrel{\eqref{eq:G}}{=} \frac{1}{2}  G \cdot H,  
\end{aligned} \]
the problem of maximizing linear interpolation's approximation error over quadratic functions in $C_\nu^{1,1}(\R^n)$ can be formulated as 
\begin{equation} \label{prob:quadratic} \everymath{\displaystyle} \begin{array}{ll}
    \max_H &G \cdot H / 2  \\
    \text{s.t.} &-\nu I \preceq H \preceq \nu I. 
\end{array} \end{equation}
The absolute sign in the objective function is again dropped due to symmetry. 

It turns out problem \eqref{prob:quadratic} can be solved analytically. 
Since $G$ is real and symmetric, it admits an eigendecomposition $G = P \Lambda P^T$, where $\Lambda \in \R^{n \times n}$ is the diagonal matrix of the eigenvalues $\lambda_1,\dots,\lambda_n$, and $P \in \R^{n \times n}$ is the orthogonal matrix (i.e. $P^{-1} = P^T$) whose columns are the corresponding eigenvectors. 
The objective function is $G\cdot H/2 = (P\Lambda P^T)\cdot H/2 = \Lambda \cdot (P^T H P)/2$. 
Since $P$ is orthogonal, the constraint in \eqref{prob:quadratic} is equivalent to $-\nu I \preceq P^T H P \preceq \nu I$, indicating all diagonal elements of $P^T H P$ are bounded between $-\nu$ and $\nu$. 
Since $\Lambda$ is diagonal, only the diagonal elements of $P^T H P$ would affect the objective function value. 
Therefore, a solution to \eqref{prob:quadratic}, denoted by $H^\star$, has the property $P^T H^\star P = \nu\, \text{sign}(\Lambda)$, where the sign function $\text{sign}(\cdot)$ acts entry-wise on a matrix, returning a matrix of the same dimension that preserves zero entries while replacing positive entries with 1 and negative entries with -1. 
This optimal solution is 
\begin{equation} \label{eq:Hstar}
    H^\star = \nu P \text{sign}(\Lambda) P^T. 
\end{equation} 
Solution \eqref{eq:Hstar} indicates the maximum approximation error achievable by the quadratic functions in $C_\nu^{1,1}(\R^n)$ is 
\begin{equation} \label{eq:phase2}
    \underline{z}(\Theta,\xx_0) :=~ G \cdot H^\star/2 = \frac{\nu}{2} \sum_{i=1}^n |\lambda_i|. 
\end{equation} 

Since the error $\underline{z}(\Theta,\xx_0)$ can be no more than the maximum error achievable by all functions in $C_\nu^{1,1}(\R^n)$, we have $\underline{z}(\Theta,\xx_0) \le z(\Theta,\xx_0)$. 
To understand how close these two values are, we fix $\Theta$ to $\{(-0.3,1), (-1.1,-0.5), (1,0)\}$ again (as in Figures~\ref{fig:numerical} and \ref{fig:numerical_sec3}) and compute the gap $z(\Theta,\xx_0) - \underline{z}(\Theta,\xx_0)$ numerically. 
We do not provide a picture here to show the result, because it would simply be a function of near zero values --- it turns out that, for this particular example, no matter where $\xx_0 \in \R^2$ is, the computed difference between $\underline{z}(\Theta,\xx_0)$ and $z(\Theta,\xx_0)$ is very small ($< 10^{-7}$). 
Considering optimization solvers typically can only find an approximate solution, it is reasonable to suspect these tiny differences come from solving \eqref{prob:f-D} numerically, and there is no actual difference between $\underline{z}(\Theta,\xx_0)$ and $z(\Theta,\xx_0)$. 
This suspicion turns out to be true, and we will prove $\underline{z}(\Theta,\xx_0)=z(\Theta,\xx_0)$ under a more general condition on $\Theta$ and $\xx_0$ later. 

\begin{figure}[ht]
    \centering
    \includegraphics[width=\linewidth]{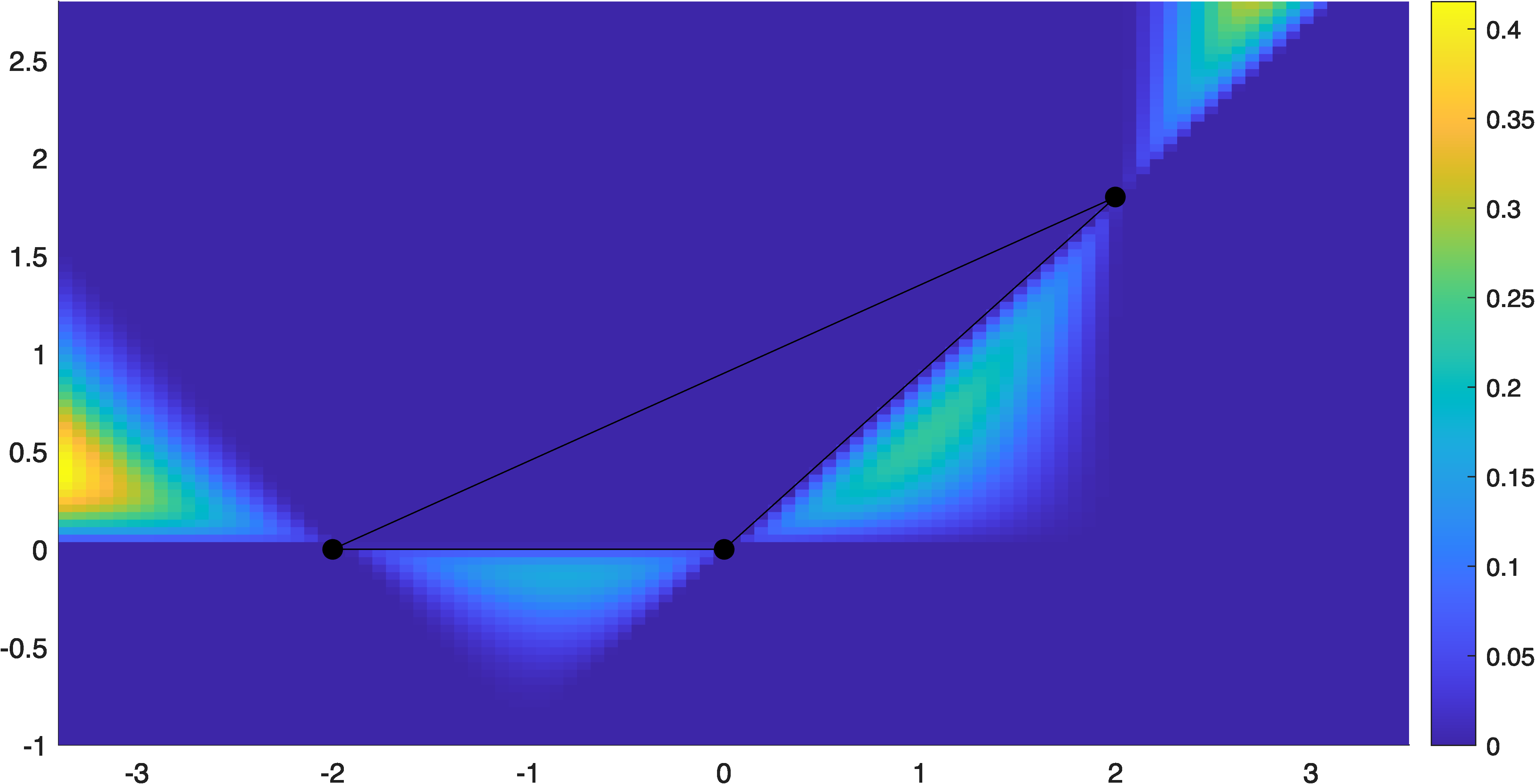}
    \caption{The difference $z(\Theta,\xx_0) - \underline{z}(\Theta,\xx_0)$ as a function of $\xx_0$. The sample set and the Lipschitz constant are chosen as $\Theta = \{(0,0), (2,1.8), (-2,0)\}$ and $\nu = 1$.}
    \label{fig:numerical_phase2_obtuse}
\end{figure}

Unfortunately, this does not mean \eqref{eq:phase2} is the closed-form solution of \eqref{prob:D}. 
Change $\Theta$ to $\{(0,0), (2,1.8), (-2,0)\}$ and then solve \eqref{prob:f-D} repeatedly to obtain $z(\Theta,\xx_0) - \underline{z}(\Theta,\xx_0)$ for each point $\xx_0$ on a $100\times100$ grid that covers the area $[-3.4,3.3]\times[-1,2.8]$ evenly. 
The result is shown in Figure~\ref{fig:numerical_phase2_obtuse}, where the value of the difference $z(\Theta,\xx_0) - \underline{z}(\Theta,\xx_0)$ at each point $\xx_0$ is indicated by the color of the pixel at $\xx_0$. 
It can be observed that $\underline{z}(\Theta,\xx_0) < z(\Theta,\xx_0)$ when $\xx_0$ is in some areas of the input space. 
These areas are made clear in Figure~\ref{fig:phase2}. 
They are the two translated cones \circnum{1} \circnum{2} and the two triangles \circnum{3} \circnum{4}. 
Mathematically, they can be described respectively as $\xx_0\in\R^2$ such that 
\begin{itemize}
    \item[\circnum{1}] $\ell_1(\xx_0) [\xx_2-\xx_1]\cdot[\xx_3-\xx_1] - \ell_2(\xx_0) [\xx_3-\xx_2]\cdot[\xx_1-\xx_2] > 0$ and $\ell_2(\xx_0)>0$;
    \item[\circnum{2}]$\ell_1(\xx_0) [\xx_2-\xx_1]\cdot[\xx_3-\xx_1] - \ell_3(\xx_0) [\xx_2-\xx_3]\cdot[\xx_1-\xx_3] > 0$ and $\ell_3(\xx_0)>0$; 
    \item[\circnum{3}] $\ell_1(\xx_0) [\xx_2-\xx_1]\cdot[\xx_3-\xx_1] - \ell_2(\xx_0) [\xx_3-\xx_2]\cdot[\xx_1-\xx_2] < 0$, $\ell_3(\xx_0)>0$, and $\ell_2(\xx_0)<0$;
    \item[\circnum{4}] $\ell_1(\xx_0) [\xx_2-\xx_1]\cdot[\xx_3-\xx_1] - \ell_3(\xx_0)[\xx_2-\xx_3]\cdot[\xx_1-\xx_3] < 0$, $\ell_2(\xx_0)>0$, and $\ell_3(\xx_0)<0$. 
\end{itemize}
We will look for the expression for $z(\Theta,\xx_0)$ when $\xx_0$ is inside these areas in the next section. 
In the remaining parts of this section, we derive a sufficient condition for $\underline{z}(\Theta,\xx_0) = z(\Theta,\xx_0)$. 

\begin{figure}[tbhp]
  \centering
    \resizebox{0.4\textwidth}{!}{\begin{tikzpicture}

\filldraw[black] (0,0) circle (2pt) node[anchor=north west] {$\xx_1$}; 
\filldraw[black] (2,1.8) circle (2pt) node[anchor=south east] {$\xx_2$}; 
\filldraw[black] (-2,0) circle (2pt) node[anchor=south] {$\xx_3$}; 

\fill[color=gray, opacity=0.7] (0,0) -- (2,1.8) -- (2,0); 
\fill[color=gray, opacity=0.7] (0,0) -- (-2,0) -- (-1.1050, -0.9945); 
\shade[left color=gray, right color=white, shading angle=156] (2,3.5) -- (2,1.8) -- (3.333,3);
\shade[left color=gray, right color=white, shading angle=-120] (-3.7,0) -- (-2,0) -- (-3.35,1.5);

\draw[thick] (-2,-1.8) -- (3.333,3);
\draw[thick] (-3,-0.45) -- (4,2.7);
\draw[thick] (-3.7,0) -- (3.333,0);
\draw[thick, dashed] (2,-0.5) -- (2,3.5); 
\draw[thick, dashed] (-3.35,1.5) -- (-0.2,-2); 

\node[circle,draw,inner sep=1pt] at (-3,0.5) {1};
\node[circle,draw,inner sep=1pt] at (2.5,2.8) {2};
\node[circle,draw,inner sep=1pt] at (-1.1,-0.4) {3};
\node[circle,draw,inner sep=1pt] at (1.3,0.5) {4};

\end{tikzpicture}}
  \caption{The areas to which if $\xx_0$ belongs, $\underline{z}(\Theta,\xx_0)$ from \eqref{eq:phase2} is not an upper bound on the function approximation error for bivariate interpolation. 
  The dashed line on the left is perpendicular to the line going through $\xx_1$ and $\xx_2$; and the one on the right is perpendicular to the line going through $\xx_3$ and $\xx_1$. } 
  \label{fig:phase2}
\end{figure}

\subsection{Certification of Upper Bound}
The maximum error \eqref{eq:phase2} provides a lower bound to the optimal value of \eqref{prob:D}, while \eqref{eq:phase1} provides an upper bound. 
By evaluating both \eqref{eq:phase1} and \eqref{eq:phase2}, one can have a reasonable estimation of sharp error bound without having to solve the QCQP \eqref{prob:f-D}. 
However, formula \eqref{eq:phase2} would be a lot more useful if there is an efficient way to check whether $\xx_0$ is in one of those areas where \eqref{eq:phase2} is not an upper bound on the approximation error. 

The existence of these areas appears to be influenced by the existence of obtuse angles at the vertices of the simplex $\conv(\Theta)$. 
Unlike triangles, which can only have up to one obtuse angle, simplices in higher dimension can have obtuse angles in many ways. 
They can have $(\xx_j-\xx_i)\cdot(\xx_k-\xx_i) < 0$ at multiple vertices $\xx_i$ and at the same time for multiple $(j,k)$ for each $\xx_i$. 
While there can only be up to four disconnected subset of $\R^2$ where \eqref{eq:phase2} is not an upper bound on the approximation error, our numerical experiments show this number can go up to at least twenty for trivariate ($n=3$) linear interpolation. 
Considering a precise description of the four shaded areas in Figure~\ref{fig:phase2} already requires four unintuitive inequalities or some wordy explanation, any description of these areas would almost certainly be extremely complicated, especially in higher dimension. 

Regardless, we have found an efficient way to check whether $\xx_0$ is in one of these areas without having to describe any of them. 
The theoretical proof that validates our approach is very technical and will be presented later in Section~\ref{sec:phase2 proofs}. 
Our approach relies on a set of parameters $\{\mu_{ij}\}_{(i,j)\in\cI_+\times\cI_-}$ that can be computed as follows. 
Recall $\Theta$ is ordered such that $\ell_1(\xx_0) \ge \ell_2(\xx_0) \ge \cdots \ge \ell_{n+1}(\xx_0)$. Let $\diag(\ell) \in \R^{(n+1)\times(n+1)}$ be the diagonal matrix containing $\ell_1(\xx_0), \dots, \ell_{n+1}(\xx_0)$. 
We now partition $\diag(\ell), G$, and $H^\star$ with respect to $\cI_+$ and $\cI_-$. 
Let $\diag(\ell_+)\in\R^{|\cI_+|\times|\cI_+|}$ be the diagonal matrix containing $\{\ell_i(\xx_0)\}_{i\in\cI_+}$, and $\diag(\ell_-) \in \R^{(|\cI_-|-1)\times(|\cI_-|-1)}$ be the diagonal matrix containing $\{\ell_i(\xx_0)\}_{i\in\cI_-\setminus\{0\}}$. 
Let $Y_+ \in \R^{|\cI_+| \times n}$ and $Y_- \in \R^{(|\cI_-|-1) \times n}$ be the first $|\cI_+|$ and the last $|\cI_-|-1$ rows of the matrix $Y$ defined in \eqref{eq:Y Phi} respectively. 
The matrix $G$ has $|\cI_+|-1$ positive eigenvalues and $|\cI_-|-1$ negative eigenvalues, as will be proved later. 
Let $\Lambda_+ \in \R^{(|\cI_+|-1) \times (|\cI_+|-1)}$ and $\Lambda_- \in \R^{(|\cI_-|-1) \times (|\cI_-|-1)}$ respectively be the submatrices of $\Lambda$ that contain the positive and negative eigenvalues of $G$, and $P_+ \in \R^{n \times (|\cI_+|-1)}$ and $P_- \in \R^{n \times (|\cI_-|-1)}$ the submatrices of $P$ with their corresponding eigenvectors. 
Then we have 
\begin{equation} \label{eq:G+-} \begin{aligned} 
G~ &\leftstackrel{\eqref{eq:G recenter}}{=} Y^T \diag(\ell) Y  = Y_+^T \diag(\ell_+) Y_+ + Y_-^T \diag(\ell_-) Y_- \\
&= P \Lambda P^T
= P_+ \Lambda_+ P_+^T + P_- \Lambda_- P_-^T
\end{aligned} \end{equation}
and 
\begin{equation} \label{eq:Hstar+-}
    H^\star = \nu P \text{sign}(\Lambda) P^T
    = \nu (P_+ P_+^T - P_- P_-^T). 
\end{equation} 

\begin{remark}
    When $|\cI_+| = 1$, the dimensions of $\Lambda_+$ and $P_+$ are $0 \times 0$ and $n \times 0$, respectively. 
    In this case, we follow the convention that both $P_+ \Lambda_+ P_+^T$ and $P_+ P_+^T$ are $\mathbf{0}_{n\times n}$, the $n\times n$ zero matrix. 
    Similarly, when $|\cI_-| = 1$, the matrices $Y_-^T \diag(\ell_-) Y_-$, $P_- \Lambda_- P_-^T$, and $P_- P_-^T$ are all equal to $\mathbf{0}_{n\times n}$. 
    It is possible that $|\cI_+| = |\cI_-| = 1$, in which case the only two Lagrange polynomials with nonzero values at $\xx_0$ are $\ell_0(\xx_0) = -1$ and $\ell_i(\xx_0) = 1$ for some $i \in\{1,\dots,n+1\}$. 
    However, this implies that $\xx_0 = \xx_i$, and the approximation error is simply $\hat{f}(\xx_0) - f(\xx_0) = \hat{f}(\xx_i) - f(\xx_i) = 0$. 
\end{remark}

We now present the definition of $\{\mu_{ij}\}_{(i,j)\in\cI_+\times\cI_-}$ and the main theorem of this section. 

\begin{theorem} \label{thm:phase2} 
    Assume $f \in C^{1,1}_\nu(\R^n)$. 
    Let $\xx_0$ be a fixed vector in $\R^n$, and $\Theta = \{\xx_1,\dots,\xx_{n+1}\}\subset \R^n$ be a fixed set of $n+1$ affinely independent vectors. 
    Consider the matrix $M \stackrel{\rm def}{=} \diag(\ell_+) Y_+ P_- (Y_- P_-)^{-1}$. 
    Let $\mu_{ij} = [M]_{i, j-n-2 + |\cI_-|}$ for all $i\in\cI_+$ and $j\in\cI_-\setminus\{0\}$, and $\mu_{i0} = \ell_i(\xx_0) - \sum_{j \in \cI_-\setminus\{0\}} \mu_{ij}$ for all $i\in\cI_+$. 
    The maximum error $\underline{z}(\Theta,\xx_0)$ defined in \eqref{eq:phase2} equals to the sharp error bound $z(\Theta,\xx_0)$ if $\mu_{ij} \ge 0$ for all $(i,j) \in \cI_+\times\cI_-$. 
\end{theorem}

\begin{remark}
    We note that $\{j-n-2 + |\cI_-|\}_{j\in\cI_-\setminus\{0\}} =$ $\{j - (n-|\cI_-|+2)\}_{j\in\cI_-\setminus\{0\}} =$ $\{1,2,\dots,|\cI_-|-1\}$. 
    The matrix $M$ is of size $|\cI_+| \times (|\cI_-|-1)$. 
    Each of its rows corresponds to an interpolation point with positive Lagrange polynomial value at $\xx_0$, while each column corresponds to an interpolation point with negative Lagrange polynomial value at $\xx_0$. 
    When $\cI_- = \{0\}$, the size of $M$ becomes $|\cI_+| \times 0$, and we simply have $\mu_{i0} = \ell_i(\xx_0)$ for all $i\in\cI_+$. 
\end{remark}

\subsection{Technical Proofs} \label{sec:phase2 proofs}
In the remaining of this section, we provide the complete proof to Theorem~\ref{thm:phase2}. 
We start with the numbers of positive and negative eigenvalues in the matrix $G$ defined in \eqref{eq:G}. 

\begin{lemma} \label{lem:Sylvester}
    The numbers of positive and negative eigenvalues in $G$ are $|\cI_+| - 1$ and $|\cI_-|-1$, respectively. If 0 is an eigenvalue of $G$, then its multiplicity is $n+2 - |\cI_+| - |\cI_-|$. 
\end{lemma}
\begin{proof} 
    Recall our definitions of $\Phi$ in Equation~\eqref{eq:Y Phi}. 
    Let $\diag(\ell)$ be the diagonal matrix containg $\ell_1(\xx_0), \dots, \ell_{n+1}(\xx_0)$. 
    Consider the matrix 
    \[\bar{G} = \sum_{i=1}^{n+1} \ell_i(\xx_0) \begin{bmatrix} 1\\ \xx_i-\xx_0 \end{bmatrix} \begin{bmatrix} 1 &(\xx_i-\xx_0)^T \end{bmatrix} = \Phi^T \diag(\ell) \Phi.
    \] 
    The first element of the first column is $\sum_{i=1}^{n+1} \ell_i(\xx_0) \stackrel{\eqref{eq:Lagrange 0}}{=} 1$, while the rest of the column is $\sum_{i=1}^{n+1} \ell_i(\xx_0) [\xx_i-\xx_0] \stackrel{\eqref{eq:Lagrange Y}}{=} \xx_0 - \sum_{i=1}^{n+1} \ell_i(\xx_0) \xx_i \stackrel{\eqref{eq:Lagrange 0}}{=} \mathbf{0}$. 
    The bottom-right $n\times n$ submatrix of $\bar{G}$ is  
    \[  \sum_{i=1}^{n+1} \ell_i(\xx_0) [\xx_i-\xx_0][\xx_i-\xx_0]^T
        = \sum_{i=0}^{n+1} \ell_i(\xx_0) [\xx_i-\xx_0][\xx_i-\xx_0]^T
        \stackrel{\eqref{eq:G recenter}}{=} G. 
    \]
    Thus, $\bar{G}$ and its eigendecomposition should be 
    \[ \bar{G} = \begin{bmatrix} 1 &\mathbf{0}^T\\ \mathbf{0} &G \end{bmatrix}
    = \begin{bmatrix} 1 &\mathbf{0}^T\\ \mathbf{0} &P \end{bmatrix}
    \begin{bmatrix} 1 &\mathbf{0}^T\\ \mathbf{0} &\Lambda \end{bmatrix}
    \begin{bmatrix} 1 &\mathbf{0}^T\\ \mathbf{0} &P^T \end{bmatrix}. 
    \]
    Then we have 
    \[ \bar\Lambda \stackrel{\rm def}{=} \begin{bmatrix} 1 &\mathbf{0}^T\\ \mathbf{0} &\Lambda \end{bmatrix}
    = \begin{bmatrix} 1 &\mathbf{0}^T\\ \mathbf{0} &P^T \end{bmatrix} 
    \Phi^T \diag(\ell) \Phi
    \begin{bmatrix} 1 &\mathbf{0}^T\\ \mathbf{0} &P \end{bmatrix}, 
    \]
    which shows $\bar\Lambda$ is congruent to $\diag(\ell)$. 
    Then by Sylvester's law of inertia \cite{sylvester1852xix} (or Theorem 4.5.8 of \cite{horn2012matrix}), the number of positive and negative eigenvalues in $\bar\Lambda$ are $|\cI_+|$ and $|\cI_-|-1$, respectively. 
    Since $\bar{G}$ shares the same eigenvalues as $G$ except an additional one that is 1, the lemma is proved. 
\end{proof}

The next lemma shows that $\{\mu_{ij}\}_{(i,j)\in\cI_+\times\cI_-}$ is well-defined by proving the invertibility of $Y_- P_-$. 
\begin{lemma} \label{lem:invertible}
    The $(|\cI_-|-1)$-dimensional square matrix $Y_- P_-$ is invertible. 
\end{lemma}
\begin{proof}
    For the purpose of contradiction, assume $Y_- P_-$ is singular. 
    That means there is a non-zero vector $\u \in \R^{|\cI_-|-1}$ such that $Y_- P_- \u = \mathbf{0}$. 
    Let $\v = P_- \u$. 
    We have $Y_- \v = \mathbf{0}$, $P_+^T \v = P_+^T P_- \u = \mathbf{0}$ and $P_-^T \v = P_-^T P_- \u = \u$. 
    Then we have the contradiction
    \[ \begin{aligned} 
    \v^T G \v &= (Y_+\v)^T \diag(\ell_+) Y_+\v + (Y_-\v)^T \diag(\ell_-) Y_-\v = (Y_+\v)^T \diag(\ell_+) Y_+\v \ge 0, \\
    \v^T G \v &= (P_+^T \v)^T \Lambda_+ P_+^T \v + (P_-^T \v)^T \Lambda_- P_-^T \v = (P_-^T \v)^T \Lambda_- P_-^T \v  = \u^T \Lambda_- \u < 0.  
    \end{aligned} \]
\end{proof}

Let $P_0 \in \R^{n \times (n+2 - |\cI_+| - |\cI_-|)}$ be the submatrix of $P$ whose columns are the eigenvectors of $G$ that correspond to the eigenvalue 0. 
The following lemma shows 
\begin{subequations} \label{eq:P0} \begin{align} 
    Y_+P_0 &= \mathbf{0} \in \R^{|\cI_+| \times (n+2 - |\cI_+| - |\cI_-|)}  \\
    \text{ and } Y_-P_0 &= \mathbf{0} \in \R^{(|\cI_-|-1) \times (n+2 - |\cI_+| - |\cI_-|)}, 
\end{align} \end{subequations}
which will be used in the proof of Lemma~\ref{lem:mu}. 
\begin{lemma}
    For $i=1,\dots,n+1$, if $\ell_i(\xx_0) \neq 0$, the vector $\xx_i - \xx_0$ is perpendicular to the null space of $G$ and $P_0^T(\xx_i - \xx_0) = \mathbf{0}$. 
\end{lemma}
\begin{proof}
    Notice the columns of $P_0$ form a basis of the null space of $G$. 
    Let $\u$ be any vector in this null space. 
    Set $\v = \diag(\ell) Y \u$. 
    Then, we have 
    \[ \begin{aligned} 
        \mathbf{1}^T \v = \mathbf{1}^T \diag(\ell) Y \u = \ell^T Y \u \stackrel{\eqref{eq:Lagrange Y}}= 0 
        \text{ and } Y^T \v = Y^T \diag(\ell) Y \u =  G\u = \mathbf{0}.
    \end{aligned} \] 
    Recall $\Phi = \begin{bmatrix} \mathbf{1} &\bY \end{bmatrix}$. 
    The above two equalities combine to give $\Phi^T \v = \mathbf{0}$. 
    Then, the invertibility of $\Phi$ implies $\v = \mathbf{0}$ or, equivalently, $\ell_i(\xx_0) (\xx_i-\xx_0)^T \u = 0$, $i=1,\dots,n+1$. 
    If $\ell_i(\xx_0) \neq 0$, we must have $(\xx_i-\xx_0)^T \u = 0$. 
\end{proof}

We develop in the following lemma the essential properties of $\{\mu_{ij}\}$. 
\begin{lemma} \label{lem:mu}
The following properties hold for the $\{\mu_{ij}\}_{(i,j)\in\cI_+\times\cI_-}$ defined in Theorem~\ref{thm:phase2}: 
\begin{align}
    \sum_{j \in \cI_-} \mu_{ij} &= \ell_i(\xx_0) &&\text{for all } i \in \cI_+, \label{eq:mu0 +}\\
    \sum_{i \in \cI_+} \mu_{ij} &= - \ell_j(\xx_0) &&\text{for all } j \in \cI_-, \label{eq:mu0 -}\\
    (\nu I-H^\star) \sum_{j \in \cI_-} \mu_{ij} \xx_j &= (\nu I-H^\star) \ell_i(\xx_0) \xx_i &&\text{for all } i \in \cI_+, \label{eq:mu1 +}\\
    (\nu I+H^\star) \sum_{i \in \cI_+} \mu_{ij}\xx_i &= -(\nu I+H^\star) \ell_j(\xx_0) \xx_j &&\text{for all } j \in \cI_-,  \label{eq:mu1 -}
\end{align}
where $H^\star$ is the optimal solution \eqref{eq:Hstar}. 
\end{lemma}
\begin{proof}
The equations in \eqref{eq:mu0 +} hold by the definition of $\mu_{ij}$ in Theorem~\ref{thm:phase2}. 

Let $\ell_+ = \big(\ell_1(\xx_0), \dots, \ell_{|\cI_+|}(\xx_0)\big)^T$ and $\ell_- =  \big(\ell_{n-|\cI_-|+3}(\xx_0), \dots, \ell_{n+1}(\xx_0)\big)^T$. 
Notice $\sum_{i \in \cI_+} \mu_{ij} = \e_j^T M^T \mathbf{1}$ for $j\in\mathcal{I}_-\backslash\{0\}$. 
The equations in \eqref{eq:mu0 -} for $j\in\mathcal{I}_-\backslash\{0\}$ can thus be written as $M^T\mathbf{1} = - \ell_-$. 
They hold true because 
\[ \begin{aligned}
    \ell_- + M^T \mathbf{1}
    = \ell_- + (P_-^T Y_-^T)^{-1} P_-^T Y_+^T \diag(\ell_+) \mathbf{1} 
    = (P_-^T Y_-^T)^{-1} P_-^T [Y_-^T \ell_- + Y_+^T \ell_+] 
    \stackrel{\eqref{eq:Lagrange Y}}{=} \mathbf{0}. 
\end{aligned} \] 
For $j = 0$, equality \eqref{eq:mu0 -} holds because 
\[ \begin{aligned}
    \sum_{i \in \cI_+} \mu_{i0} 
    &= \sum_{i \in \cI_+} \bigg( \ell_i(\xx_0) - \sum_{j \in \cI_-\setminus\{0\}} \mu_{ij} \bigg)
    = \sum_{i \in \cI_+} \ell_i(\xx_0) - \sum_{j \in \cI_-\setminus\{0\}} \sum_{i \in \cI_+} \mu_{ij} \\
    &\leftstackrel{\eqref{eq:mu0 -}}{=} \sum_{i \in \cI_+} \ell_i(\xx_0) + \sum_{j \in \cI_-\setminus\{0\}} \ell_j(\xx_0)
    = 1 = - \ell_0(\xx_0). 
\end{aligned} \] 

Notice $P_-^T (Y_-^T M^T - Y_+^T \diag(\ell_+)) = \mathbf{0}$ by the definition of $M$ in Theorem~\ref{thm:phase2} and $\nu I - H^\star \stackrel{\eqref{eq:Hstar+-}}{=} \nu (P_+P_+^T + P_0P_0^T + P_-P_-^T) - \nu (P_+P_+^T - P_-P_-^T) = \nu P_0P_0^T + 2\nu P_-P_-^T$. 
Following these two equations, we have for all $i \in \cI_+$, 
\[ \begin{aligned}
    (\nu I-H^\star) \left[\sum_{j \in \cI_-} \mu_{ij} \xx_j - \ell_i(\xx_0) \xx_i\right] 
    &\leftstackrel{\eqref{eq:mu0 +}}{=} (\nu I-H^\star) \left[\sum_{j \in \cI_-} \mu_{ij} (\xx_j-\xx_0) - \ell_i(\xx_0) [\xx_i-\xx_0]\right] \\
    &= (\nu I-H^\star) \left[Y_-^T M^T\e_i - Y_+^T \diag(\ell_+)\e_i\right]  \\
    &{\stackrel{\eqref{eq:P0}}{=}} 2\nu P_- P_-^T (Y_-^T M^T - Y_+^T \diag(\ell_+)) \e_i \\
    &= 2\nu P_- \mathbf{0} \e_i = \mathbf{0}, 
\end{aligned} \] 
which proves \eqref{eq:mu1 +}. 

To prove \eqref{eq:mu1 -}, we use $G$ and its eigendecomposition. 
Let $Y_0 \in \R^{(n+2-|\cI_+|-|\cI_-|) \times n}$ be the submatrix of $Y$ that contains rows $(\xx_i-\xx_0)^T$ with $\ell_i(\xx_0) = 0$. 
The eigenvalue matrix $\Lambda$ follows $\Lambda = P^T Y^T \diag(\ell) Y P$. 
This equation can be written as 
\[ \begin{aligned}  
    \Lambda &= \begin{bmatrix} \Lambda_- &\mathbf{0} &\mathbf{0}\\ \mathbf{0} &\mathbf{0} &\mathbf{0}\\ \mathbf{0} &\mathbf{0} &\Lambda_+ \end{bmatrix} \\
    &= \begin{bmatrix} P_-^T Y_-^T &P_-^T Y_0^T &P_-^T Y_+^T\\ P_0^T Y_-^T &P_0^T Y_0^T &P_0^T Y_+^T\\ P_+^T Y_-^T &P_+^T Y_0^T &P_+^T Y_+^T \end{bmatrix}
    \begin{bmatrix} \diag(\ell_-) \\ ~ &\mathbf{0}\\ ~ &~ &\diag(\ell_+) \end{bmatrix}
    \begin{bmatrix} Y_- P_- &Y_-P_0 &Y_- P_+\\ Y_0 P_- &Y_0P_0 &Y_0 P_+\\ Y_+P_- &Y_+P_0 &Y_+ P_+ \end{bmatrix}, 
\end{aligned} \]
where the (3,1)-block equality is $\mathbf{0} = P_+^T Y_-^T \diag(\ell_-) Y_- P_- + P_+^T Y_+^T \diag(\ell_+) Y_+ P_-$, so 
\[ P_+^T Y_-^T \diag(\ell_-) + P_+^T Y_+^T \diag(\ell_+) Y_+ P_- (Y_- P_-)^{-1} = P_+^T (Y_-^T \diag(\ell_-) + Y_+^T M) = \mathbf{0}. 
\]
Then with $\nu I+H^\star \stackrel{\eqref{eq:Hstar+-}}{=} \nu(P_+P_+^T + P_0P_0^T + P_-P_-^T) + \nu(P_+P_+^T - P_-P_-^T) = \nu P_0P_0^T + 2\nu P_+ P_+^T$, we obtain
\begin{equation} \label{eq:mu1 - proof} \begin{aligned} 
    (\nu I + H^\star) (Y_-^T \diag(\ell_-) + Y_+^T M)
    \stackrel{\eqref{eq:P0}}{=}&\ 2\nu P_+ P_+^T (Y_-^T \diag(\ell_-) + Y_+^T M) \\
    ={}&\ 2\nu P_+ \mathbf{0} = \mathbf{0}, 
\end{aligned} \end{equation}  
which proves \eqref{eq:mu1 -} for all $j \in \cI_- \setminus \{0\}$. 
With \eqref{eq:mu1 - proof}, we also have 
\[ \begin{aligned} 
    (\nu I+H^\star) \left( \ell_0(\xx_0) \xx_0 + \sum_{i\in\cI_+} \mu_{i0} \xx_i \right) 
    &\leftstackrel{\eqref{eq:mu0 -}}{=} (\nu I+H^\star) \sum_{i\in\cI_+} \mu_{i0} (\xx_i-\xx_0) \\
    &= (\nu I+H^\star) \sum_{i\in\cI_+} \left(\ell_i(\xx_0) - \sum_{j\in\cI_-\setminus\{0\}} \mu_{ij}\right) (\xx_i-\xx_0) \\
    &= (\nu I+H^\star) Y_+^T (\ell_+ - M\mathbf{1}) \\
    &\leftstackrel{\eqref{eq:mu1 - proof}}{=} (\nu I+H^\star) (Y_+^T \ell_+ + Y_-^T \ell_-)
    \stackrel{\eqref{eq:Lagrange Y}}{=} \mathbf{0}, 
\end{aligned} \]
which proves \eqref{eq:mu1 -} for $j = 0$. 
\end{proof}

The function $\psi$ is defined and proved non-positive in Lemma~\ref{lem:psi}. 
It will be used to prove Theorem~\ref{thm:phase2} in conjunction with the parameters $\{\mu_{ij}\}$. 
\begin{lemma} \label{lem:psi}
Assume $f \in C^{1,1}_\nu(\R^n)$. 
For any $\u, \v \in \R^n$ and any matrix $H \in \R^{n \times n}$, we have 
\begin{equation} \label{eq:Lipscthiz stronger H} \begin{aligned} 
    \psi(\u,\v,H) \stackrel{\rm def}{=} &f(\u) -  f(\v) - \frac{1}{2\nu} [(\nu I-H)(\u-\v)] \cdot Df(\u) \\
    &- \frac{1}{2\nu} [(\nu I+H) (\u-\v)] \cdot Df(\v) 
    - \frac{1}{4\nu} \|H (\u - \v)\|^2 - \frac{\nu}{4} \|\u - \v\|^2 \le 0. 
\end{aligned} \end{equation} 
\end{lemma}
\begin{proof}
For the purpose of contradiction, assume \eqref{eq:Lipscthiz stronger H} is false. Then we have 
\[ \begin{aligned}
    - f(\u) <& - f(\v) - \frac{1}{2\nu} [(\nu I+H) (\u-\v)] \cdot Df(\u) \\
    &- \frac{1}{2\nu} [(\nu I-H)(\u-\v)] \cdot Df(\v) - \frac{1}{4\nu} \|H (\u - \v)\|^2 - \frac{\nu}{4} \|\u - \v\|^2. 
\end{aligned} \]
Swap $\u$ and $\v$ in this inequality and then add \eqref{eq:Lipschitz stronger} to the result. This leads to the contradiction 
$\frac{1}{4\nu} \|H (\v - \u) + (Df(\v) - Df(\u))\|^2 < 0$. 
\end{proof}

Finally, we prove the main result of this section, Theorem~\ref{thm:phase2}, which states $\underline{z}(\Theta,\xx_0)$ from \eqref{eq:phase2} is a sharp bound on the approximation error $|\hat{f}(\xx_0) - f(\xx_0)|$ when $\{\mu_{ij}\}$ are all non-negative. 
\begin{proof}[proof of Theorem~\ref{thm:phase2}]
Let $H^\star$ be the matrix defined in \eqref{eq:Hstar}. 
When $\mu_{ij} \ge 0$ for all $(i,j) \in \cI_+\times\cI_-$, the following inequality holds 
\begin{equation} \label{eq:phase2 summation}
    \sum_{i \in \cI_+} \sum_{j \in \cI_-} \mu_{ij} \psi(\xx_i,\xx_j,H^\star) \stackrel{\eqref{eq:Lipscthiz stronger H}}{\le} 0. 
\end{equation}
The zeroth-order term in the summation \eqref{eq:phase2 summation} is 
\[ \begin{aligned}
    \sum_{i \in \cI_+} \sum_{j \in \cI_-} \mu_{ij} (f(\xx_i) - f(\xx_j)) 
    &= \left[\sum_{i \in \cI_+} \sum_{j \in \cI_-} \mu_{ij} f(\xx_i)\right] - \left[\sum_{i \in \cI_+} \sum_{j \in \cI_-} \mu_{ij} f(\xx_j)\right] \\
    &\leftstackrel{\eqref{eq:mu0 +}\eqref{eq:mu0 -}}{=} \left[ \sum_{i \in \cI_+} \ell_i(\xx_0) f(\xx_i)\right] + \left[\sum_{j \in \cI_-} \ell_j(\xx_0) f(\xx_j)\right] \\
    &\leftstackrel{\eqref{eq:Lagrange m}}{=} \hat{f}(\xx_0) - f(\xx_0). 
\end{aligned} \] 
The sum of the first-order terms is $-1/(2\nu)$ multiplied by 
\[ \begin{aligned} 
    &\sum_{i \in \cI_+} \sum_{j\in\cI_-} \mu_{ij} \big( [(\nu I-H^\star)(\xx_i-\xx_j)] \cdot Df(\xx_i) + [(\nu I+H^\star)(\xx_i-\xx_j)] \cdot Df(\xx_j) \big) \\
    &= \left[\sum_{i \in \cI_+} \sum_{j\in\cI_-} \mu_{ij} [(\nu I-H^\star)\xx_i] \cdot Df(\xx_i) \right] 
    - \left[\sum_{i \in \cI_+} \sum_{j\in\cI_-} \mu_{ij} [(\nu I+H^\star)\xx_j] \cdot Df(\xx_j) \right] \\
    &\quad -\left[\sum_{i \in \cI_+} \sum_{j\in\cI_-} \mu_{ij}  [(\nu I-H^\star)\xx_j] \cdot Df(\xx_i) \right] + \left[\sum_{i \in \cI_+} \sum_{j\in\cI_-} \mu_{ij} [(\nu I+H^\star)\xx_i] \cdot Df(\xx_j) \right] \\
    &= \left[\sum_{i \in \cI_+} \ell_i(\xx_0) [(\nu I-H^\star)\xx_i] \cdot Df(\xx_i)  \right] 
    + \left[\sum_{j\in\cI_-} \ell_j(\xx_0) [(\nu I+H^\star)\xx_j] \cdot Df(\xx_j) \right] \\
    &\quad -\left[\sum_{i \in \cI_+} \ell_i(\xx_0) [(\nu I-H^\star)\xx_i] \cdot Df(\xx_i) \right] - \left[\sum_{j\in\cI_-} \ell_j(\xx_0) [(\nu I+H^\star)\xx_j] \cdot Df(\xx_j) \right] 
    = \mathbf{0}, 
\end{aligned} \]
where the second equality holds because of \eqref{eq:mu0 +}, \eqref{eq:mu0 -}, \eqref{eq:mu1 +}, and \eqref{eq:mu1 -} respectively for the four terms. 
Notice $H^{\star T} H^\star = \nu^2I$. 
The constant term in the summation \eqref{eq:phase2 summation} is $-1/2$ multiplies
\[ \begin{aligned} 
    &\sum_{i \in \cI_+} \sum_{j\in\cI_-} \mu_{ij} \left(\frac{1}{2\nu}\|H^\star(\xx_i-\xx_j)\|^2 + \frac{\nu}{2} \|\xx_i-\xx_j\|^2\right) \\
    &= \nu \left[ \sum_{i \in \cI_+} \sum_{j\in\cI_-} \mu_{ij} (\xx_i-\xx_j) \cdot \xx_i \right] - \nu \left[ \sum_{i \in \cI_+} \sum_{j\in\cI_-} \mu_{ij} (\xx_i-\xx_j) \cdot \xx_j \right] \\
    &\stackrel{\mathmakebox[\widthof{=}]{\scriptsize \begin{array}{c}\eqref{eq:mu0 +}\\\eqref{eq:mu0 -}\end{array}}}{=} \sum_{i \in \cI_+} \nu \left(\ell_i(\xx_0) \xx_i- \sum_{j\in\cI_-} \mu_{ij} \xx_j\right) \cdot \xx_i  - \sum_{j\in\cI_-} \nu \left(\sum_{i \in \cI_+} \mu_{ij} \xx_i + \ell_j(\xx_0) \xx_j\right) \cdot \xx_j \\
    &\leftstackrel{\mathmakebox[\widthof{=}]{\scriptsize \begin{array}{c}\eqref{eq:mu1 +}\\ \eqref{eq:mu1 -}\end{array}}}{=} \sum_{i \in \cI_+} \left[H^\star \left(\ell_i(\xx_0) \xx_i- \sum_{j\in\cI_-} \mu_{ij} \xx_j\right)\right] \cdot \xx_i + \sum_{j\in\cI_-} \left[H^\star \left(\sum_{i \in \cI_+} \mu_{ij} \xx_i + \ell_j(\xx_0) \xx_j\right)\right] \cdot \xx_j \\
    &= \left[ \sum_{i \in \cI_+} \ell_i(\xx_0) [H^\star \xx_i] \cdot \xx_i \right] + \left[ \sum_{j\in\cI_-} \ell_j(\xx_0) [H^\star \xx_j] \cdot \xx_j \right] 
    = G \cdot H^\star. 
\end{aligned} \]
Thus, the summation \eqref{eq:phase2 summation} is $\hat{f}(\xx_0) - f(\xx_0) \le \frac{1}{2} G \cdot H^\star = \underline{z}(\Theta,\xx_0)$. 
\end{proof}

\section{Sharp Error Bounds for Bivariate Extrapolation} \label{sec:phase3}
We investigate in this section the sharp error bounds when $\xx$ is in the four areas shown in Figure~\ref{fig:phase2}. 
This investigation is not just for the completeness of our analysis of the sharp error bound, but also to understand what type of function can be more difficult for linear interpolation to approximate than the quadratics. 

Before we start, we need the following proposition, which shows how the approximation error changes when $\xx_0$ is swapped with a point in $\Theta$. 
\begin{proposition} \label{thm:swap}
    Let $k\in\{1,\dots,n+1\}$ be a fixed index. 
    Assume there is an affinely independent sample set $\Theta_1 = \{\xx_1,\dots,\xx_{n+1}\}$ and a point $\xx_0\in\R^n$ such that $\Theta_2 := \Theta_1\setminus\{\xx_k\}\cup\{\xx_0\}$ is also affinely independent.  
    Let $\ell_k$ be the Lagrange polynomial (with respect to $\Theta_1$ not $\Theta_2$) corresponding to $\xx_k$. 
    Let $\hat{f}$ and $\hat{h}$ be the affine functions that interpolate some function $f:\R^n \rightarrow \R$ on $\Theta_1$ and $\Theta_2$, respectively. 
    The following two statements hold. 
    \begin{enumerate}
        \item The function approximation error of $\hat{h}$ at $\xx_k$ is the error of $\hat{f}$ at $\xx_0$ divided by $-\ell_k(\xx_0)$, i.e., $\hat{h}(\xx_k) - f(\xx_k) = (\hat{f}(\xx_0) - f(\xx_0)) / (-\ell_k(\xx_0))$. 
        \item If $|\hat{f}(\xx_0) - f(\xx_0)| = z(\Theta_1,\xx_0)$ holds for some $f\in C_\nu^{1,1}(\R^n)$, i.e. $f$ is an optimal solution to \eqref{prob:D}, then $|\hat{h}(\xx_k) - f(\xx_k)| = z(\Theta_2,\xx_k)$ holds for the same $f$. 
    \end{enumerate}
\end{proposition} 

\begin{proof}
Firstly, $-\ell_k(\xx_0) \neq 0$ because otherwise, by \eqref{eq:Lagrange Y}, vector $\xx_0$ would be an affine combination of $\{\xx_i\}_{i \in \{1,\dots,n+1\}\setminus\{k\}}$, violating the premise that $\Theta_2$ is affinely independent. 

If we divide $\hat{f}(\xx_0) - f(\xx_0) = \sum_{i=0}^{n+1} \ell_i(\xx_0) f(\xx_i)$ by $-\ell_k(\xx_0)$, the coefficient before $f(\xx_i)$ becomes $\alpha_i = -\ell_i(\xx_0)/\ell_k(\xx_0)$ for all $i = 0,1,\dots,n+1$. 
Since $\alpha_k = -\ell_k(\xx_0)/\ell_k(\xx_0) = -1$, $\sum_{i=0}^{n+1} \alpha_i=$ $-\ell_k(\xx_0)^{-1} \sum_{i=0}^{n+1}\ell_i(\xx_0) \stackrel{\eqref{eq:Lagrange 0}}{=}$ $0$, and $\sum_{i=0}^{n+1} \alpha_i \xx_i =$ $-\ell_k(\xx_0)^{-1} \sum_{i=0}^{n+1}\ell_i(\xx_0) \xx_i \stackrel{\eqref{eq:Lagrange Y}}{=} \mathbf{0}$, the coefficients $\{\alpha_i\}_{i=0,i\neq k}^{n+1}$ are the values of the Lagrange polynomials with respect to $\Theta_2$ at $\xx_k$. 
Thus, the quotient $(\hat{f}(\xx_0) - f(\xx_0)) / (-\ell_k(\xx_0))$ is exactly $\hat{h}(\xx_k) - f(\xx_k)$. 

It is assumed in the second statement that $f$ is an optimal solution to \eqref{prob:D}. 
The same $f$ must also be an optimal solution to the problem of finding the largest $|\hat{h}(\xx_k) - f(\xx_k)|$, since this optimization problem is simply \eqref{prob:D} with its objective function divided by the constant $-\ell_k(\xx_k)$. The sign of $-\ell_k(\xx_k)$ does not matter for the same reason the absolute sign can be ignored --- due to the symmetry of $C_\nu^{1,1}(\R^n)$. 
\end{proof}

Figure~\ref{fig:phase3 triangle} shows the case where $\xx_0$ is a fixed point inside the shaded triangle $\circnum{4}$ in Figure~\ref{fig:phase2}. 
This case is the same as when $\xx_0$ is inside the shaded triangle $\circnum{3}$ except the roles of $\xx_2$ and $\xx_3$ are swapped. 
Thus, if an expression of $z(\Theta,\xx_0)$ is found for the case when $\xx_0$ is in $\circnum{4}$, then we have also found it for the case when $\xx_0$ is in $\circnum{3}$. 
Similarly, if an expression of $z(\Theta,\xx_0)$ is found for the case when $\xx_0$ is in the translated cone $\circnum{2}$ in Figure~\ref{fig:phase2}, then we have also found it for the case when $\xx_0$ is in $\circnum{1}$. 
This reduces the cases that need to be studied to the two in Figure~\ref{fig:phase3}. 
Furthermore, after we obtain an expression of $z(\Theta,\xx_0)$ for the case in Figure~\ref{fig:phase3 triangle}, an expression for the case in Figure~\ref{fig:phase3 cone} can be obtained by switching the roles of $\xx_0$ and $\xx_2$ and apply Proposition~\ref{thm:swap}. 
Therefore, we only need to find an expression of $z(\Theta,\xx_0)$ for one of these four cases, and the rest can be easily derived from it. 

\begin{figure}[tbhp]
     \centering
     \subfloat[\raggedright When $\xx_0$ is in the open triangle such that $\ell_1(\xx_0)[\xx_2-\xx_1{]} \cdot [\xx_3-\xx_1{]} - \ell_3(\xx_0) [\xx_2-\xx_3{]} \cdot [\xx_1-\xx_3{]} < 0$, $\ell_2(\xx_0)>0$, and $\ell_3(\xx_0)<0$]{\label{fig:phase3 triangle}\resizebox{0.4\linewidth}{!}{\begin{tikzpicture}
\filldraw[black] (0,0) circle (2pt) node[anchor=north east] {$\xx_1$}; 
\filldraw[black] (2,1.8) circle (2pt) node[anchor=south east] {$\xx_2$}; 
\filldraw[black] (-2,0) circle (2pt) node[anchor=south] {$\xx_3$};  
\filldraw[black] (0.9,0) circle (2pt) node[anchor=north west] {$\tilde\w$}; 
\filldraw[black] (1.3889,0.8) circle (2pt) node[anchor=north west] {$\xx_0$}; 

\draw[thick] (0,0) -- (3.333,3);
\draw[thick] (2,1.8) -- (-2,0);
\draw[thick] (-2,0) -- (3.333,0);
\draw[thick, dashed] (2,-0.5) -- (2,3.5);
\draw[thick, dashed] (2,1.8) -- (0.9,0);

\end{tikzpicture}}}
     \hspace{0.1\linewidth}
     \subfloat[\raggedright When $\xx_0$ is in the open cone such that $\ell_1(\xx_0) [\xx_2-\xx_1{]} \cdot [\xx_3-\xx_1{]} - \ell_3(\xx_0) [\xx_2-\xx_3{]} \cdot [\xx_1-\xx_3{]} > 0$ and $\ell_3(\xx_0)>0$]{\label{fig:phase3 cone}\resizebox{0.4\linewidth}{!}{ \begin{tikzpicture}
\filldraw[black] (0,0) circle (2pt) node[anchor=north east] {$\xx_1$}; 
\filldraw[black] (2,1.8) circle (2pt) node[anchor=south east] {$\xx_2$}; 
\filldraw[black] (-2,0) circle (2pt) node[anchor=south] {$\xx_3$};  
\filldraw[black] (0.9,0) circle (2pt) node[anchor=north west] {$\tilde\w$}; 
\filldraw[black] (2.4278,2.5) circle (2pt) node[anchor=south west] {$\xx_0$}; 

\draw[thick] (0,0) -- (3.333,3);
\draw[thick] (2,1.8) -- (-2,0);
\draw[thick] (-2,0) -- (3.333,0);
\draw[thick, dashed] (2,-0.5) -- (2,3.5);
\draw[thick, dashed] (2.4278,2.5) -- (0.9,0);

\end{tikzpicture}}}
    \caption{Two configurations of $\Theta$ and $\xx_0$ where $\underline{z}(\Theta,\xx_0)$ defined in \eqref{eq:phase2} is less than the sharp error bound $z(\Theta,\xx_0)$ for bivariate extrapolation.}
    \label{fig:phase3}
\end{figure}

From now on, we focus on the case in Figure~\ref{fig:phase3 triangle}, which can be defined mathematically as $\ell_2(\xx_0)>0, \ell_3(\xx_0)<0$, and $\ell_1(\xx_0)[\xx_2-\xx_1] \cdot [\xx_3-\xx_1] - \ell_3(\xx_0)[\xx_2-\xx_3] \cdot [\xx_1-\xx_3] < 0$.
The following lemma shows the point $\tilde\w$, as defined in \eqref{eq:phase3 w}, is the intersection of the line going through $\xx_1$ and $\xx_3$ and the line going through $\xx_0$ and $\xx_2$. 

\begin{lemma}
Assume $-\ell_0(\xx_0)-\ell_2(\xx_0) \stackrel{\eqref{eq:Lagrange 0}}{=} \ell_1(\xx_0)+\ell_3(\xx_0)\neq 0$ for some affinely independent $\Theta = \{\xx_1,\xx_2,\xx_3\} \subset\R^2$ and $\xx_0\in\R^2$. 
Let 
\begin{equation} \label{eq:phase3 w}
    \tilde\w = \frac{-\ell_0(\xx_0)\xx_0+\ell_1(\xx_0)\xx_1-\ell_2(\xx_0)\xx_2+\ell_3(\xx_0)\xx_3}{-\ell_0(\xx_0)+\ell_1(\xx_0)-\ell_2(\xx_0)+\ell_3(\xx_0)}. 
\end{equation}
Then 
\[ \tilde\w = \frac{\ell_1(\xx_0)\xx_1+\ell_3(\xx_0)\xx_3}{\ell_1(\xx_0)+\ell_3(\xx_0)} = \frac{\ell_0(\xx_0)\xx_0+\ell_2(\xx_0)\xx_2}{\ell_0(\xx_0)+\ell_2(\xx_0)}, 
\]
and 
\begin{subequations} \label{eq:phase3 w y} \begin{align}
    \ell_0(\xx_0)[\xx_0-\tilde\w] + \ell_2(\xx_0)[\xx_2-\tilde\w] &= \mathbf{0}, \\
    \ell_1(\xx_0)[\xx_1-\tilde\w] + \ell_3(\xx_0)[\xx_3-\tilde\w] &= \mathbf{0}. 
\end{align} \end{subequations}
\end{lemma}
\begin{proof}
These equalities are direct results of \eqref{eq:Lagrange 0} and \eqref{eq:Lagrange Y}. 
\end{proof}

We define in the following lemma a square matrix $\tilde{H}^\star$ that plays a role similar to $H^\star$ in \eqref{eq:Hstar} but is asymmetric. 
\begin{lemma}
Assume for some affinely independent $\Theta = \{\xx_1,\xx_2,\xx_3\} \subset\R^2$ and $\xx_0\in\R^2$ that $\ell_2(\xx_0)>0, \ell_3(\xx_0)<0$, and $\ell_1(\xx_0)[\xx_2-\xx_1] \cdot [\xx_3-\xx_1] - \ell_3(\xx_0) [\xx_2-\xx_3] \cdot [\xx_1-\xx_3] < 0$. 
Let 
\begin{equation} \label{eq:Hstar 3}
    \tilde{H}^\star = \tilde{P} \begin{bmatrix} +\nu &0\\ 0 &-\nu \end{bmatrix} \tilde{P}^{-1} \text{ with } \tilde{P} = \begin{bmatrix} \xx_2-\xx_0 &\xx_1-\xx_3 \end{bmatrix}. 
\end{equation}
Let $\tilde\w$ be defined as \eqref{eq:phase3 w}. 
Then 
\begin{equation} \label{eq:phase3 Hstar eigvector} \begin{aligned} 
    \tilde{H}^\star(\xx_i-\tilde\w) &= \nu(\xx_i-\tilde\w) \text{ for } i\in\{0,2\}, \\
    \tilde{H}^\star(\xx_i-\tilde\w) &= -\nu(\xx_i-\tilde\w) \text{ for } i\in\{1,3\}. 
\end{aligned} \end{equation} 
\end{lemma}
\begin{proof}
It is clear from Figure~\ref{fig:phase3 triangle} that the assumption guarantees the invertibility of $\tilde{P}$ and $-\ell_0(\xx_0)-\ell_2(\xx_0) = \ell_1(\xx_0)+\ell_3(\xx_0) \neq 0$. 
Notice by the definition of $\tilde{H}^\star$, we have $\tilde{H}^\star(\xx_2-\xx_0) = \nu(\xx_2-\xx_0)$ and $\tilde{H}^\star(\xx_1-\xx_3) = -\nu(\xx_1-\xx_3)$.
The lemma holds true because $\xx_i-\w$ is parallel to $\xx_2-\xx_0$ for $i\in\{0,2\}$ and to $\xx_1-\xx_3$ for $i\in\{1,3\}$. 
\end{proof}

Now we are ready to show $G \cdot \tilde{H}^\star/2$ equals to $z(\Theta,\xx_0)$ for the case in Figure~\ref{fig:phase3 triangle}. 
\begin{theorem} \label{thm:phase3}
Assume $f \in C^{1,1}_\nu(\R^2)$. 
Let $\Theta = \{\xx_1,\xx_2,\xx_3\}\subset \R^2$ be a fixed set of affinely independent vectors such that $(\xx_2-\xx_1) \cdot (\xx_3-\xx_1) < 0$. 
Let $\hat{f}$ be the affine function that interpolates $f$ on $\Theta$. 
Let $\xx_0$ be a fixed point in $\R^2$ such that the Lagrange polynomials satisfy $\ell_2(\xx_0)>0, \ell_3(\xx_0)<0$, and $\ell_1(\xx_0) [\xx_2-\xx_1] \cdot [\xx_3-\xx_1] - \ell_3(\xx_0) [\xx_2-\xx_3] \cdot [\xx_1-\xx_3] < 0$. 
Let $G$ and $\tilde{H}^\star$ be the matrices defined in \eqref{eq:G} and \eqref{eq:Hstar 3}. 
Then the function approximation error of $\hat{f}$ at $\xx_0$ is bounded as 
\begin{equation} \label{eq:phase3}
    |\hat{f}(\xx_0) - f(\xx_0)| \le \frac{1}{2} G \cdot \tilde{H}^\star. 
\end{equation}
\end{theorem}

\begin{proof}
We only provide the proof for the case when $\hat{f}(\xx_0)-f(\xx_0) \ge 0$. 
We use the function $\psi$ defined in \eqref{eq:Lipscthiz stronger H} again. 
Since $\ell_3(\xx_0) < 0$, $[\xx_2-\xx_1] \cdot [\xx_3-\xx_1] < 0$, and 
\[ \begin{aligned}
    0 &> \ell_1(\xx_0) [\xx_2-\xx_1] \cdot [\xx_3-\xx_1] - \ell_3(\xx_0) [\xx_2-\xx_3] \cdot [\xx_1-\xx_3] \\
    &\leftstackrel{\eqref{eq:Lagrange 0}}{=}  (1-\ell_2(\xx_0)-\ell_3(\xx_0)) [\xx_2-\xx_1] \cdot [\xx_3-\xx_1] - \ell_3(\xx_0) [\xx_2-\xx_3] \cdot [\xx_1-\xx_3] \\
    &= (1-\ell_2(\xx_0)) [\xx_2-\xx_1] \cdot [\xx_3-\xx_1] - \ell_3(\xx_0) \|\xx_1-\xx_3\|^2, 
\end{aligned} \] 
we have $1-\ell_2(\xx_0) > 0$. Thus, the following inequalities hold: 
\begin{subequations} \label{eq:phase3 sum} \begin{align}
    (1-\ell_2(\xx_0)) \psi(\xx_1, \xx_0, \tilde{H}^\star) &\le 0, \\
    \ell_2(\xx_0) \psi(\xx_2, \xx_0, \tilde{H}^\star) &\le 0, \\ 
    -\ell_3(\xx_0) \psi(\xx_1, \xx_3, \tilde{H}^\star) &\le 0. 
\end{align} \end{subequations}
Similar to the previous proofs, we add these inequalities together. 
The sum of their zeroth-order terms is 
\[ \begin{aligned} 
&(1-\ell_2(\xx_0)) [f(\xx_1) - f(\xx_0)] + \ell_2(\xx_0) [f(\xx_2) - f(\xx_0)] - \ell_3(\xx_0) [f(\xx_1) - f(\xx_3)] \\
&= (1-\ell_2(\xx_0)-\ell_3(\xx_0)) f(\xx_1) + \ell_2(\xx_0) f(\xx_2) + \ell_3(\xx_0) f(\xx_3) - f(\xx_0) 
\stackrel{\eqref{eq:Lagrange m}\eqref{eq:Lagrange 0}}{=} \hat{f}(\xx_0) - f(\xx_0). 
\end{aligned} \]
The sum of their first-order terms is $-1/(2\nu)$ multiplies 
\[ \begin{aligned} 
&\hspace{-0.5em} (1-\ell_2(\xx_0)) \left\{ [(\nu I-\tilde{H}^\star) (\xx_1-\xx_0)] \cdot Df(\xx_1) + [(\nu I+\tilde{H}^\star) (\xx_1-\xx_0)] \cdot Df(\xx_0) \right\} \\
&\hspace{-0.5em} + \ell_2(\xx_0) \left\{ [(\nu I-\tilde{H}^\star) (\xx_2-\xx_0)] \cdot Df(\xx_2) + [(\nu I+\tilde{H}^\star) (\xx_2-\xx_0)] \cdot Df(\xx_0) \right\} \\
&\hspace{-0.5em} -\ell_3(\xx_0) \left\{ [(\nu I-\tilde{H}^\star) (\xx_1-\xx_3)] \cdot Df(\xx_1) + [(\nu I+\tilde{H}^\star) (\xx_1-\xx_3)] \cdot Df(\xx_3) \right\} \\
&= \{(\nu I-\tilde{H}^\star) [(1-\ell_2(\xx_0))(\xx_1-\xx_0) - \ell_3(\xx_0)(\xx_1-\xx_3)]\} \cdot Df(\xx_1) \\
&\quad + \ell_2(\xx_0) [(\nu I-\tilde{H}^\star) (\xx_2-\xx_0)] \cdot Df(\xx_2) - \ell_3(\xx_0) [(\nu I+\tilde{H}^\star) (\xx_1-\xx_3)] \cdot Df(\xx_3) \\
&\quad + \left\{ (\nu I+\tilde{H}^\star) [(1-\ell_2(\xx_0))(\xx_1-\xx_0) + \ell_2(\xx_0)(\xx_2-\xx_0)] \right\} \cdot Df(\xx_0) \\ 
&\leftstackrel{\eqref{eq:Lagrange 0}\eqref{eq:Lagrange Y}}{=} \ell_2(\xx_0) [(\nu I-\tilde{H}^\star) (\xx_0 - \xx_2)] \cdot Df(\xx_1) + \ell_2(\xx_0) [(\nu I-\tilde{H}^\star) (\xx_2-\xx_0)] \cdot Df(\xx_2) \\
&\quad - \ell_3(\xx_0) [(\nu I+\tilde{H}^\star) (\xx_1-\xx_3)] \cdot Df(\xx_3) + \ell_3(\xx_0) [(\nu I+\tilde{H}^\star) (\xx_1 - \xx_3)] \cdot Df(\xx_0) \\ 
&\leftstackrel{\eqref{eq:phase3 Hstar eigvector}}{=} \mathbf{0}. 
\end{aligned} \]
Let $\tilde\w$ be defined as \eqref{eq:phase3 w}. 
The sum of the constant terms is $-1/2$ multiplied by 
\[ \begin{aligned} 
&\hspace{-1em} (1-\ell_2(\xx_0)) \left[\frac{1}{2\nu} \|\tilde{H}^\star (\xx_1-\xx_0)\|^2 + \frac{\nu}{2} \|\xx_1-\xx_0\|^2 \right] + \ell_2(\xx_0) \left[\frac{1}{2\nu} \|\tilde{H}^\star (\xx_2-\xx_0)\|^2 \right. \\
&\hspace{-1em} \left.+ \frac{\nu}{2} \|\xx_2-\xx_0\|^2 \right] - \ell_3(\xx_0) \left[\frac{1}{2\nu} \|\tilde{H}^\star (\xx_1-\xx_3)\|^2 + \frac{\nu}{2} \|\xx_1-\xx_3\|^2 \right] \\
&\leftstackrel{\eqref{eq:phase3 Hstar eigvector}}{=} (1-\ell_2(\xx_0))\left\{ -\tilde{H}^\star(\xx_1-\tilde{\w}) \cdot (\xx_1-\tilde{\w}) + \tilde{H}^\star(\xx_0-\w) \cdot (\xx_0-\tilde{\w}) \right\} \\
&\qquad + \ell_2(\xx_0) \tilde{H}^\star(\xx_2-\xx_0) \cdot (\xx_2-\xx_0) + \ell_3(\xx_0) \tilde{H}^\star(\xx_1-\xx_3) \cdot (\xx_1-\xx_3) \\
&\leftstackrel{\eqref{eq:Lagrange 0}}{=} \tilde{H}^\star[\ell_3(\xx_0)(\xx_1-\xx_3)-(\ell_1(\xx_0)+\ell_3(\xx_0))(\xx_1-\tilde{\w})] \cdot (\xx_1-\tilde{\w}) \\
&\qquad + \tilde{H}^\star[(1-\ell_2(\xx_0))(\xx_0-\tilde{\w}) - \ell_2(\xx_0)(\xx_2-\xx_0)] \cdot (\xx_0-\tilde{\w}) \\
&\qquad - \ell_3(\xx_0) \tilde{H}^\star(\xx_1-\xx_3) \cdot (\xx_3-\tilde{\w}) + \ell_2(\xx_0) \tilde{H}^\star(\xx_2-\xx_0) \cdot (\xx_2-\tilde{\w}) \\
&\leftstackrel{\eqref{eq:Lagrange 0}\eqref{eq:Lagrange Y}}{=} 0 + 0 - \ell_3(\xx_0) [\tilde{H}^\star(\xx_1-\tilde{\w}) - \tilde{H}^\star(\xx_3-\tilde{\w})] \cdot (\xx_3-\tilde{\w}) \\
&\qquad + \ell_2(\xx_0) [\tilde{H}^\star(\xx_2-\w) - \tilde{H}^\star(\xx_0-\tilde{\w})] \cdot (\xx_2-\tilde{\w}) \\ 
&\leftstackrel{\eqref{eq:phase3 w y}}{=} \sum_{i=0}^{3} \ell_i(\xx_0) \tilde{H}^\star (\xx_i-\tilde{\w}) \cdot (\xx_i-\tilde{\w})
\stackrel{\eqref{eq:G recenter}}{=} G \cdot \tilde{H}^\star. 
\end{aligned} \]
Thus, the sum of the inequalities in \eqref{eq:phase3 sum} is \eqref{eq:phase3} when $\hat{f}(\xx_0) - f(\xx_0) \ge 0$. 
\end{proof} 

We show in Theorem~\ref{thm:phase3 sharp} the upper bound \eqref{eq:phase3} can be achieved by a piecewice quadratic, and \eqref{eq:phase3} is therefore sharp. 
\begin{theorem} \label{thm:phase3 sharp} 
Under the setting of Theorem~\ref{thm:phase3}, the bound \eqref{eq:phase3} is sharp and can be achieved by 
\[ f(\u) = \left\{ \begin{aligned}
    &\frac{\nu}{2} \|\u-\tilde{\w}\|^2 - \frac{\nu[(\xx_1-\xx_3) \cdot (\u-\tilde{\w})]^2}{\|\xx_1-\xx_3\|^2}  &&\text{if } (\u-\tilde{\w}) \cdot (\xx_1-\xx_3) \le 0,  \\
    &\frac{\nu}{2} \|\u-\tilde{\w}\|^2 &&\text{if } (\u-\tilde{\w}) \cdot (\xx_1-\xx_3) \ge 0,
\end{aligned} \right. 
\]
where $\tilde{\w}$ is defined in \eqref{eq:phase3 w}. 
\end{theorem}
\begin{proof}
The function approximation error for this piecewise quadratic function is 
\[ \begin{aligned}
    &\hat{f}(\xx_0) - f(\xx_0)
    = \sum_{i=0}^{n+1} \ell_i(\xx_0) f(\xx_i) \\
    &= \frac{\nu}{2} \sum_{i=0}^3 \!\ell_i(\xx_0) \|\xx_i-\tilde{\w}\|^2 \!-\! \frac{\nu\ell_1(\xx_0) [(\xx_1\!-\!\xx_3) \!\cdot\! (\xx_1\!-\!\tilde{\w})]^2}{\|\xx_1-\xx_3\|^2} \!-\! \frac{\nu\ell_3(\xx_0) [(\xx_1\!-\!\xx_3) \!\cdot\! (\xx_3\!-\!\tilde{\w})]^2}{\|\xx_1-\xx_3\|^2} \\
    &= \frac{\nu}{2} \sum_{i=0}^3 \ell_i(\xx_0) \|\xx_i-\tilde{\w}\|^2 - \nu\ell_1(\xx_0)\|\xx_1-\tilde{\w}\|^2 - \nu\ell_3(\xx_0)\|\xx_3-\tilde{\w}\|^2 \\
    &= \frac{\nu}{2} \left(\ell_0(\xx_0)\|\xx_0-\tilde{\w}\|^2 - \ell_1(\xx_0)\|\xx_1-\tilde{\w}\|^2 + \ell_2(\xx_0)\|\xx_2-\tilde{\w}\|^2 - \ell_3(\xx_0)\|\xx_3-\tilde{\w}\|^2\right) \\
    &\leftstackrel{\eqref{eq:phase3 Hstar eigvector}}{=} \frac{1}{2} \sum_{i=0}^3 \ell_i(\xx_0) \|\xx_i-\tilde{\w}\|_{\tilde{H}^\star} 
    \stackrel{\eqref{eq:G recenter}}{=} \frac{1}{2} G \cdot \tilde{H}^\star. 
\end{aligned} \]
Now we prove $f \in C_\nu^{1,1}(\R^n)$. 
Firstly, it is clear that $f$ is continuous on $\R^2$ and differentiable on the two half spaces $\{\u:~ (\u-\tilde{\w}) \cdot (\xx_1-\xx_3) < 0\}$ and $\{\u:~ (\u-\tilde{\w}) \cdot (\xx_1-\xx_3) > 0\}$. 
Then given any $\u$ such that $(\u-\tilde{\w}) \cdot (\xx_1-\xx_3) = 0$, it can be calculated for any $\v\in\R^2$ that 
\begin{equation*} \begin{aligned} 
    |f(\u+\v) &- f(\u) - \nu (\u-\tilde{\w})\cdot\v| \\
    &= \left\{ \begin{aligned} 
        &-\frac{\nu}{2}\|\v\|^2 - \frac{\nu [(\xx_1-\xx_3)\cdot \v]^2}{\|\xx_1-\xx_3\|^2} &&\text{if }  (\u+\v-\tilde{\w}) \cdot (\xx_1-\xx_3) \le 0,  \\
        &-\frac{\nu}{2}\|\v\|^2 &&\text{if }  (\u+\v-\tilde{\w}) \cdot (\xx_1-\xx_3) \ge 0. 
    \end{aligned} \right.
\end{aligned} \end{equation*} 
Thus 
\[ \lim_{\v\rightarrow\mathbf{0}} \frac{|f(\u+\v) - f(\u) - \nu(\u-\tilde{\w})\cdot\v|}{\|\v\|} = 0, 
\]
which shows $f$ is differentiable with gradient $\nu(\u-\tilde{\w})$ on $\{\u:~ (\u-\tilde{\w}) \cdot (\xx_1-\xx_3) = 0\}$. 
The condition \eqref{eq:Lipschitz} is clearly satisfied if $\u_1$ and $\u_2$ are in the same half space. 
Now assume $(\u_1-\tilde{\w})\cdot(\xx_1-\xx_3) < 0$ and $(\u_2-\tilde{\w})\cdot(\xx_1-\xx_3) > 0$. 
Then, we have 
\[ \begin{aligned} 
&\|Df(\u_1) - Df(\u_2)\|^2 \\
&= \|\nu(\u_1-\tilde{\w}) - 2\nu\left[(\xx_1-\xx_3) \cdot (\u_1-\tilde{\w})/\|\xx_1-\xx_3\|^2\right] (\xx_1-\xx_3)- \nu(\u_2-\tilde{\w})\|^2 \\
&= \nu^2\|\u_1-\u_2\|^2 + 4\nu^2 [(\u_1-\tilde{\w}) \cdot (\xx_1-\xx_3)] [(\u_2-\tilde{\w})\cdot(\xx_1-\xx_3)] /\|\xx_1-\xx_3\|^2 \\
&< \nu^2\|\u_1-\u_2\|^2, 
\end{aligned} \]
which shows \eqref{eq:Lipschitz} always holds.
Therefore $f \in C_\nu^{1,1}(\R^n)$. 
\end{proof}

\section{Application in Derivative-Free Optimization} \label{sec:app}
In this section, we demonstrate the application of our theories in DFO by providing a complexity analysis of a basic simplicial search method similar to the original one proposed in \cite{spendley1962sequential}, as described in Algorithm~\ref{alg}. 
\begin{algorithm}[ht]
    \caption{A Basic Simplicial Search Method} \label{alg}
    \begin{algorithmic}[1]
      \State \textbf{input: } the objective function $f:\R^n\to\R$, a starting point $\mathbf{c}_0\in\R^n$, and the radius $\delta > 0$ 
      \State Construct a regular simplex with center $\mathbf{c}_0$ and radius $\delta$. Let its $n+1$ vertices be the starting sample set $\Theta_0 \subset \R^n$, and evaluate $f$ at these points. 
        \For{$k = 0,1,\dots$}
          \State Sort and label the points in $\Theta_k$ as $\{\xx_i\}_{i=1}^{n+1}$ such that $f(\xx_1) \le f(\xx_2) \le \cdots \le f(\xx_{n+1})$. 
          \State Compute the reflection point $\xx_0 = - \xx_{n+1} + \frac{2}{n} \sum_{i=1}^n \xx_i$, and evaluate $f(\xx_0)$.
          \State $\Theta_{k+1} \gets \Theta_k \setminus \{\xx_{n+1}\} \cup \{\xx_0\}$. 
        \EndFor
    \end{algorithmic}
\end{algorithm}

\subsection{Regular Simplex} 
Algorithm~\ref{alg} starts with a regular simplex --- a simplex in which all edges have the same length. 
While regular simplices are classical geometric objects, existing treatments appear in specialized literature (e.g., \cite{coxeter1973regular}) where the results are considerably more advanced than required here.
For completeness and clarity, we therefore include in this subsection a self-contained discussion of regular simplices. 
Recall $J$ represents the all one matrix. 
We start with the following formal definition of regular simplex, which is much more complicated than the classical one, but is more practical in the upcoming analysis. 
\begin{definition}[Regular Simplex] \label{def:regular simplex}
    An $n$-dimensional simplex with vertices $\Theta = \{\xx_1,\dots,\xx_{n+1}\} \subset \R^n$ is said to be a regular simplex with center $\bc\in\R^n$ and radius $\delta>0$ if $\bc = \frac{1}{n+1} \sum_{i=1}^{n+1} \xx_i$ and there exists an $(n+1)$-dimensional orthogonal matrix $\bQ$ (i.e. $\bQ^{-1} = \bQ^T$) such that $\bA = \bB_\delta \bQ$, where 
    \begin{equation} \label{eq:regular simplex}
        \bA = \begin{bmatrix}
            0 &(\xx_1 - \bc)^T \\ 
            0 &(\xx_2 - \bc)^T \\ 
            \vdots &\vdots \\
            0 &(\xx_{n+1} - \bc)^T 
        \end{bmatrix}, 
        \bB = \bI_{n+1} - \frac{\bJ_{n+1}}{n+1}, 
        \text{ and }
        B_\delta = \sqrt{\frac{n+1}{n}} \delta B. 
    \end{equation}
\end{definition}

This definition is not intuitive. 
However, the equality 
\[ (\xx_i - \bc)\cdot(\xx_j - \bc) = [\bA\bA^T]_{ij} = [\bB_\delta\bQ\bQ^T\bB_\delta^T]_{ij} = [\bB_\delta \bB_\delta^T]_{ij} = 
    \left\{ \begin{aligned}
        &\delta^2 &&\text{if } i = j, \\
        - &\delta^2 / n &&\text{if } i \neq j
    \end{aligned} \right. 
\] 
reveals the essential properties 
\begin{enumerate}
    \item the distance from the center to any vertex $\|\xx_i - \bc\| = \delta, ~\forall i$; 
    \item the length of any edge $\|\xx_i - \xx_j\| = \|(\xx_i-\bc) - (\xx_j-\bc)\| = \big[\|\xx_i-\bc\|^2 + \|\xx_j-\bc\|^2 - 2 (\xx_i - \bc)\cdot(\xx_j - \bc)\big]^{1/2} = (\delta^2 + \delta^2 + 2\delta^2/n)^{1/2} = \sqrt{\frac{2n+2}{n}} \delta, ~\forall i\neq j$. 
\end{enumerate} 

Moreover, this definition implies a method for constructing the regular simplex required at the beginning of Algorithm~\ref{alg}, as explained in the following proposition. 
\begin{proposition}[Construction of Regular Simplex]
    An $n$-dimensional regular simplex with center $\bc$ and radius $\delta$ can be constructed by first using the Gram–Schmidt process to decompose the matrix $\bB_\delta$ in \eqref{eq:regular simplex} into an orthogonal matrix $\bar{\bQ}$ and an upper triangular matrix $\bR \in \R^{(n+1)\times(n+1)}$ such that $\bar{\bQ} \bR = \bB_\delta$, then setting the $i$th vertex of the simplex as $\xx_i = \bc + \begin{bmatrix} \bI_n &\mathbf{0} \end{bmatrix} \bR \e_i$ for $i = 1,\dots,n+1$. 
\end{proposition}
\begin{proof}
    First notice the matrix $B$ would map any vector parallel to $\mathbf{1}$ to $\mathbf{0}$, i.e., $\bB\mathbf{1} = \mathbf{0}$, and any vector $\v \in \R^{n+1}$ to $\v$ if $\v \cdot \mathbf{1} = 0$. 
    With $B_\delta\mathbf{1} = \sqrt{\frac{n+1}{n}} \delta B\mathbf{1} = \mathbf{0}$, we can prove the first condition in Definition~\ref{def:regular simplex}: 
    \begin{equation*} \begin{aligned} 
        \frac{1}{n+1} \sum_{i=1}^{n+1} \xx_i 
        &= \frac{1}{n+1} \sum_{i=1}^{n+1} \left[\bc + \begin{bmatrix} \bI_n &\mathbf{0} \end{bmatrix} \bR \e_i\right] 
        = \bc + \frac{1}{n+1} \begin{bmatrix} \bI_n &\mathbf{0} \end{bmatrix} \bR \mathbf{1} \\ 
        &= \bc + \frac{1}{n+1} \begin{bmatrix} \bI_n &\mathbf{0} \end{bmatrix} \bar{\bQ}^{-1} \bB_\delta \mathbf{1} = \bc. 
    \end{aligned} \end{equation*}
    Secondly, for any $\u = (u_1, u_2, \dots, u_n, 0)^T \neq \mathbf{0}$, the vector $B \u = \u + \frac{\mathbf{1}^T\u}{n+1} \mathbf{1}$ is never $\mathbf{0}$, because either $\mathbf{1}^T\u=0$, and the vector equals to $\u \neq \mathbf{0}$; or $\mathbf{1}^T\u \neq 0$, and $(n+1)$-st element of the vector equals to $\frac{\mathbf{1}^T\u}{n+1} \neq 0$. 
    This proves the first $n$ columns of $B$ are linearly independent. 
    Since we also have $B\mathbf{1} = \mathbf{0}$, the square matrix $\bB$ of order $n+1$ has rank $n$. 
    Thus, the upper triangular matrix $\bR$ of order $n+1$ resulted from the Gram–Schmidt process also has rank $n$, and its last element must be zero, making the last row of $\bR$ zero. 
    The matrix $\bA$ from \eqref{eq:regular simplex} then satisfies 
    \begin{equation*} \begin{aligned} 
        \bA &= 
        \begin{bmatrix}
            \mathbf{0} &\bR^T \begin{pmatrix} \bI_n\\ \mathbf{0}^T \end{pmatrix} 
        \end{bmatrix}
        = \bR^T \begin{bmatrix} \mathbf{0} &\bI_n\\ 1 &\mathbf{0}^T\end{bmatrix}
        = (\bar{\bQ} \bR)^T \bar{\bQ} \begin{bmatrix} \mathbf{0} &\bI_n\\ 1 &\mathbf{0}^T\end{bmatrix} \\
        &= \bB_\delta^T \bar{\bQ} \begin{bmatrix} \mathbf{0} &\bI_n\\ 1 &\mathbf{0}^T\end{bmatrix}
        = \bB_\delta \bar{\bQ} \begin{bmatrix} \mathbf{0} &\bI_n\\ 1 &\mathbf{0}^T\end{bmatrix}, 
    \end{aligned} \end{equation*}
    where the last equality holds because $\bB_\delta$ is symmetric. 
    Since the matrix $\bar{\bQ} \begin{bmatrix} \mathbf{0} &\bI_n\\ 1 &\mathbf{0}^T\end{bmatrix}$ is orthogonal, the second condition in Definition~\ref{def:regular simplex} is also met. 
\end{proof} 

We now derive the properties held by a regular simplex that will be used in the subsequent complexity analysis. 

\begin{lemma}[Properties of Regular Simplex] \label{lem:regular simplex}
    Let $\Theta = \{\xx_1,\dots,\xx_{n+1}\} \subset \R^n$ be any set of points whose convex hull is an $n$-dimensional regular simplex with radius $\delta>0$. 
    Let $\bc := \frac{1}{n+1}\sum_{i=1}^{n+1} \xx_i$ and $\xx_0 := - \xx_{n+1} + \frac{2}{n} \sum_{i=1}^n \xx_i$. 
    Then, the following statements are true. 
    \begin{enumerate}
        \item The convex hull of $\Theta \setminus \{\xx_{n+1}\} \cup \{\xx_0\} = \{\xx_0, \xx_1, \dots, \xx_n\}$ is also an $n$-dimensional regular simplex with radius $\delta$. 
        \item The eigenvalues $\{\lambda_i\}_{i=1}^n$ of the matrix $G$ defined in \eqref{eq:G} satisfy 
            \begin{equation}
                \lambda_1 = - \frac{2(n+1)^2}{n^2} \delta^2, \lambda_2 = \cdots = \lambda_n = \frac{2(n+1)}{n^2} \delta^2. 
            \end{equation}
        \item The parameters $\{\mu_{ij}\}$ defined in Theorem~\ref{thm:phase2} satisfy 
            \begin{equation}
                \mu_{ij} = 1/n \text{ for all } (i,j) \in \cI_+ \times \cI_-. 
            \end{equation}
        \item For any real numbers $a_1,\dots, a_{n+1}$, we have 
            \begin{equation} \label{eq: NM g variance}
                \left\|\sum_{i=1}^{n+1} a_i D\ell_i(\xx_0)\right\|^2 
                = \frac{n}{\delta^2} \cdot \frac{1}{n+1} \sum_{i=1}^{n+1} \bigg(a_i - \frac{1}{n+1}\sum_{j=1}^{n+1} a_j\bigg)^2. 
            \end{equation}
    \end{enumerate}
\end{lemma}
\begin{proof}
    Denote $\bp = \begin{pmatrix} -\mathbf{1}_n \\ n \end{pmatrix} \in \R^{n+1}$. 
    Let $A,B, B_\delta$ be defined as in \eqref{eq:regular simplex}. 
    Then, there exists an orthogonal matrix $Q$ such that $A = B_\delta Q$. 

    (1) 
    The simplex $\text{conv}(\{\xx_0, \xx_1, \dots, \xx_n\})$ has a new center 
    \begin{equation*} \begin{aligned} 
        \tilde\bc &= \frac{1}{n+1} \sum_{i=0}^{n} \xx_i 
        = \frac{1}{n+1} \left( - \xx_{n+1} + \frac{2}{n} \sum_{i=1}^n \xx_i + \sum_{i=1}^{n} \xx_i \right) \\
        &= \frac{-1}{n+1} \xx_{n+1} + \frac{n+2}{n(n+1)} \sum_{i=1}^n \xx_i, 
    \end{aligned} \end{equation*}
    which implies 
    $   \xx_0 - \tilde\bc 
        = (- \xx_{n+1} + \frac{2}{n} \sum_{i=1}^n \xx_i) - (\frac{-1}{n+1} \xx_{n+1} + \frac{n+2}{n(n+1)} \sum_{i=1}^n \xx_i)
        = \frac{-n}{n+1} \xx_{n+1} + \frac{1}{n+1} \sum_{i=1}^n \xx_i. 
    $
    Define the symmetric matrix of order $n+1$
    \begin{equation*}
        \bS = 
        \begin{bmatrix}
            \bI_n - \frac{n+2}{n(n+1)}\bJ_n  &\frac{1}{n+1} \mathbf{1} \\
            \frac{1}{n+1} \mathbf{1}^T &\frac{-n}{n+1} 
        \end{bmatrix} 
        = B  - \frac{2 \bp\bp^T}{n(n+1)} 
    \end{equation*}
    and notice $\bS \mathbf{1} = \mathbf{0}$. 
    Recall the symmetric matrix $\bB$ would map any vector parallel to $\mathbf{1}$ to $\mathbf{0}$ and any vector $\v \in \R^{n+1}$ to $\v$ if $\v \cdot \mathbf{1} = 0$.  
    This implies $B_\delta \mathbf{1} = \sqrt{\frac{n+1}{n}}\delta B \mathbf{1} = \mathbf{0}$ and that $B_\delta J_{n+1}$ is a zero matrix. 
    It also implies $BB = B$ and $B \bp = \bp$, which further implies $SB = S = BS$ and $SB_\delta = B_\delta S$. 
    Then, the matrix 
    \begin{equation*} \begin{aligned} 
        \begin{bmatrix}
            0 &(\xx_1 - \tilde\bc)^T \\ 
            0 &(\xx_2 - \tilde\bc)^T \\ 
            \vdots &\vdots \\
            0 &(\xx_n - \tilde\bc)^T \\ 
            0 &(\xx_0 - \tilde\bc)^T 
        \end{bmatrix}
        &= \bS \bA
        ~ = \bS \bB_\delta \bQ 
        = \bB_\delta \bS \bQ \\
        &= \bB_\delta \left(\bS + \frac{\bJ_{n+1}}{n+1} \right) \bQ
        = \bB_\delta \left(I_{n+1} - \frac{2\bp\bp^T}{n(n+1)} \right) \bQ, 
    \end{aligned} \end{equation*}
    where, with $\bp^T\bp = n(n+1)$, it can be shown that $\left(I_{n+1} - \frac{2\bp\bp^T}{n(n+1)} \right) \bQ$ is orthogonal: 
    \begin{equation*} \begin{aligned} 
        &Q^T\left(I_{n+1} - \frac{2\bp\bp^T}{n(n+1)} \right)^T \left(I_{n+1} - \frac{2\bp\bp^T}{n(n+1)} \right) Q \\
        &= Q^T \left(I_{n+1} - \frac{2\bp\bp^T}{n(n+1)} - \frac{2\bp\bp^T}{n(n+1)} + \frac{4\bp\bp^T}{n(n+1)}\right) Q
         = I_{n+1}. 
    \end{aligned} \end{equation*}
    Thus, by Definition~\ref{def:regular simplex}, the convex hull of $\{\xx_0, \xx_1, \dots, \xx_n\}$ is a regular simplex with radius $\delta$. 

    (2) 
    Since $\xx_0 = - \xx_{n+1} + \frac{2}{n} \sum_{i=1}^n \xx_i$ is defined in Algorithm~\ref{alg} with the normalized barycentric coordinates with respect to $\{\xx_1,\dots,\xx_{n+1}\}$, we have 
    \[ \ell_1(\xx_0) = \cdots = \ell_n(\xx_0) = \frac{2}{n} 
        \text{ and } \ell_{n+1}(\xx_0) = -1. 
    \]
    This can also be verified by checking that the above set of values indeed satisfy \eqref{eq:Lagrange 0} and \eqref{eq:Lagrange Y} and recalling that the Lagrange polynomials are unique. 

    With these $\ell_i(\xx_0)$, consider the matrix $\tilde{\bG}$ defined below, which admits the eigendecomposition on the right hand-side of the equation, 
    \begin{equation*}
        \tilde{\bG}
        := \begin{bmatrix} 0 &\mathbf{0}^T\\ \mathbf{0} &\bG \end{bmatrix} 
        = \begin{bmatrix}
            1 &\mathbf{0}^T \\
            \mathbf{0} &\bP
        \end{bmatrix}
        \begin{bmatrix}
            0 &\mathbf{0}^T \\
            \mathbf{0} &\Lambda
        \end{bmatrix}
        \begin{bmatrix}
            1 &\mathbf{0}^T \\
            \mathbf{0} &\bP^T
        \end{bmatrix}, 
    \end{equation*}
    where $\Lambda$ and $P$ are the matrices containing the eigenvalues and eigenvectors of $G$, same as their definitions for equation \eqref{eq:Hstar}. 
    This matrix has a zero eigenvalue associated with the direction $\e_1$, and its other $n$ eigenvalues are identical to the eigenvalues of $G$. 

    With the equality  
    \begin{equation*} \begin{aligned} 
        -\xx_0 + \bc 
        &= \xx_{n+1} - \frac{2}{n} \sum_{i=1}^n \xx_i + \frac{1}{n+1} \sum_{i=1}^{n+1} \xx_i 
        = \frac{n+2}{n+1} \xx_{n+1} - \frac{n+2}{n(n+1)} \sum_{i=1}^n \xx_i \\ 
        &= \frac{n+2}{n+1} (\xx_{n+1} - \bc) - \frac{n+2}{n(n+1)} \sum_{i=1}^n (\xx_i - \bc) \\
        &= \frac{n+2}{n(n+1)} \begin{bmatrix}
            (\xx_1 - \bc)
            &\cdots 
            &(\xx_{n+1} - \bc)
        \end{bmatrix} 
        \begin{pmatrix} -\mathbf{1}\\ n \end{pmatrix}, 
    \end{aligned} \end{equation*}
    which implies $\begin{pmatrix} 0\\ -\xx_0 + \bc \end{pmatrix}= \frac{n+2}{n(n+1)} A^T \bp$, we can show that $\tilde{G}$ can also be written as 
    \begin{equation*} \begin{aligned} 
        \tilde{G}
        &= \begin{bmatrix} 0 &\mathbf{0}^T\\ \mathbf{0} &\sum_{i=0}^{n+1} \ell_i(\xx_0) \xx_i \xx_i^T \end{bmatrix}
        \stackrel{\eqref{eq:G recenter}}{=} \begin{bmatrix} 0 &\mathbf{0}^T\\ \mathbf{0} &\sum_{i=0}^{n+1} \ell_i(\xx_0) (\xx_i - \bc) (\xx_i - \bc)^T \end{bmatrix} \\
        &= \begin{bmatrix} 0 &\mathbf{0}^T\\ \mathbf{0} &\ell_0(\xx_0) (-\xx_0 + \bc) (-\xx_0 + \bc)^T \end{bmatrix} + A^T \diag(\ell_1(\xx_0), \dots, \ell_{n+1}(\xx_0)) A \\
        &= - \frac{(n+2)^2}{n^2(n+1)^2} A^T \bp \bp^T A + A^T \diag\left(\frac{2}{n}, \dots, \frac{2}{n}, -1\right) A \\
        &= A^T \left(\frac{2}{n} I_{n+1} - \frac{n+2}{n} \e_{n+1} \e_{n+1}^T - \frac{(n+2)^2}{n^2 (n+1)^2} \bp \bp^T \right) A. 
    \end{aligned} \end{equation*} 

    Let $D = B - \frac{n+2}{n(n+1)} \bp\bp^T$. 
    Since $A = B_\delta Q$ holds with an orthogonal $Q$, $BB=B$, $B\e_{n+1} = \frac{\bp}{n+1}$, and $B\bp=\bp$, we can show that $\tilde{G}$ is similar to $\frac{2(n+1)\delta^2}{n^2} D$:
    \begin{equation*} \begin{aligned} 
        \tilde{G} 
        &= Q^T B_\delta^T \left(\frac{2}{n} I_{n+1} - \frac{n+2}{n} \e_{n+1} \e_{n+1}^T - \frac{(n+2)^2}{n^2 (n+1)^2} \bp \bp^T \right) B_\delta Q\\
        &= \frac{(n+1)\delta^2}{n} Q^T \left( 
            \frac{2}{n} B - \frac{n+2}{n} \cdot \frac{\bp\bp^T}{(n+1)^2} 
            - \frac{(n+2)^2}{n^2 (n+1)^2} \bp\bp^T
        \right) Q \\
        &= \frac{(n+1) \delta^2}{n} Q^{-1} ~\frac{2}{n}\left(B - \frac{n+2}{n(n+1)} \bp \bp^T \right) Q \\
        & = Q^{-1} \left( \frac{2(n+1) \delta^2}{n^2}D \right) Q. 
    \end{aligned} \end{equation*}
    
    Matrix $D = B - \frac{n+2}{n(n+1)} \bp\bp^T = I_{n+1} - \frac{J_{n+1}}{n+1} - \frac{n+2}{n(n+1)} \bp\bp^T$ has eigenvalues $0, -n-1, 1$ with multiplicities $1,1,n-1$, respectively, because 
    \begin{itemize}
        \item $D \mathbf{1} = \mathbf{1} - \frac{\|\mathbf{1}\|^2}{n+1} \mathbf{1} - \mathbf{0} = \mathbf{0}$, 
        \item $D\bp = I_{n+1} \bp - \frac{J_{n+1}}{n+1} \bp - \frac{n+2}{n(n+1)} \left\|\bp\right\|^2 \bp = \bp - \mathbf{0} - (n+2) \bp = -(n+1) \bp$, 
        \item and for any $\v$ such that $\mathbf{1} \cdot \v = \bp \cdot \v = 0$ there is $D \v = \v - \mathbf{0} - \mathbf{0} = \v$. 
    \end{itemize}
    As similar matrices share the same eigenvalues, $\tilde{G}$ has the same eigenvalues as $\frac{2(n+1)\delta^2}{n^2} D$, and they are $0, - 2(n+1)^2 \delta^2 / n^2$, and $2(n+1)\delta^2 / n^2$ with multiplicities $1,1$, and $n-1$, respectively. 
    Combining this with the relationship between the eigenvalues of $\tilde{G}$ and $G$, claim (2) is proved. 

    (3) 
    Recall the definition $M = \diag(\ell_+) Y_+ P_- (Y_- P_-)^{-1}$ in Theorem~\ref{thm:phase2}. 
    With the knowledge of $\{\ell_i(\xx_0)\}_{i=1}^{n+1}$ (at the beginning of the proof of part (2)), we have $\diag(\ell_+) = \diag(\ell_1(\xx_0),\dots,\ell_n(\xx_0)) = \frac{2}{n} I_n$, 
    \[  Y_+ = \begin{bmatrix}
            (\xx_1-\xx_0)^T \\
            (\xx_2-\xx_0)^T \\
            \vdots \\
            (\xx_n-\xx_0)^T 
        \end{bmatrix}, 
        \text{ and } 
        Y_- = (\xx_{n+1} - \xx_0)^T. 
    \]
    To find the expression for $P_-$, we follow the proof of part (2). 
    Since 
    \begin{align*}
        \tilde{G} Q^{-1} \bp
        &= Q^{-1} Q \tilde{G} Q^{-1} \bp
        = Q^{-1} \frac{2(n+1) \delta^2}{n^2}D \bp \\ 
        &= Q^{-1} \frac{2(n+1) \delta^2}{n^2} (-n-1) \bp
        = - \frac{2(n+1)^2 \delta^2}{n^2} Q^{-1} \bp, 
    \end{align*}
    the eigenvector associated with the only negative eigenvalue $- 2(n+1)^2 \delta^2 / n^2$ of $\tilde{G}$ is $Q^{-1} \bp$. 
    Hence $P_- = \sqrt{\frac{1}{n(n+1)}} Q^{-1} \bp$, after scaling the vector to unit length. 

    Subsequently, 
    \begin{equation*} \begin{aligned} 
        \begin{bmatrix}
            Y_+ \\ Y_- 
        \end{bmatrix} P_- 
        &= \Big( A + \mathbf{1} \begin{bmatrix} 0 &(\bc - \xx_0)^T\end{bmatrix} \Big) P_- 
        = \Big( A + \mathbf{1} \frac{n+2}{n(n+1)} \begin{bmatrix} -\mathbf{1}^T &n \end{bmatrix} A \Big) P_-  \\
        &= \Big( I + \frac{n+2}{n(n+1)} \mathbf{1} \bp^T \Big) A P_-  
        = \Big( I + \frac{n+2}{n(n+1)} \mathbf{1} \bp^T \Big) A \sqrt{\frac{1}{n(n+1)}} Q^{-1} \bp \\
        &= \Big( I + \frac{n+2}{n(n+1)} \mathbf{1} \bp^T \Big) \frac{\delta}{n} \bp 
        = \frac{\delta}{n} \big( \bp + (n+2) \mathbf{1} \big) 
        = \frac{(n+1)\delta}{n}\begin{bmatrix} \mathbf{1}\\ 2 \end{bmatrix}, 
    \end{aligned} \end{equation*}
    where the second equality holds by the definition of $\bc$ and $\xx_0$, and the fifth equality holds because \vspace{-\baselineskip} 
    \begin{equation*} 
        \sqrt{\frac{1}{n(n+1)}} A Q^{-1} \bp
        = \sqrt{\frac{1}{n(n+1)}} B_\delta \bp
        = \frac{\delta}{n} B \bp 
        = \frac{\delta}{n} \bp. 
    \end{equation*}
    Thus, the matrix $M$ from Theorem~\ref{thm:phase2}, which is of size $n \times 1$ in this case, is 
    \begin{equation*} 
        M = \diag(\ell_+) Y_+ P_- (Y_- P_-)^{-1} 
        = \frac{2}{n} \Big(\frac{(n+1)\delta}{n} \mathbf{1}\Big) \Big(\frac{(n+1)\delta}{n} \cdot 2 \Big)^{-1} 
        = \frac{1}{n} \mathbf{1}. 
    \end{equation*}
    This proves $\mu_{ij} = 1/n$ for all $i \in \cI_+$ and $j \in \cI_- \setminus\{0\}$. 
    Finally, we also have $\mu_{i0} = \ell_i(\xx_0) - \sum_{j \in \cI_- \setminus\{0\}} \mu_{ij} = 2/n - 1/n = 1/n$ for all $i \in \cI_+$. 
    
    (4) 
    Again, following the proof for part (2), since 
    \begin{align*}
        \tilde{G} Q^{-1} \mathbf{1}
        = Q^{-1} Q \tilde{G} Q^{-1} \mathbf{1} 
        = Q^{-1} \frac{2(n+1) \delta^2}{n^2}D \mathbf{1} 
        = \mathbf{0}, 
    \end{align*}
    the eigenvector associated with the only zero eigenvalue of $\tilde{G}$ is $Q^{-1} \mathbf{1}$. 
    As we also know from the definition of $\tilde{G}$ that $\tilde{G} \e_1 = \mathbf{0}$, it can be deduced that 
    \begin{equation*}
        Q^T \mathbf{1} = Q^{-1} \mathbf{1} = \sqrt{n+1} \e_1 
        \text{ or, equivalently, }
        Q \e_1 = \sqrt{\frac{1}{n+1}} \mathbf{1}. 
    \end{equation*}
    Recall, from the definition of $\Phi$ and the Lagrange polynomials in Section~\ref{sec:preliminaries}, the affine function $\ell_i$ has derivative $D\ell_i(\u) = \begin{bmatrix} \mathbf{0} &I_n \end{bmatrix} \Phi^{-1} \e_i$ for all $\u\in\R^n$ and $i=1,\dots,n+1$. 
    The matrix $Q$ can be written as 
    \begin{equation*} \begin{aligned}
        \Phi &= A + \mathbf{1} \begin{bmatrix} 1 &(\bc - \xx_0)^T \end{bmatrix} 
        = A + \mathbf{1} \left(\e_1^T + \frac{n+2}{n(n+1)} \begin{bmatrix}-\mathbf{1}^T &n \end{bmatrix} A \right) \\
        &= \left(I_{n+1} + \frac{n+2}{n(n+1)} \mathbf{1} \bp^T \right) B_\delta \bQ + \mathbf{1} \e_1^T \bQ^T \bQ \\
        &= \sqrt{\frac{n+1}{n}} \delta \left( B + \frac{n+2}{n(n+1)} \mathbf{1} \bp^T \right)\bQ + \frac{\mathbf{1} \mathbf{1}^T}{\sqrt{n+1}} \bQ \\
        &= \sqrt{\frac{n+1}{n}} \delta \left[ I_{n+1} + \frac{1}{n+1} \mathbf{1}\left(\frac{n+2}{n} \bp^T - \mathbf{1}^T + \frac{\sqrt{n}}{\delta} \mathbf{1}^T\right) \right]\bQ, 
    \end{aligned} \end{equation*}
    so, by the Sherman–Morrison formula, its inverse is 
    \begin{equation*} \begin{aligned}
        \Phi^{-1} &= \sqrt{\frac{n}{n+1}} \cdot \frac{1}{\delta} Q^T \left(I_{n+1} - \frac{\frac{1}{n+1}\mathbf{1} \left(\frac{n+2}{n} \bp^T - \mathbf{1}^T + \frac{\sqrt{n}}{\delta} \mathbf{1}^T\right)}{1 + \left(\frac{n+2}{n} \bp^T - \mathbf{1}^T + \frac{\sqrt{n}}{\delta} \mathbf{1}^T\right) \frac{1}{n+1}\mathbf{1}} \right) \\
        &= \sqrt{\frac{n}{n+1}} \cdot \frac{1}{\delta} Q^T \left(I_{n+1} - \frac{\delta}{\sqrt{n}(n+1)}\mathbf{1} \left(\frac{n+2}{n} \bp^T - \mathbf{1}^T + \frac{\sqrt{n}}{\delta} \mathbf{1}^T\right) \right). 
    \end{aligned} \end{equation*}
    Then, \vspace{-\baselineskip} 
    \begin{equation*} \begin{aligned}
        \begin{bmatrix} 0\\ D\ell_i(\xx_0) \end{bmatrix} 
        &= \begin{bmatrix} 0 &\mathbf{0}^T\\ \mathbf{0} &I_n \end{bmatrix} \Phi^{-1} \e_i 
        = \sqrt{\frac{n}{n+1}} \cdot \frac{1}{\delta} \left(\bQ^T - \frac{\e_1\mathbf{1}^T}{\sqrt{n+1}} \right) \e_i, 
    \end{aligned} \end{equation*}
    where the second equality holds because 
    \[ \begin{bmatrix} 0 &\mathbf{0}^T\\ \mathbf{0} &I_n \end{bmatrix} \bQ^T = \bQ^T - \begin{bmatrix} \e_1^T\\ \mathbf{0} \end{bmatrix} \bQ^T = \bQ^T - \frac{\e_1\mathbf{1}^T}{\sqrt{n+1}} 
    \]
    and \vspace{-\baselineskip} 
    \[  \left(\bQ^T - \frac{\e_1\mathbf{1}^T}{\sqrt{n+1}}\right) \mathbf{1} 
        = \sqrt{n+1} \e_1 - \sqrt{n+1} \e_1 
        = \mathbf{0}. 
    \] 
    Thus, \vspace{-\baselineskip} 
    \begin{equation*} \begin{aligned}
        D&\ell_i(\xx_0) \cdot D\ell_j(\xx_0) 
        = \begin{bmatrix} 0\\ D\ell_i(\xx_0) \end{bmatrix} \cdot \begin{bmatrix} 0\\ D\ell_j(\xx_0) \end{bmatrix} \\ 
        &= \frac{n}{(n+1)\delta^2} \e_i^T \left(Q^T - \frac{\e_1\mathbf{1}^T}{\sqrt{n+1}} \right)^T \left(Q^T - \frac{\e_1\mathbf{1}^T}{\sqrt{n+1}} \right) \e_j \\
        &= \frac{n}{(n+1)\delta^2} \e_i^T \left(\bQ\bQ^T - \frac{\mathbf{1}\e_1^T}{\sqrt{n+1}}\bQ^T - \bQ\frac{\e_1\mathbf{1}^T}{\sqrt{n+1}} + \frac{\bJ_{n+1}}{n+1} \right) \e_j \\
        &= \frac{n}{(n+1)\delta^2} \e_i^T \left(\bI_{n+1} - \frac{\bJ_{n+1}}{n+1}\right) \e_j. 
    \end{aligned} \end{equation*}
    Finally, for any vector $\mathbf{a} = (a_1,\dots, a_{n+1})^T$, we have 
    \begin{equation*} \begin{aligned}
        \left\|\sum_{i=1}^{n+1} a_i D\ell_i(\xx_0)\right\|^2 
        &= \frac{n}{(n+1)\delta^2} \a^T \left(\bI_{n+1} - \frac{\bJ_{n+1}}{n+1}\right) \a \\ 
        &= \frac{n}{(n+1)\delta^2} \sum_{i=1}^{n+1} \bigg(a_i - \frac{1}{n+1}\sum_{j=1}^{n+1} a_j\bigg)^2. 
    \end{aligned} \end{equation*}
\end{proof}

The first property in the lemma above shows that a regular simplex maintains its regularity and size after a reflection operation.  
Since Algorithm~\ref{alg} starts with a regular simplex with radius $\delta$ and only uses the reflection operation, we can conclude that the simplex in any of its iterations is regular and has radius $\delta$.  
Therefore, the other three properties in this lemma would hold across all iterations. 

\subsection{Complexity Analysis} 
For any iteration $k$, let $\hat{f}$ be the affine function that interpolates $f$ on $\Theta_k$. 
While simplicial search methods do not use linear interpolation per se, the new point $\xx_0$ is chosen based on $n+1$ affinely independent points (vertices of $\Theta_k$) and their function values. 
The function value $f(\xx_0)$ would be $\hat{f}(\xx_0)$ plus the approximation error $f(\xx_0) - \hat{f}(\xx_0)$. 
This allows us to bound the range of $f(\xx_0)$ with the theories in Section~\ref{sec:phase2}.

\begin{lemma} \label{lem: NM f bound}
    Assume $f \in C^{1,1}_\nu(\R^n)$. In any iteration of Algorithm~\ref{alg}, the function value at the reflection point $\xx_0$ is always bounded as 
    \begin{equation} \label{eq: NM f bound}
        \bigg| f(\xx_0) + f(\xx_{n+1}) - \frac{2}{n}\sum_{i=1}^n f(\xx_i) \bigg| \le \frac{2n+2}{n} \nu \delta^2. 
    \end{equation}
\end{lemma}
\begin{proof}
    Let $\hat{f}$ be the affine function that interpolates $f$ on $\{\xx_1,\dots,\xx_{n+1}\}$. 
    Then $\hat{f}(\xx_0) = -f(\xx_{n+1}) + \frac{2}{n}\sum_{i=1}^n f(\xx_i)$. 
    Since it has been established in Lemma~\ref{lem:regular simplex} that $\mu_{ij} = 1/n\ge0$ for all $(i,j)\in\cI_+\times\cI_-$, $\underline{z}(\Theta,\xx_0)$ defined in \eqref{eq:phase2} is a sharp upper bound on $|\hat{f}(\xx_0) - f(\xx_0)|$ according to Theorem~\ref{thm:phase2}, meaning 
    \[  |\hat{f}(\xx) - f(\xx)| \le \frac{\nu}{2} \sum_{i=1}^n |\lambda_i| = \frac{\nu}{2} \cdot \frac{2(n+1)}{n^2} \delta^2 \big[|-(n+1)| + (n-1)\cdot1\big] = \frac{2(n+1)}{n} \nu \delta^2, 
    \]
    where the first equality is a result of Lemma~\ref{lem:regular simplex} (2). 
\end{proof}

We link the function values at the simplex points $\Theta_k$ to the gradient of the affine function $\hat{f}$ in the following lemma. 
\begin{lemma} \label{lem: NM f g}
    Assume $f \in C^{1,1}_\nu(\R^n)$. For any iteration $k$ of Algorithm~\ref{alg}, let $\hat{f}$ be the affine function that interpolates $f$ on $\Theta_k$, and $\bc_k$ the center of $\Theta_k$. Then 
    \begin{equation} \label{eq: NM f g}
        \|D\hat{f}(\bc_k)\| \le \frac{n}{\delta} \Big[f(\xx_{n+1}) - \frac{1}{n+1} \sum_{i=1}^{n+1} f(\xx_i)\Big]. 
    \end{equation}
\end{lemma}
\begin{proof}
    We first show that given any set of real numbers $\{a_i\}_{i=1}^{n+1} \subset \R$ such that $a_1 \le a_2 \le \cdots \le a_{n+1}$, it always holds that $a_{n+1} \ge \mu + \sigma/\sqrt{n}$, where $\mu = (n+1)^{-1} \sum_{i=1}^{n+1} a_i$ and $\sigma^2= (n+1)^{-1} \sum_{i=1}^{n+1} (a_i - \mu)^2$ are their average and variance. 

    Since $\mu$ is the average, there must be a $1 \le k \le n$ such that $a_1 \le \cdots \le a_k \le \mu \le a_{k+1} \le \cdots \le a_{n+1}$. 
    By the definition of variance, 
    \[ \begin{aligned}
        (n+1)\sigma^2 &= \sum_{i=1}^k (a_i - \mu)^2 + \sum_{i=k+1}^{n+1} (a_i - \mu)^2 \\
        &\le \Big[\sum_{i=1}^k (a_i - \mu)\Big]^2 + \sum_{i=k+1}^{n+1} (a_i - \mu)^2 
        = \Big[\sum_{i=k+1}^{n+1} (a_i - \mu)\Big]^2 + \sum_{i=k+1}^{n+1} (a_i - \mu)^2 \\
        &\le (n-k+1+1) \sum_{i=k+1}^{n+1} (a_i - \mu)^2 
        \le (n-k+2)(n-k+1)(a_{n+1}-\mu)^2, 
    \end{aligned} \]
    where the first inequality holds because $\{a_i-\mu\}_{i=1}^k$ share the same sign, the second equality is due to $\sum_{i=1}^{n+1} (a_i - \mu) = 0$, and the second inequality is due to $\bJ_{n-k+1} \preceq (n-k+1) \bI_{n-k+1}$. 
    Then, since $k$ is at least 1, we have $(n+1)\sigma^2 \le (n+1)n (a_{n+1} - \mu)^2$, which implies $a_{n+1} \ge \mu + \sigma/\sqrt{n}$. 

    Now, let $\hat{f}$ be the affine function that interpolates $f$ on $\Theta_k$. 
    Its gradient follows 
    \[ \begin{aligned} 
        \|D\hat{f}(\bc_k)\|^2 
        &= \|D\hat{f}(\xx_0)\|^2 
        \stackrel{\eqref{eq:Lagrange m}}{=} \bigg\|\sum_{i=1}^{n+1} f(\xx_i) D\ell_i(\xx_0)\bigg\|^2 \\
        &\leftstackrel{\eqref{eq: NM g variance}}{=} \frac{n}{\delta^2} \cdot \frac{1}{n+1} \sum_{i=1}^{n+1} \bigg(f(\xx_i) - \frac{1}{n+1}\sum_{j=1}^{n+1} f(\xx_j)\bigg)^2, 
    \end{aligned} \]
    implying the variance of $\{f(\xx_i)\}_{i=1}^{n+1}$ is $\|D\hat{f}(\bc_k)\|^2 \delta^2 / n$. 
    Thus, $f(\xx_{n+1}) \ge \frac{1}{n+1}\sum_{j=1}^{n+1} f(\xx_j) + \|D\hat{f}(\bc_k)\| \delta / n$, which can be rearranged to \eqref{eq: NM f g}. 
\end{proof}

The following lemma provides a bound on the difference between the gradient $Df(\bc_k)$ and its estimate $D\hat{f}(\bc_k)$ from the linear model that interpolates $\Theta_k$. 
\begin{lemma} \label{lem: NM g bound}
    Assume $f \in C^{1,1}_\nu(\R^n)$. For any iteration $k$ of Algorithm~\ref{alg}, let $\hat{f}$ be the affine function that interpolates $f$ on $\Theta_k = \{\xx_1,\dots,\xx_{n+1}\}$, and $\bc_k$ the center of $\Theta_k$. Then 
    \begin{equation} \label{eq: NM g bound}
        \|Df(\bc_k) - D\hat{f}(\bc_k)\|^2 \le \frac{n}{4} \nu^2 \delta^2. 
    \end{equation}
\end{lemma}
\begin{proof}
    Let $t$ be the first-order Taylor expansion of $f$ at $\bc_k$: 
    \[ t(\u) = f(\bc_k) + Df(\bc_k) \cdot (\u-\bc_k) \text{ for all } \u\in\R^n. 
    \]
    Then, \vspace{-\baselineskip}
    \begin{equation} \begin{aligned}
        \|&Df(\bc_k) - D\hat{f}(\bc_k)\|^2
        = \|Df(\bc_k) - D\hat{f}(\xx_0)\|^2 
        \leftstackrel{\eqref{eq:Lagrange m}}{=} \bigg\|Df(\bc_k) - \sum_{i=1}^{n+1} f(\xx_i) D\ell_i(\xx_0)\bigg\|^2 \\
        &= \bigg\|D(f-t)(\bc_k) - \sum_{i=1}^{n+1} (f-t)(\xx_i) D\ell_i(\xx_0)\bigg\|^2 
        = \bigg\|\sum_{i=1}^{n+1} (f-t)(\xx_i) D\ell_i(\xx_0)\bigg\|^2 \\
        &\leftstackrel{\eqref{eq: NM g variance}}{=} \frac{n}{\delta^2} \cdot \frac{1}{n+1}\sum_{i=1}^{n+1} \bigg(a_i - \frac{1}{n+1}\sum_{j=1}^{n+1} a_j\bigg)^2, 
    \end{aligned} \end{equation} 
    where, for $i=1,\dots,n+1$, $a_i$ represents $(f-t)(\xx_i)$, and $|a_i|\le\nu\delta^2/2$ holds by \eqref{eq:Lipschitz quadratic}. As the variance of real numbers inside an interval of length $b$ is bounded by $b^2/4$, we have \eqref{eq: NM g bound}.
\end{proof}


Finally, we can establish that the norm of the gradient at the center of the simplex $\Theta_k$ decreases to a neighborhood of zero at a rate of $\mathcal{O}(1/k)$. 
\begin{theorem} \label{thm:NM}
    Assume $f \in C^{1,1}_\nu(\R^n)$ and $f(\u) \ge f^\star$ for all $\u\in\R^n$. 
    Let $\bc_k$ be the center of $\Theta_k$ for iteration $k=0,1,\dots$ of Algorithm~\ref{alg}. We have for any $k\ge1$ 
    \begin{equation} \label{eq: MN convergence}
        \frac{1}{k} \sum_{t=0}^{k-1} \|Df(\bc_t)\| \le \frac{n^2}{2\delta k} \cdot \Big[\frac{1}{n+1}\sum_{\u\in\Theta_0} f(\u) - f^\star\Big] + \Big(n + \frac{\sqrt{n}}{2}\Big) \nu \delta. 
    \end{equation}
\end{theorem}
\begin{proof}
    Let $\Theta_{k} = \{\xx^{(k)}_i\}_{i=1}^{n+1}$ be the vertices of the simplex in iteration $k$, and $\xx^{(k)}_0$ their reflection point. 
    By Lemma~\ref{lem: NM f bound}, 
    \begin{equation} \label{eq: NM descent} \begin{aligned} 
        \sum_{\u\in\Theta_{k+1}} f(\u) 
        &= \sum_{\u\in\Theta_{k}} f(\u) - f(\xx^{(k)}_{n+1}) + f(\xx^{(k)}_0) \\
        &\leftstackrel{\eqref{eq: NM f bound}}{\le} \sum_{\u\in\Theta_{k}} f(\u) - f(\xx_{n+1}^{(k)}) + \Big[-f(\xx^{(k)}_{n+1}) + \frac{2}{n}\sum_{i=1}^n f(\xx^{(k)}_i) + \frac{2n+2}{n} \nu \delta^2\Big] \\
        &= \sum_{\u\in\Theta_{k}} f(\u) - \frac{2n+2}{n} \Big[f(\xx^{(k)}_{n+1}) - \frac{1}{n+1} \sum_{i=1}^{n+1} f(\xx^{(k)}_i) \Big] + \frac{2n+2}{n} \nu \delta^2. 
    \end{aligned} \end{equation}
    By telescoping, we have 
    \[
        \sum_{\u \in \Theta_{k}} f(\u) \le \sum_{\u\in\Theta_0} f(\u) - \frac{2n+2}{n} \sum_{t=0}^{k-1} \Big[f(\xx^{(t)}_{n+1}) - \frac{1}{n+1} \sum_{i=1}^{n+1} f(\xx^{(t)}_i) \Big] + k\frac{2n+2}{n} \nu \delta^2. 
    \]
    Since $\sum_{\u \in \Theta_{k}} f(\u) \ge (n+1)f^\star$, the above inequality implies
    \begin{equation} \label{eq: NM convergence f}
        \frac{1}{k} \sum_{t=0}^{k-1} \Big[f(\xx^{(t)}_{n+1}) - \frac{1}{n+1} \sum_{i=1}^{n+1} f(\xx^{(t)}_i) \Big] \le \frac{n}{2 k} \cdot \Big[\frac{1}{n+1}\sum_{\u\in\Theta_0} f(\u) - f^\star\Big] + \nu \delta^2. 
    \end{equation}
    Finally, use Lemmas~\ref{lem: NM f g} and \ref{lem: NM g bound} to obtain 
    \[ \begin{aligned} 
        f(\xx^{(t)}_{n+1}) - \frac{1}{n+1} \sum_{i=1}^{n+1} f(\xx^{(t)}_i)~
        \leftstackrel{\eqref{eq: NM f g}}{\ge} \frac{\delta}{n} \|D\hat{f}(\bc_t)\| 
        &\ge \frac{\delta}{n} \big(\|Df(\bc_t)\| - \|Df(\bc_t) - D\hat{f}(\bc_t)\| \big) \\
        &\leftstackrel{\eqref{eq: NM g bound}}{\ge} \frac{\delta}{n} \Big(\|Df(\bc_t)\| - \frac{\sqrt{n}}{2}\nu\delta \Big), 
    \end{aligned} \]
    which can be combined with \eqref{eq: NM convergence f} to obtain \eqref{eq: MN convergence}. 
\end{proof}

If the Lipschitz constant $\nu$ is known, we can select the size of the simplex and a stopping criterion to obtain a solution of desired accuracy. 
Following Theorem~\ref{thm:NM}, the next shows the number of iterations Algorithm~\ref{alg} takes to find an $\varepsilon$-stationary point. 
\begin{theorem} \label{thm:NM2}
    Assume $f \in C^{1,1}_\nu(\R^n)$ and $f(\u) \ge f^\star$ for all $\u\in\R^n$. 
    Given a desired accuracy $\varepsilon > 0$, if we set the simplex radius $\delta = \frac{2\varepsilon}{5n\nu}$ and add stopping criterion $\big[f(\xx_{n+1}) - \frac{1}{n+1} \sum_{i=1}^{n+1} f(\xx_i) \big] \le 2\nu\delta^2$ before the computation of the reflection point, then Algorithm~\ref{alg} would terminate in iteration $k$ with 
    \begin{equation} \|Df(\bc_k)\| \le \varepsilon
        \quad \text{ and } \quad 
        k \le \frac{25n^3\nu}{8\varepsilon^2}\left[\frac{1}{n+1}\sum_{\u\in\Theta_0} f(\u) - f^\star\right]. 
    \end{equation}
\end{theorem}
\begin{proof}
    With our choice of $\delta$, the difference between the true gradient at the center and the gradient of the model $\hat{f}$ follows 
    \[ \|Df(\bc_k) - D\hat{f}(\bc_k)\| \stackrel{\eqref{eq: NM g bound}}{\le} \frac{\sqrt{n}}{2} \nu \delta \le \frac{n}{2} \nu \delta = \frac{\varepsilon}{5}; 
    \]
    and when the algorithm terminates, 
    \begin{equation} \label{eq: NM stop g}
        \|D\hat{f}(\bc_k)\| 
        \stackrel{\eqref{eq: NM f g}}{\le} \frac{n}{\delta} \Big[f(\xx_{n+1}) - \frac{1}{n+1} \sum_{i=1}^{n+1} f(\xx_i)\Big]
        \le 2n\nu\delta = \frac{4}{5}\varepsilon. 
    \end{equation}
    Therefore, $\|Df(\bc_k)\| \le \|D\hat{f}(\bc_k)\| + \|Df(\bc_k) - D\hat{f}(\bc_k)\| \le \varepsilon$. 
    Moreover, the inequality \eqref{eq: NM convergence f} implies by iteration $\bar{k} = \Big\lceil \frac{25n^3\nu}{8\epsilon^2}\big[\frac{1}{n+1}\sum_{\u\in\Theta_0} f(\u) - f^\star\big] \Big\rceil$, where $\lceil \cdot \rceil$ denotes the ceiling function that returns the least integer greater than or equal to the input, 
    \[ \begin{aligned} 
        \frac{1}{\bar{k}} \sum_{t=0}^{\bar{k}-1} \Big[f(\xx^{(t)}_{n+1}) - \frac{1}{n+1} \sum_{i=1}^{n+1} f(\xx^{(t)}_i) \Big] 
        &\le \frac{n}{2\bar{k}} \cdot \Big[\frac{1}{n+1}\sum_{\u\in\Theta_0} f(\u) - f^\star\Big] + \nu \delta^2 \\
        &\le \frac{4\varepsilon^2}{25n^2\nu}  + \nu \delta^2 
        = 2 \nu \delta^2. 
    \end{aligned} \]
    Then, since $f(\xx^{(t)}_{n+1}) - \frac{1}{n+1} \sum_{i=1}^{n+1} f(\xx^{(t)}_i) \ge 0$ for all $t$, the stopping criterion must be satisfied at least once in iterations $0, \dots, \bar{k}-1$. 
\end{proof}

The stopping criterion in Theorem~\ref{thm:NM2} can also be changed to $\|D\hat{f}(\bc_k)\| \le 4\varepsilon/5$, as in \eqref{eq: NM stop g}. 
It is easier to satisfy and would not weaken the theoretical guarantees on the algorithmic complexity or the solution accuracy, but requires the computation of $D\hat{f}$. 

To the best of our knowledge, the analysis in this section is the first to show the complexity of any simplicial search method. 
The two main theorems provide quantitative insights into the effects of dimension, nonlinearity, simplex size, and initial sample on the algorithm’s performance. 
Our theories on the approximation error of linear interpolation play a crucial role in bounding the function value at the new point generated in each iteration (Lemma~\ref{lem: NM f bound}), which is a prerequisite for the important inequality \eqref{eq: NM descent}. 

\subsection{Comparison with Similar Results} 
Theorems~\ref{thm:NM} and \ref{thm:NM2} show the ability of Algorithm~\ref{alg} in finding near-stationary points when optimizing (possibly) nonconvex functions with Lipschitz continuous gradient. 
This is different from most analyses of simplicial search methods, which aim at showing the convergence of the points $\Theta_k$. 
Our analysis instead resembles that of noisy gradient methods. 

For example, the randomized derivative-free method proposed in \cite{nesterov2017random} starts with a point $\xx_0\in\R^n$ and then applies the following update rule 
\begin{equation} \label{eq: GFO}
    \xx_{k+1} = \xx_k - h_k g_\mu(\xx_k) 
    \quad \text{ with } 
    g_\mu(\xx_k) = \frac{f(\xx_k + \mu \u) - f(\xx_k)}{\mu} \u, 
\end{equation}
where $\mu > 0$ is considered to be a smoothing parameter, $\u$ is a random vector following the standard multivariate normal distribution, and $h_k > 0$ is the step size at iteration $k$. 
It is shown in Section 7 of \cite{nesterov2017random} that, when $f \in C_\nu^{1,1}(\R^n)$ and $h_0 = h_1 = \cdots = \frac{1}{4(n+4)\nu}$, the following inequality holds: 
\begin{equation} \label{eq: GFO convergence}
    \frac{1}{k} \sum_{t=0}^{k-1} \mathbb{E} \|D f(\xx_t)\|^2 \le \frac{16(n+4)\nu}{k} \Big( f(\xx_0) - f^\star \Big) + \bigg[ \frac{3(n+4)^2}{2} + \frac{(n+6)^3}{2} \bigg] \nu^2 \mu^2. 
\end{equation}
It bears much resemblance to \eqref{eq: MN convergence} from our Theorem~\ref{thm:NM}: 
\begin{equation} \tag{\ref{eq: MN convergence}}
    \frac{1}{k} \sum_{t=0}^{k-1} \|Df(\bc_t)\| \le \frac{n^2}{2\delta k} \cdot \Big[\frac{1}{n+1}\sum_{\u\in\Theta_0} f(\u) - f^\star\Big] + \Big(n + \frac{\sqrt{n}}{2}\Big) \nu \delta. 
\end{equation}
Both \eqref{eq: GFO convergence} and \eqref{eq: MN convergence} provide bounds on averaged stationarity measures. 
In \eqref{eq: GFO convergence}, stationarity is measured by the expected squared gradient norm, namely $\mathbb{E}\|Df(\xx_t)\|^2$. 
In \eqref{eq: MN convergence}, it is measured by the gradient norm $\|Df(\bc_t)\|$. 
On the right-hand sides, the first terms decrease at the rate $\mathcal{O}(1/k)$ and are proportional to the gap between an initial function value and $f^\star$. 
For \eqref{eq: MN convergence}, this initial value is represented by the average of the function values at the vertices of the initial simplex. 
The second terms determine the size of the neighborhood around zero to which these averaged stationarity measures can be driven.

However, there are also several important differences between \eqref{eq: GFO convergence} and \eqref{eq: MN convergence} as the mechanisms of the two algorithms are fundamentally different. 
Notice the expected length of $\u$ is approximately $\sqrt{n}$, so, when $n$ is large and $h_k = \frac{1}{4(n+4)\nu} \approx \frac{1}{4n\nu}$, we have 
\[
    h_k g_\mu(\xx_k) 
    \approx \frac{1}{4\nu} \cdot \frac{f(\xx_k + \mu \u) - f(\xx_k)}{\|\mu \u\|} \cdot \frac{\u}{\|\u\|}. 
\]
Therefore, the update rule in \eqref{eq: GFO} can be roughly seen as performing a gradient descent step restricted to the one-dimensional subspace $\operatorname{span}\{\u\}$, with an effective step size of order $\frac{1}{\nu}$, while using finite difference to approximate the projection of $\nabla f(\xx_k)$ onto this subspace.
From this perspective, we can see that $x_k$ is only updated in a random 1-dimensional subspace of $\R^n$ each time, therefore the $\mathcal{O}(1/k)$ term in \eqref{eq: GFO convergence} is similar to that of the regular gradient descent but has an extra $n$ factor. 

Meanwhile, the $\mathcal{O}(1/k)$ term in \eqref{eq: MN convergence} is proportional to $n^2$. 
Intuitively speaking, the first $n$ factor comes from the potential misalignment between the negative (approximated) gradient and the direction the simplex is moving $\bc_{k+1} - \bc_k$. 
The reflection operation changing only one of the $n+1$ vertices and moving the center only by $\|\bc_{k+1} - \bc_k\| = \frac{2}{n} \delta$ puts another $n$ in the numerator. 
It also puts a $\delta$ in the denominator because the movement of the simplex is limited by the radius $\delta$.

The sizes of the neighborhoods represented by the second terms on the right-hand sides of \eqref{eq: GFO convergence} and \eqref{eq: MN convergence} are affected by the dimension $n$, the Lipschitz constant $\nu$, and either $\mu$ or $\delta$. 
From the definition of $g_\mu$, we can interpret $\mu$ as the size of a finite difference step, as opposed to $\delta$, which is the simplex radius. 
These two parameters are conceptually different, but they both represent approximation errors in their corresponding algorithms that would prevent the gradient from converging to zero. 
However, by setting them to functions of $\varepsilon$, it is possible to obtain an $\varepsilon$-stationary point for any $\varepsilon > 0$. 
Indeed, another result from Section 7 of \cite{nesterov2017random} is that by choosing $\mu \le \mathcal{O}(\frac{\varepsilon}{n^{3/2} \nu})$, one can obtain an $\varepsilon$-stationary point (in expectation) in $\mathcal{O}(n/\varepsilon^2)$ iterations. 

In the randomized method \eqref{eq: GFO}, $\mu$ controls the finite-difference error, while the step size $h_k$ is chosen separately. 
In Algorithm~\ref{alg}, however, the same parameter $\delta$ controls both the accuracy of the local approximation and the size of the movement of the simplex. 
As a result, when $\delta$ is set to $\frac{2\varepsilon}{5n\nu}$ in Theorem~\ref{thm:NM2} so that the size of the neighborhood can be reduced to below $\varepsilon$, the $\mathcal{O}(1/k)$ term in \eqref{eq: MN convergence} becomes proportional to $n^3$, leading to the $\mathcal{O}(n^3 / \varepsilon^2)$ complexity. 
This means the performance of Algorithm~\ref{alg} is likely to be much poorer on high-dimensional problems compared to algorithm \eqref{eq: GFO}. 
This outcome aligns with empirical evidence: algorithm \eqref{eq: GFO} has been successfully employed in high-dimensional machine learning tasks (see, e.g., \cite{chen2017zoo,malladi2023fine}), whereas simplicial search methods generally struggle to scale effectively with dimensionality.

\section{Discussion} \label{sec:discussion}
We presented a QCQP formulation of the sharp bound $z(\Theta,\xx_0)$ on the function approximation error of linear interpolation. 
The mathematical expression of $z(\Theta,\xx_0)$ is then determined under several conditions. 
They include an upper bound \eqref{eq:phase1} that improves the existing ones to better cover the extrapolation case, a sharp bound \eqref{eq:phase2} when $f$ is quadratic, and a sharp bound \eqref{eq:phase3} for bivariate extrapolation. 
The two bounds \eqref{eq:phase2} and \eqref{eq:phase3} together provide the sharp error bound for bivariate linear interpolation under any configuration of $\xx_0$ and an affinely independent $\Theta$. 
These bounds can provide an important theoretical foundation for the design and analysis of derivative-free optimization methods and any other numerical methods that employs linear interpolation.

We proposed to compute $\{\mu_{ij}\}$ and check their signs to determine whether $\underline{z}(\Theta,\xx_0)$ from \eqref{eq:phase2} equals to $z(\Theta,\xx_0)$ and proved in Theorem~\ref{thm:phase2} that $\{\mu_{ij}\}$ being all non-negative is a sufficient condition for them to be equal.
We want to mention that one of our numerical experiments seems to indicate that it is also a necessary condition. 
This experiment involves generating many different $\Theta$ and $\xx$ with various $n$ and calculating the corresponding $\{\mu_{ij}\}$. 
From this experiment, we also observed some geometric patterns of the signs of $\{\mu_{ij}\}$, which we present in the following conjecture. 
\begin{conjecture}
    Assume $f \in C^{1,1}_\nu(\R^n)$. 
    Let $\xx_0$ be a fixed vector in $\R^n$, and $\Theta = \{\xx_1,\dots,\xx_{n+1}\}\subset \R^n$ be a fixed set of $n+1$ affinely independent vectors. 
    Let $\hat{f}$ be the affine function that interpolates $f$ on $\Theta$. 
    Let $\{\mu_{ij}\}_{i \in \cI_+,~ j \in \cI_-}$ be the set of parameters defined in Theorem~\ref{thm:phase2}. Then the following statements are true. 
    \begin{enumerate}
        \item The parameters $\{\mu_{ij}\}$ are all non-negative regardless of the position of $\xx_0$ if there is no obtuse angle at the vertices of the simplex $\conv(\Theta)$, i.e.,  
        \begin{equation} \label{eq:acute simplex} 
            (\xx_j-\xx_i) \cdot (\xx_k-\xx_i) \ge 0 \text{ for all } i,j,k = 1,2,\dots,n+1. 
        \end{equation}
        \item If there is at least one obtuse angle at the vertices of the simplex $\conv(\Theta)$, then there is a non-empty subset of $\R^n$ such that if $\xx_0$ is in this subset, then there is at least one negative element in $\{\mu_{ij}\}$. 
    \end{enumerate}
\end{conjecture} 

A general formula for the sharp bound on the function approximation error of linear interpolation and extrapolation remains an open question. 
It would appear $G\cdot H/2$ is a good candidate, since both sharp bounds \eqref{eq:phase2} and \eqref{eq:phase3} can be written in this form, but the matrix $H$ depends on the geometry of $\Theta$ and $\xx_0$. 
Using $G\cdot H/2$ as the general formula, we would need five different definitions of $H$ just for the bivariate case (\eqref{eq:Hstar} and four variants of \eqref{eq:Hstar 3} that corresponds to the four shaded areas in Figure~\ref{fig:phase2}). 
Note that the matrix $H$ is tied to $\{\mu_{ij}\}$ in \eqref{eq:mu1 +} and \eqref{eq:mu1 -}, and we believe even when there are negatives in $\{\mu_{ij}\}$, they are still tied in the same manner to a version of $\{\mu_{ij}\}$ that is modified to be all non-negative. 
In fact, \eqref{eq:mu0 +} - \eqref{eq:mu1 -} all hold true under the setting of Theorem~\ref{thm:phase3} if $H$ is defined as the matrix $\tilde{H}^\star$ from \eqref{eq:Hstar 3} and $\{\mu_{ij}\}$ is defined as
\[ \begin{aligned} 
\mu_{10} &= 1-\ell_2(\xx_0), 
&\mu_{13} &= -\ell_3(\xx_0), 
&\mu_{20} &= \ell_2(\xx_0), 
&\mu_{23} &= 0, 
\end{aligned} \]
which are the coefficients in \eqref{eq:phase3 sum}. 
Considering the difficulty in analyzing the signs of $\{\mu_{ij}\}$, it is unlikely for $G\cdot H/2$ to be suitable for this general formula. 
Whether there even exists a concise analytical form to the sharp error bound that can fit all the geometric configurations of $\Theta$ and $\xx_0$ is still unclear to us. 

Another topic worth investigating is the gradient approximation error of linear interpolation. 
This, however, might be much more difficult. 
If the error is measured by the Euclidean norm of the difference, then the problem of finding its sharp bound can be formulated in a way similar to \eqref{prob:D}: 
\begin{equation} \everymath{\displaystyle} \begin{array}{ll}
    \max_{f} &\|D\hat{f}(\xx_0) - Df(\xx_0)\| \\
    \text{s.t. } &f \in C_\nu^{1,1}(\R^n). 
\end{array}
\end{equation}
No matter how the constraint is handled, the objective of this problem is to maximize a convex function, making it a much less tractable nonconvex problem. 
The problem can become tractable if the error is measured differently, but whether the measure is meaningful depends on the application. 
For those who are interested in this error, some results can be found in \cite{handscornb1995errors, cao2005error, berahas2022theoretical}.

\section*{Acknowledgment}
We would like to acknowledge the help from Dr. Xin Shi and Yunze Sun in solving \eqref{prob:quadratic} and the help from Ruisong Zhou in proving Lemma~\ref{lem: NM f g}. 
We would like to thank Dr. Shuonan Wu and Dr. Katya Scheinberg for carefully reading this paper and providing their suggestions. 
We are grateful to the reviewers of Optimization Methods and Software for their thorough reviews and constructive comments. 

\bibliographystyle{plain}
\bibliography{references}

@book{DFO_book,
  title={Introduction to derivative-free optimization},
  author={Conn, Andrew R and Scheinberg, Katya and Vicente, Luis N},
  year={2009},
  publisher={SIAM}
}

@incollection{powell1994direct,
  title={A direct search optimization method that models the objective and constraint functions by linear interpolation},
  author={Powell, Michael JD},
  booktitle={Advances in Optimization and Numerical Analysis},
  pages={51--67},
  year={1994},
  publisher={Springer}
}

@article{nelder1965simplex,
  title={A simplex method for function minimization},
  author={Nelder, John A and Mead, Roger},
  journal={The computer journal},
  volume={7},
  number={4},
  pages={308--313},
  year={1965},
  publisher={Oxford University Press}
}

@book{Nesterov_book,
  title={Lectures on convex optimization},
  author={Nesterov, Yurii},
  volume={137},
  year={2018},
  publisher={Springer Optimization and Its Applications}
}

@book{davis1975book,
  title={Interpolation and Approximation},
  author={Davis, Philip J},
  year={1975},
  publisher={Courier Corporation}
}

@article{ciarlet1972general,
  title={General {L}agrange and {H}ermite interpolation in {$\R^n$} with applications to finite element methods},
  author={Ciarlet, Philippe G and Raviart, Pierre-Arnaud},
  journal={Archive for Rational Mechanics and Analysis},
  volume={46},
  number={3},
  pages={177--199},
  year={1972},
  publisher={Springer}
}

@article{powell2001lagrange,
  title={On the {L}agrange functions of quadratic models that are defined by interpolation},
  author={Powell, MJD},
  journal={Optimization Methods and Software},
  volume={16},
  number={1-4},
  pages={289--309},
  year={2001},
  publisher={Taylor \& Francis}
}

@techreport{handscornb1995errors,
  title={Errors of linear interpolation on a triangle},
  author={Handscornb, DC},
  institution = {Oxford University Computing Laboratory}, 
  year={1995}
}

@article{waldron1998error,
  title={The error in linear interpolation at the vertices of a simplex},
  author={Waldron, Shayne},
  journal={SIAM Journal on Numerical Analysis},
  volume={35},
  number={3},
  pages={1191--1200},
  year={1998},
  publisher={SIAM}
}

@article{stampfle2000optimal,
  title={Optimal estimates for the linear interpolation error on simplices},
  author={St{\"a}mpfle, Martin},
  journal={Journal of Approximation Theory},
  volume={103},
  number={1},
  pages={78--90},
  year={2000},
  publisher={Elsevier}
}

@article{cao2005error,
  title={On the error of linear interpolation and the orientation, aspect ratio, and internal angles of a triangle},
  author={Cao, Weiming},
  journal={SIAM Journal on Numerical Analysis},
  volume={43},
  number={1},
  pages={19--40},
  year={2005},
  publisher={SIAM}
}

@article{berahas2022theoretical,
  title={A theoretical and empirical comparison of gradient approximations in derivative-free optimization},
  author={Berahas, Albert S and Cao, Liyuan and Choromanski, Krzysztof and Scheinberg, Katya},
  journal={Foundations of Computational Mathematics},
  volume={22},
  number={2},
  pages={507--560},
  year={2022},
  publisher={Springer}
}

@article{drori2014performance,
  title={Performance of first-order methods for smooth convex minimization: {A} novel approach},
  author={Drori, Yoel and Teboulle, Marc},
  journal={Mathematical Programming},
  volume={145},
  number={1-2},
  pages={451--482},
  year={2014},
  publisher={Springer}
}

@article{taylor2017smooth,
  title={Smooth strongly convex interpolation and exact worst-case performance of first-order methods},
  author={Taylor, Adrien B and Hendrickx, Julien M and Glineur, Fran{\c{c}}ois},
  journal={Mathematical Programming},
  volume={161},
  pages={307--345},
  year={2017},
  publisher={Springer}
}

@article{taylor2017exact,
  title={Exact worst-case performance of first-order methods for composite convex optimization},
  author={Taylor, Adrien B and Hendrickx, Julien M and Glineur, Fran{\c{c}}ois},
  journal={SIAM Journal on Optimization},
  volume={27},
  number={3},
  pages={1283--1313},
  year={2017},
  publisher={SIAM}
}

@book{horn2012matrix,
  title={Matrix Analysis},
  author={Horn, Roger A and Johnson, Charles R},
  year={2012},
  publisher={Cambridge university press}
}

@article{sylvester1852xix,
  title={XIX. A demonstration of the theorem that every homogeneous quadratic polynomial is reducible by real orthogonal substitutions to the form of a sum of positive and negative squares},
  author={Sylvester, James Joseph},
  journal={The London, Edinburgh, and Dublin Philosophical Magazine and Journal of Science},
  volume={4},
  number={23},
  pages={138--142},
  year={1852},
  publisher={Taylor \& Francis}
}

@article{spendley1962sequential,
  title={Sequential application of simplex designs in optimisation and evolutionary operation},
  author={Spendley, WGRFR and Hext, George R and Himsworth, Francis R},
  journal={Technometrics},
  volume={4},
  number={4},
  pages={441--461},
  year={1962},
  publisher={Taylor \& Francis}
}

@article{conn1997convergence,
  title={On the convergence of derivative-free methods for unconstrained optimization},
  author={Conn, Andrew R and Scheinberg, Katya and Toint, Ph L},
  journal={Approximation theory and optimization: tributes to MJD Powell},
  pages={83--108},
  year={1997}
}

@article{conn2008geometry,
  title={Geometry of interpolation sets in derivative free optimization},
  author={Conn, Andrew R and Scheinberg, Katya and Vicente, Lu{\'\i}s N},
  journal={Mathematical Programming},
  volume={111},
  pages={141--172},
  year={2008},
  publisher={Springer}
}

@article{conn2009global,
  title={Global convergence of general derivative-free trust-region algorithms to first- and second-order critical points},
  author={Conn, Andrew R and Scheinberg, Katya and Vicente, Lu{\'\i}s N},
  journal={SIAM Journal on Optimization},
  volume={20},
  number={1},
  pages={387--415},
  year={2009},
  publisher={SIAM}
}

@article{bandeira2012computation,
  title={Computation of sparse low degree interpolating polynomials and their application to derivative-free optimization},
  author={Bandeira, Afonso S and Scheinberg, Katya and Vicente, Lu{\'\i}s Nunes},
  journal={Mathematical Programming},
  volume={134},
  number={1},
  pages={223--257},
  year={2012},
  publisher={Springer}
}

@article{bandeira2014convergence,
  title={Convergence of trust-region methods based on probabilistic models},
  author={Bandeira, Afonso S and Scheinberg, Katya and Vicente, Luis Nunes},
  journal={SIAM Journal on Optimization},
  volume={24},
  number={3},
  pages={1238--1264},
  year={2014},
  publisher={SIAM}
}

@article{gratton2018complexity,
  title={Complexity and global rates of trust-region methods based on probabilistic models},
  author={Gratton, Serge and Royer, Cl{\'e}ment W and Vicente, Lu{\'\i}s N and Zhang, Zaikun},
  journal={IMA Journal of Numerical Analysis},
  volume={38},
  number={3},
  pages={1579--1597},
  year={2018},
  publisher={Oxford University Press}
}

@article{chen2018stochastic,
  title={Stochastic optimization using a trust-region method and random models},
  author={Chen, Ruobing and Menickelly, Matt and Scheinberg, Katya},
  journal={Mathematical Programming},
  volume={169},
  pages={447--487},
  year={2018},
  publisher={Springer}
}

@article{blanchet2019convergence,
  title={Convergence rate analysis of a stochastic trust-region method via supermartingales},
  author={Blanchet, Jose and Cartis, Coralia and Menickelly, Matt and Scheinberg, Katya},
  journal={INFORMS Journal on Optimization},
  volume={1},
  number={2},
  pages={92--119},
  year={2019},
  publisher={INFORMS}
}

@article{cao2023first,
  title={First- and second-order high probability complexity bounds for trust-region methods with noisy oracles},
  author={Cao, Liyuan and Berahas, Albert S and Scheinberg, Katya},
  journal={Mathematical Programming},
  pages={1--52},
  year={2023},
  publisher={Springer}
}

@article{woods1985interactive,
  title={An interactive approach for solving multi-objective optimization problems},
  author={Woods, Daniel John},
  year={1985}
}

@article{mckinnon1998convergence,
  title={Convergence of the {N}elder-{M}ead Simplex method to a nonstationary Point},
  author={McKinnon, Ken IM},
  journal={SIAM Journal on Optimization},
  volume={9},
  number={1},
  pages={148--158},
  year={1998},
  publisher={SIAM}
}

@article{lagarias2012convergence,
  title={Convergence of the restricted {N}elder-{M}ead algorithm in two dimensions},
  author={Lagarias, Jeffrey C and Poonen, Bjorn and Wright, Margaret H},
  journal={SIAM Journal on Optimization},
  volume={22},
  number={2},
  pages={501--532},
  year={2012},
  publisher={SIAM}
}

@article{tseng1999fortified,
  title={Fortified-descent simplicial search method: A general approach},
  author={Tseng, Paul},
  journal={SIAM Journal on Optimization},
  volume={10},
  number={1},
  pages={269--288},
  year={1999},
  publisher={SIAM}
}

@article{larson2019derivative,
  title={Derivative-free optimization methods},
  author={Larson, Jeffrey and Menickelly, Matt and Wild, Stefan M},
  journal={Acta Numerica},
  volume={28},
  pages={287--404},
  year={2019},
  publisher={Cambridge University Press}
}

@book{coxeter1973regular,
  title={Regular polytopes},
  author={Coxeter, Harold Scott Macdonald},
  year={1973},
  publisher={Courier Corporation}
}

@article{nesterov2017random,
  title={Random gradient-free minimization of convex functions},
  author={Nesterov, Yurii and Spokoiny, Vladimir},
  journal={Foundations of Computational Mathematics},
  volume={17},
  number={2},
  pages={527--566},
  year={2017},
  publisher={Springer}
}

@inproceedings{chen2017zoo,
  title={Zoo: Zeroth order optimization based black-box attacks to deep neural networks without training substitute models},
  author={Chen, Pin-Yu and Zhang, Huan and Sharma, Yash and Yi, Jinfeng and Hsieh, Cho-Jui},
  booktitle={Proceedings of the 10th ACM workshop on artificial intelligence and security},
  pages={15--26},
  year={2017}
}

@article{malladi2023fine,
  title={Fine-tuning language models with just forward passes},
  author={Malladi, Sadhika and Gao, Tianyu and Nichani, Eshaan and Damian, Alex and Lee, Jason D and Chen, Danqi and Arora, Sanjeev},
  journal={Advances in Neural Information Processing Systems},
  volume={36},
  pages={53038--53075},
  year={2023}
}
\end{document}